\documentclass[11pt]{article}

\usepackage{graphicx}
\usepackage{epstopdf}
\usepackage{subfigure}
\usepackage{amsfonts,amssymb,amsmath}
\usepackage{makeidx}

\usepackage{url}
\usepackage{helvet}
\usepackage{multicol}

\evensidemargin \oddsidemargin

\makeindex

\newcommand{\uQ}{{\mathbf{Q}}}

\newcommand{\uP}{{\mathbf{P}}}
\newcommand{\uG}{{\mathbf{G}}}

\newcommand{\uK}{{\mathbf{K}}}
\newcommand{\uC}{{\mathbf{C}}}
\newcommand{\uTau}{{\mathbf{T}}}

\newtheorem{theorem}{\quad Theorem}

\newtheorem{lemma}{\quad Lemma}

\graphicspath{{.//Art/}}

\usepackage{graphicx}
\newcommand{\BibTeX}{{\rm B\kern-.05em{\sc i\kern-.025em b}\kern-.08em\ignorespaces
  T\kern-.1667em\lower.7ex\hbox{E}\kern-.125emX}}

\begin{document}

%%%%%%%%%%%%%%%%%%%%%%%%%%%%%%%%%%%%%%%%%%%%%%%%%%%%%%%%%%%%%%%%%%
%                         THE TOP MATTER                         %
%%%%%%%%%%%%%%%%%%%%%%%%%%%%%%%%%%%%%%%%%%%%%%%%%%%%%%%%%%%%%%%%%%
%
\title{Pollution control for a spatially structured \\ economic growth system}%

\author{
Sebastian Ani\c{t}a
\thanks{Faculty of Mathematics, ``Alexandru Ioan Cuza'' University of Ia\c{s}i, Ia\c{s}i 700506, Romania; 
``Octav Mayer'' Institute of Mathematics of the Romanian Academy,  Ia\c{s}i 700506, Romania.}
\and Vincenzo Capasso
\thanks{ADAMSS (Advanced Applied  Mathematical and Statistical Sciences), Universit\`a  degli Studi di Milano, Italy.}
\and Simone Scacchi
\thanks{Dipartimento di Matematica, Universit\`a degli Studi di Milano, Italy.}
}

\maketitle

\begin{abstract}
In this paper investigations by the same authors  on environmental
issues concerning the control of the pollution produced by   human
activities have been extended to include costs related to environmental interventions. 
  The proposed model consists of  a spatially
structured dynamic economic growth model which takes into account
the level of pollution induced by production,  a possible
taxation based on the amount of produced pollution, and possible environmental interventions.
 It has been 
analyzed an optimal harvesting control problem with an objective
function composed of four  terms, namely the intertemporal utility
of the decision maker, the space-time average of the level of
pollution in the habitat,  the disutility due to the imposition
of taxation and the cost of environmental interventions.

A specific novelty in the model proposed   here is the localization of the possible interventions  to a subregion of 
the whole habitat.

Computational experiments have been carried out to exemplify the outcomes of the proposed model.
\end{abstract}


%%%%%%%%%%%%%%%%%%%%%%%%%%%%%%%%%%%%%%%%%%%%%%%%%%%%%%%%%%%%%%%%%%
% \section*{Introduction}                                        %
% \section{About the "head" of your paper}                       %
% \subsection{Your private macros (the preamble)}                %
% \subsection{The top matter}                                    %
% \section{About the "body" of your paper}                       %
% \subsection{The environments}                                  %
% \subsection{specific macros}                                   %
% \subsection{Cross reference and bibliography}                  %
% \subsection{Including postscript files}                        %
%%%%%%%%%%%%%%%%%%%%%%%%%%%%%%%%%%%%%%%%%%%%%%%%%%%%%%%%%%%%%%%%%%
%                              THE TEXT                          %
%%%%%%%%%%%%%%%%%%%%%%%%%%%%%%%%%%%%%%%%%%%%%%%%%%%%%%%%%%%%%%%%%%
\tableofcontents
\section{Introduction}\label{intro}

As  a  possible  application of the methods presented in \cite{SA_VK_SS_2024}, in this paper we introduce an environmental issue
concerning the control of the pollution produced by   human
activities.

Standard macroeconomics and environmental economics
\index{environmental economics} have been two completely
independent research areas for a long time. The seminal paper by
Grossman and Kruger \cite{grossman95} has stimulated considerable
interest in the relationship between economic development and
environmental pollution.  The role of taxation \index{taxation}
for controlling pollution  is the key issue of the report
\cite{IMF13} (see also references therein)  issued by the
International Monetary Fund.   In view of these policy
developments,   some recent works tend to develop a global theory
combining the  above mentioned two branches of literature (see
\cite{BrTa10} and \cite{bertinelli12}).

Here  we   present  an extension of our  own results concerning the impact of
pollution diffusion on the economic and environmental joint
dynamics, and  consider the optimal control problem  of the relevant  system subject to capital consumption, to taxation based on the amount of produced pollution, and possible environmental interventions.
A key issue of this paper is the localization of the possible interventions  to a subregion of 
the whole habitat.

In \cite{ACKL13} and \cite{ACKL15} the authors introduced a
spatial economic growth model \index{economic growth!spatial
	model} subject to pollution diffusion\index{pollution!diffusion}.
In particular, in \cite{ACKL13} the authors analyzed the
large-time behavior of a spatially structured economic growth
model coupling physical capital accumulation and pollution
diffusion: this model extends other results in the literature
along different directions. Furthermore, a negative feedback to
the production function was added, in order to describe the
(negative) influence of pollution on capital accumulation. In
\cite{ACKL15} the authors modified the model presented in
\cite{ACKL13} by introducing an optimal harvesting control problem
by  an objective function composed of three terms, namely the
intertemporal utility of the decision maker, the space-time
average of the level of pollution in the habitat, and the
disutility due to the imposition of taxation. In \cite{AKKT2017}
two controls were  introduced, namely the level of consumption
and, in addition, the level of taxation on physical capital. This
includes taxation on machinery, buildings, computers, and any
equipment that help to turn the raw materials into finished
products or services and because of these activities they have a
huge impact on pollution emissions.

In the literature other
contributions which explore the spatial dimension in environmental
and resource economics can be found in \cite{BrXe08},
\cite{BrXe10},  \cite{Xe10}, and \cite{AtXe12}. Of interest it might be  also \cite{BV},  concerning 
the optimal control of age dependent systems.

%  spatial  models

The inclusion of a spatial dimension in economic analysis has
recently increased relevance and interest. The first studies in
economic geography go back to Beckerman \cite{Beck} and Puu
\cite{Pu82}, who study regional problems based simply on flow
equations. Starting from these works, a new economic geography
arises adopting general equilibrium models to analyze the
peculiarities of local and global markets, and the mobility of
production factors (see \cite{Kr91}, \cite{FukrVe99}). More
recently, this geographical approach has been introduced in
economic growth models \index{geographical economics} to study the
implications between accumulation and diffusion of capital on
economic dynamics (see \cite{Br04},  \cite{CaZo04},
\cite{BoCaZo09}). The Solow model \cite{So56} with a continuous
spatial dimension has been extensively studied.   In \cite{CaZo04}
the authors analyze problems of convergence across regions when
capital is mobile, while \cite{Br04} considers the case in which
both capital and labor are mobile. The Ramsey model \cite{Ramsey}
has been extended to a spatial dimension in \cite{Br04} and
\cite{BoCaZo09}, respectively in average and total utilitarianism
version. In \cite{CaEnLa10}, alongside the classical Cobb--Douglas
production function, \index{Cobb--Douglas!production function}
\index{production function!Cobb--Douglas} the authors have
proposed a larger class of convex--concave production functions,
\index{production function!convex-concave}thus allowing a larger
variety of dynamical scenarios, including the existence of poverty
traps. For  a recent review  see \cite{ABV_2019}.

In addition to the
impact of taxation on the  production  as emerged  in
\cite{AKKT2017}, in this paper  we consider the role of additional costs due to possible environmental interventions.  

For reaction models with nonlocal reaction terms in mathematical biology the reader may  refer to \cite{SA_VK_SS_2024},  \cite{V1} and  \cite{V2}. 

In Section \ref{sec:2} the  optimal control problem
is presented, by  acting on both consumption,   taxation, and environmental interventions. Since
the  search  of the optimal control is not an easy analytical
task, we have proposed   a computational algorithm for the search
by the  gradient method. To this aim, in the same Section
\ref{sec:2} the evaluation of the directional derivative of the
cost functional with respect to the three  controls has been
presented. In Section \ref{numerics}  a numerical
algorithm that allows to solve the system    (\ref{dynamics_regional_pollution_treatment_new}), (\ref{kopt}), (\ref{popt}) of four forward-backward
reaction-diffusion equations, subject to the relevant initial/final and boundary conditions,  has been reported, together  with
a set of computational experiments (see  Section \ref{sec:3}).  A discussion on the outcomes of the computational experiments is presented in Section \ref{concluding_remarks}.

\section{A spatially structured model}\label{sec:1}

By following our previous literature on the subject   (see \cite{AKKT2017}, and references therein), we shall consider the following spatially structured model. As far as the evolution of the stock capital is concerned, let  $k(x,t)$ denote the capital stock held by the representative
household located at $x \in \Omega $  (a nonempty bounded domain
in $\mathbb R^n$, $n=1,2$ with a sufficiently smooth boundary), at
time $t \geq 0$. The time evolution of the physical capital
\index{capital} is modelled by

\begin{equation} \label{eq:solow}
{{\partial k}\over {\partial t}}(x,t) = d_1 \Delta k(x,t) + A(x,t) f
\big(k(x,t)\big) - \delta_1 k(x,t) - c(x, t) k(x,t),
\end{equation}
for $(x,t)\in \Omega\times (0,+\infty).$

Here  $A(x,t)$ denotes the technological progress
\index{technological progress} available at $x \in \Omega,$ and
time $t \geq 0$. The parameter $\delta_1$ takes into account  the
depreciation of the capital. $c(x, t)$ describes the level of
consumption per unit of physical capital $k(x,t);$  the quantity $
(1-c(x, t))k(x,t)$  describes the saved capital after
consumption.

Equation (\ref{eq:solow})  is complemented by an initial  capital
distribution, $k(x,0)$, and  suitable boundary conditions. If we
assume that there is no capital flow through the boundary
$\partial\Omega$, we impose homogeneous  Neumann boundary
conditions, i.e.
\begin{equation*}
{{\partial k}\over {\partial n}}(x, t)=0, \quad x \in
\partial\Omega, \ t>0.
\end{equation*}

The production function $f$ can be taken of the following form
\index{production function} (see \cite{CaEnLa10}, and references
therein)
\begin{equation*}
f(r)={{\alpha _1r^{\gamma}}\over{1+\alpha _2r^{\gamma }}} ,
%\label{production}
\end{equation*}
where $\alpha _1\in (0, +\infty ), \ \alpha _2\in [0,+\infty ), \
\gamma \in [1,+\infty )$.

As far as  pollution diffusion \index{pollution!production} is
concerned, if we denote by $p(x,t)$ the pollution stock faced by a
representative household located at $x \in \Omega,$  at time $t
\ge 0,$  we assume the following model

\begin{equation}
{{\partial p}\over {\partial t}}(x,t)=d_2\Delta p(x,t)+\theta \int_{\Omega }
f(k(x',t))\varphi(x',x)dx' -\delta_2 p(x,t).
\label{pollution_dynamics}
\end{equation}

The first term, including the Laplace operator,  may model the random dispersal of the pollutant in the relevant habitat. Via the integral term  in Equation (\ref{pollution_dynamics}),  we
may see the impact  of the economy on pollution production.  At
any point $x \in \Omega$  pollution  is   generated  due to the
production $f(k(x,t)).$ In our  model we have supposed that any
point $x' \in \Omega$ may  have an impact on  enhancing pollution
at any other point $x \in \Omega,$  by  any  other possible
mechanism than random dispersal already described by the Laplace
operator. Such mechanisms are  described by a kernel $\varphi $
which satisfies the following hypotheses: $\varphi \in L^{\infty
}(\Omega \times \Omega )$, and $\varphi (x', x)\geq 0$ a.e. $(x',
x) \in \Omega \times \Omega. $  As a limiting case  one may reduce
the kernel $\varphi$  to a Dirac Delta function, to mean the
impact  of production on pollution at the same point $x.$  The
parameter $\theta$ takes into account the ``efficiency'' of the
economic system with respect to pollution production.  The
parameter $\delta_2$ accounts for the natural decay  of pollution.

Equation  (\ref{pollution_dynamics}) has to be  complemented by
suitable initial  and boundary conditions.

Actually we have  required that the model takes into account the negative feedback due to the impact of pollution on production, since  it is expected that the relevant financial agents tend to reduce taxation and other indirect costs due   to pollution, by a self control; as a consequence  in \cite{AKKT2017} we have modified Equation (\ref{eq:solow}) as follows

%\begin{eqnarray}
%\partial _tk(x,t)&=\displaystyle d_1\Delta k(x,t)+\frac{A}{g(p(x, t))}
%f(k(x,t))  \nonumber  \\
% \ & \ \ \ \ \ -\delta_1 k(x,t)-c(x,t)k(x,t).
%  \label{production_dynamics}
%\end{eqnarray}

\begin{equation}  \label{production_dynamics}
{{\partial k}\over {\partial t}}(x,t)=\displaystyle d_1\Delta k(x,t)+\frac{A(x,t)}{g(p(x, t))}
f(k(x,t)) -\delta_1 k(x,t)-c(x,t)k(x,t).
\end{equation}

All  parameters, $ A, d_1, d_2, \delta_1, \delta_2,$  which
characterize a specific economic/\newline ecological system are
assumed strictly positive. \index{economic/ecological system}

The function $g: [0, + \infty) \rightarrow (0, + \infty)$ models
the negative impact of pollution on the production;  usual
assumptions on $g$ are

\begin{itemize}
	\item[i)] $g$  is continuously differentiable and monotonically
	increasing;
	
	\item[ii)]$g(0) =1$ and $ \displaystyle \lim_{p \rightarrow
		+\infty} g(p) =  +\infty.$
\end{itemize}

As a typical choice for the function $g,$   later we shall adopt
the expression

$$ g(p) = 1+ \chi p^2, \quad p \ge 0.$$

In the sequel  it will  be kept $\chi=1,$  as a   typographical simplification.

\section{An optimal control problem}\label{sec:2}  \index{economic/ecological system!optimal control}

We plan to study our eco-economic system on a
finite time horizon $[0,T]$. We  shall  denote by

$$  Q_T =  \Omega \times (0, T); \quad \Sigma_T = \partial \Omega \times (0, T).   $$

In order  to control pollution emissions by production processes,  
local authorities may impose taxes proportional to the level of
production. A possible approach might be to impose taxes related to interventions on pollution abatement policies implemented only in the areas of interest, which we have indicate by a subregion $\omega \subset \Omega$ of the whole habitat. We shall denote by $\tau(x,t)$ the environmental taxation rate. 
By taking all this into account, we need to modify System (\ref{production_dynamics})-(\ref{pollution_dynamics}) as follows (see \cite{BrTa10}, \cite{AKKT2017}):

\begin{equation}   \left\{ \begin{array}{ll}
\displaystyle {{\partial k}\over {\partial t}}(x,t)&=\displaystyle
d_1\Delta k(x,t)+\frac{A (x, t)}{1+ p(x,t)^2}
f((1- \mathbb{I}_{\omega}(x)\tau(x,t))k(x,t)) \\
~ & ~ \\
\ & \ \ \ \ \ -\delta_1 k(x,t)-c(x,t)k(x,t),  \\
~ & ~ \\
\displaystyle {{\partial p}\over {\partial t}}(x,t)&=\displaystyle
d_2\Delta p(x,t)+\theta \int_{\Omega }
f((1- \mathbb{I}_{\omega}(x')\tau(x',t))k(x',t))\varphi(x',x)dx' \\
~ & ~ \\
\ & \ \ \ \ \ -\delta_2 p(x,t)  -  \mathbb{I}_{\omega}(x) \xi (x,t) p(x,t),
\end{array}
\right. \label{dynamics_regional_pollution_treatment_new}
\end{equation}
for $(x,t)\in Q_T,$  where  $\omega$  is a nontrivial  open subset
of  $\Omega.$  System (\ref{dynamics_regional_pollution_treatment_new}) shall be complemented by initial and boundary conditions to be precised later.

The term $(1-\tau(x,t))k(x,t)$ accounts for the amount of the
physical capital left to be used in the gross domestic production.
An increase in the level of taxation allows the implementation of practices leading to  an abatement of the level of pollution
emissions, so that  we may  interpret  it as a ``green'' taxation
policy. \index{green taxation policy}
By $\xi(x,t)$ we have denoted the rate of pollution abatement by suitable environmental interventions.

In this case  the optimal control  problem would be

$$\underset{\substack{    (c,\tau, \xi )\in \widetilde{U},\\  \omega \subset \Omega}} {Maximize} \ \left \{
\int_0^T\int_{\Omega} c(x,t)k(x,t)dx \ dt -  \beta_0
\int_0^T\int_{\Omega} p(x,t) dx \ dt \right . $$
$$- \beta_1\int_0^T\int_{\omega } \tau(x,t) k(x,t) dx \ dt   $$
$$ \left . - \beta_2\int_0^T\int_{\omega } \xi (x,t) p(x,t) dx \ dt \right \} ,$$
subject to  System (\ref{dynamics_regional_pollution_treatment_new}) and the initial and boundary conditions given below.

\par
It is   clear  that the solutions of the above problems would
depend upon the geographical distribution  of the potentially
polluting production activities; in this respect we may expect a
crucial role played by the functional parameter $A(x, t).$ This
fact would confirm the relevance of spatially structured models
for geographical economics.

For simplicity we may
consider homogeneous Neumann boundary conditions
\begin{equation}
{{\partial k}\over {\partial n}}(x,t)={{\partial p}\over {\partial n}}(x,t)=0, \qquad (x,t)\in \Sigma_T,
\label{obc}
\end{equation}
and  initial conditions
\begin{equation}
k(x,0)=k_0(x), \ p(x,0)=p_0(x), \qquad x\in \Omega , \label{oic}
\end{equation}
where $T>0$ is fixed and $k_0, p_0\in L^{\infty }(\Omega )$, $k_0(x)\geq 0, \ p_0(x)\geq 0$
a.e. $x\in \Omega $.

We   impose the constraints

\begin{equation}
1-c(x, t) - \tau(x, t)\ge s, \ \  c(x, t), \tau(x, t)\ge 0, \quad \mbox{\rm a.e. } (x,t)\in Q_T,
\label{constraints}
\end{equation}
where  $s>0$  is the  saving factor. In addition, we may assume   that  $0 \le \xi (x,t) \le
L,$ a.e. $(x,t)\in \omega \times (0,T)$, for a given constant $L.$

If we denote by $U=\{(c,\tau )\in L^{\infty}(Q_{T})\times L^{\infty}(Q_{T}): \
(c,\tau ) $ satisfies $(\ref{constraints})\},$  the set of
controls is   $\widetilde{U}:= U \times  \{ \xi \in L^{\infty} (\omega
\times (0, T)): \ 0 \le \xi (x,t) \le L, \ \mbox{\rm a.e. } (x,t)\in \omega \times (0,T) \} $  and  $\beta_0, \beta_1, \beta_2$  are
positive constants;  $(k,p)$ is the solution to
(\ref{dynamics_regional_pollution_treatment_new})-(\ref{oic}) corresponding to $Q_T$ and
$\Sigma_T$. 
The objective function is composed of  three terms that we hereby
describe:

\begin{itemize}
	\item[$T1:$] $\, \,$ the term $\displaystyle \int_0^T\int_{\Omega}
	c(x,t)k(x,t)dx \, dt$, to be maximized, describes the
	intertemporal utility of the decision maker,
	
	\item[$T2:$] $\, \,$ the second term $\displaystyle \beta_0
	\int_0^T\int_{\Omega} p(x,t) dx \, dt$, to be minimized,
	represents an average of the total level of pollution across space
	and time,
	
	\item[$T3:$] $\, \,$ the third term $\beta_1 \displaystyle
	\int_0^T\int_{\Omega} \tau(x,t) k(x,t) dx \, dt$, to be minimized,
	takes into account   the disutility due to imposition of taxation.
	
	\item[$T4:$] $\, \,$ the fourth term $ \beta_2 \displaystyle
	\int_0^T\int_{\Omega} \xi(x,t) p(x,t) dx \, dt$, to be minimized,
	takes into account the costs due to environmental interventions.
\end{itemize}

The costate equations for the previous optimal control problem
are:

%\large
\begin{equation}\label{kopt}
\begin{array}{ll}
\displaystyle
\frac{\partial \lambda_{k}}{\partial t}(x,t)&= \displaystyle
-d_1 \Delta \lambda_{k}(x,t) -c(x,t)+\beta_1 \mathbb{I}_{\omega }(x)\tau(x,t) \vspace{2mm}\\
~ & \displaystyle
-\lambda_{k}(x,t)\left[{\dfrac{A(x,t)(1-\mathbb{I}_{\omega }(x)\tau(x,t))}{1+ p(x,t)^2}} f'((1-\mathbb{I}_{\omega }(x)\tau(x,t))k(x,t)) - \delta_1 -c(x,t)\right] \vspace{2mm}\\
~ &\displaystyle -(1-\mathbb{I}_{\omega }(x)\tau(x,t)) f'((1-\mathbb{I}_{\omega }(x)\tau(x,t))k(x,t))
\theta \int_\Omega  \lambda_{p} (x',t)\varphi(x,x')dx',
\end{array}
\end{equation}

\begin{equation}\label{popt}
\begin{array}{ll}
\displaystyle
\frac{\partial \lambda_{p}}{\partial t}(x,t)&= \displaystyle -d_2  \Delta \lambda_{p}(x,t)  +\beta_0+\beta _2\mathbb{I}_{\omega }(x)\xi (x,t) \vspace{2mm}\\
~ & \displaystyle +\lambda_{k}(x,t)\frac{2A(x,t)p(x,t)}{(1+p(x,t)^2)^2}f((1-\mathbb{I}_{\omega }(x)\tau(x,t))k(x,t))  +\delta_2\lambda_{p}(x,t)\vspace{2mm} \\
~ & \displaystyle +\mathbb{I}_{\omega }(x)\xi (x,t)\lambda _p(x,t)
\end{array}
\end{equation}

\noindent
together with homogeneous Neumann boundary conditions and zero
final data.

As in \cite{ROM_2015} it follows the well posedness of the initial/boundary problem for System (\ref{dynamics_regional_pollution_treatment_new}). Moreover it
follows in the same manner that there exists at least one optimal control
$(c^*,\tau^*, \xi^*)$.

Let  us now proceed  with the construction of a computational algorithm for the search of an optimal control by the gradient method.

For the sake of clarity we shall use the notation $k^u$, $p^u$, $\lambda_k^u$, $\lambda _p^u$
for $k, p, \lambda _k, \lambda _p$, where $u=(x,\tau ,\xi )\in \tilde{U}$ when some danger of confusion occurs. Denote by
$$I(u)= \int_0^T\int_{\Omega} c(x,t)k^{u}(x,t)dx \, dt-\beta_0 \int_0^T\int_{\Omega} p^{u}(x,t) dx \, dt  $$
$$- \beta_1\int_0^T\int_{\Omega} \mathbb{I}_{\omega }(x)\tau(x,t)k^{u}(x,t) dx \, dt-\beta_2\int_0^T\int_{\Omega} \mathbb{I}_{\omega }(x)\xi (x,t)p^{u}(x,t) dx \, dt$$
the value of the cost functional corresponding to $u=(c,\tau ,\xi)$.

Consider an arbitrary $v=(v_1,v_2,v_3)\in L^{\infty }(Q_{T})\times L^{\infty }(\tilde{Q}_{T})\times L^{\infty }(\tilde{Q}_{T})$ (where $\tilde{Q}_T=\omega \times (0,T)$) such that $u+\varepsilon v\in \tilde{U}$ for any $\varepsilon >0$ sufficiently small. Let $(z_k,z_p)$ be the solution to
\begin{equation}\label{eq-var}
\left\{ \begin{array}{ll}
\displaystyle
\frac{\partial z_k}{\partial t}=&\displaystyle d_1\Delta z_k+\frac{A}{1+ (p^{u})^2}\bigg[f'((1-\mathbb{I}_{\omega }\tau )k^{u})(1-\mathbb{I}_{\omega }\tau )z_k \\
~ & ~ \\
\ & \ \ \ \ \ \displaystyle -f'((1-\mathbb{I}_{\omega }\tau )k^{u})\mathbb{I}_{\omega }k^{u}v_2 \bigg]  -\delta_1 z_k-c z_k-v_1k^{u} \\
\ & \ \\
\ & \ \ \ \ \ \displaystyle -{{2Ap^{u}z_p}\over {(1+(p^{u})^2)^2}}f((1-\mathbb{I}_{\omega }\tau )k^{u}), \\
\ & \ \\
\displaystyle
\frac{\partial z_p}{\partial t}=& \displaystyle d_2\Delta z_p+\theta \int_{\Omega }f'((1-\mathbb{I}_{\omega }(x')\tau(x',t))k^{u}(x',t)) \\
\ & \ \\
\ & \ \ \ \ \ \displaystyle \times [(1-\mathbb{I}_{\omega }(x')\tau(x',t))z_k(x',t)-\mathbb{I}_{\omega }(x')v_2(x',t)k^{u}(x',t)]\varphi(x',x)dx' \\
~ & ~ \\
~ &   \ \ \ \ \ -\delta_2 z_p-\mathbb{I}_{\omega }\xi z_p-\mathbb{I}_{\omega }v_3p^u,
\end{array}
\right.
\end{equation}
for $(x,t)\in Q_{T}$ together with homogeneous Neumann boundary
conditions and zero initial data. One may prove in a standard way
the following auxiliary result \vspace{3mm}

\noindent
\begin{lemma}	
The following convergences hold for any $u\in \tilde{U}$
\begin{itemize}
	\item[(i)] $k^{u+\varepsilon v}\longrightarrow
	k^{u}, \quad p^{u+\varepsilon v}\longrightarrow p^{u}$ in $L^2(Q_{T})$, as
	$\varepsilon \rightarrow 0+$;  \vspace{2mm}
	\item[(ii)] $\displaystyle
	\frac{k^{u+\varepsilon v}-k^{u}}{\varepsilon } \longrightarrow z_k$, \ \ \ $\displaystyle
	\frac{p^{u+\varepsilon v}-p^{u}}{\varepsilon } \longrightarrow z_p$ in $L^2(Q_{T})$, as
	$\varepsilon \rightarrow 0+$.
\end{itemize}
\end{lemma}	

One may compute now the directional derivative of $I$:

\vspace{3mm}

\noindent
\begin{theorem} 
$$dI(u)(v)=\int_{Q_T}k^{u}(x,t)[v_1(x,t)(1-\lambda _k(x,t))-\mathbb{I}_{\omega }(x)v_2(x,t)\beta^*(x,t)]dx \, dt,$$
$$-\int_{Q_T}p^{u}(x,t)v_3(x,t)\mathbb{I}_{\omega }(x)(\beta _2+\lambda _p(x,t))dx \, dt,$$
where
$$
\beta^*(x,t)= \beta_1 $$
$$+ f'((1-\mathbb{I}_{\omega }(x)\tau (x,t))k^{u}(x,t))
\left(\displaystyle{\frac{A(x,t)\lambda_k (x, t)}{1+p^{u}(x,t)^2}}  +
\theta \int_\Omega \lambda_p(x',t)\varphi(x,x') dx'\right).$$
\end{theorem}

\it Sketch of the proof. \rm For any $\varepsilon >0$ sufficiently
small we have that
$$ I(u+\varepsilon v)-I(u)=\int_{Q_T}[(c+\varepsilon v_1)k^{u+\varepsilon v}-ck^{u}]dx \, dt$$
$$-\int_{Q_T}\{ \beta _0[p^{u+\varepsilon v}-p^{u}]
+\beta_1\mathbb{I}_{\omega }[(\tau +\varepsilon v_2)k^{u+\varepsilon v}-\tau k^{u} ]\} dx \, dt $$ 
$$-\int_{Q_T}\{ \beta _2\mathbb{I}_{\omega }[(\xi +\varepsilon v_3)p^{u+\varepsilon v}-\tau p^{u} ]\} dx \, dt .$$ 
By dividing the above by  $\varepsilon $ and passing to the limit
($\varepsilon \rightarrow 0+$), then using $(i), (ii)$ in the
lemma we get that

\begin{equation}\label{patru}
\begin{array}{ll}
dI(u)(v)&=\displaystyle \int_{Q_T}[v_1k^{u}+cz_k]dx \ dt -\beta_0\int_{Q_T }z_p \,dx \, dt \vspace{2mm} \\
~ & \displaystyle \displaystyle -\beta_1\int_{Q_T }\mathbb{I}_{\omega }[v_2k^{u}+\tau z_k]dx \, dt -\beta_2\int_{Q_T }\mathbb{I}_{\omega }[v_3p^{u}+\xi z_p]dx \, dt .
\end{array}
\end{equation}

If we multiply  Equation (\ref{kopt}) by $z_k$ and integrate over
$Q_{T}$ and use (\ref{eq-var}), then we get that
\begin{eqnarray}
&&-\int_{Q_T }\bigg\{ -d_1\nabla \lambda_k (x,t)\cdot
\nabla z_k (x,t) \nonumber\\
\nonumber\\
&& +\frac{A(x,t)\lambda_k (x,t)}{1+p^2(x,t)}[f'((1-\mathbb{I}_{\omega }(x)\tau (x,t) )k
(x,t))(1-\mathbb{I}_{\omega }(x)\tau (x,t) )z_k (x,t)\nonumber\\
\nonumber\\
\nonumber &&-f'\!((1-\mathbb{I}_{\omega }(x)\tau (x,t))k(x,t))k(x,t)\mathbb{I}_{\omega }(x)v_2(x,t)]\!
\nonumber\\
\nonumber\\
\nonumber &&- \lambda _k (x,t)(\delta _1z_k (x,t)\!+\!cz_k (x,t)\!+\!v_1(x,t)k(x,t))\nonumber\\
\nonumber\\
&&-\dfrac{2A(x,t)p(x,t) z_p(x,t) \lambda_k (x,t)} {(1+p^2
	(x,t))^2}f((1-\mathbb{I}_{\omega }(x)\tau
(x,t))k(x,t))\bigg\} dx \, dt\nonumber\\
\nonumber\\
&&=\int_0^T\int_{\Omega }\bigg\{ d_1\nabla \lambda_k (x,t)\cdot
\nabla z_k (x,t) - cz_k (x,t)+   \beta_1\mathbb{I}_{\omega }(x)\tau(x,t)
z_k(x,t)\nonumber\\
\nonumber\\
&&\!-\!\lambda_k (x,t)z_k(x,t)\bigg[ \frac{A(x,t)(1\!-\mathbb{I}_{\omega }(x)\tau (x,t))}
{1+p^2(x,t)}f'((1\!-\mathbb{I}_{\omega }(x)\tau (x,t))k(x,t))\!-\!\delta _1\!-\!c
\bigg] \nonumber\\
\nonumber\\
&&-z_k(x,t)(1\!-\mathbb{I}_{\omega }(x)\tau (x,t))\nonumber\\
\nonumber\\
&&\times f'((1-\tau (x,t))k(x,t))\theta \int_{\Omega }\lambda
_p(x',t)\varphi (x,x')dx' \bigg\}\ dx \, dt.\nonumber
\end{eqnarray}

It follows that
\begin{eqnarray}
&&\displaystyle \int_{Q_T}(cz_k(x,t)-\beta_1\mathbb{I}_{\omega }(x)\tau (x,t)z_k(x,t))dx \, dt \nonumber \vspace {2mm}\\
&&\displaystyle =-\int_{Q_T }\int_{\Omega }z_k(x,t)(1\!-\mathbb{I}_{\omega }(x)\tau (x,t))f'((1\!-\mathbb{I}_{\omega }(x)\tau (x,t))k(x,t))\theta \lambda _p(x',t)\varphi (x,x')dx' \, dx \, dt \nonumber \vspace {2mm}\\
&&\displaystyle -\int_{Q_T }\!\! \bigg[
\frac{A(x,t)\lambda _k(x,t)}{(1+p^2(x,t))^2}f'((1-\mathbb{I}_{\omega }(x)\tau
(x,t))k(x,t))k(x,t)\mathbb{I}_{\omega }(x)v_2(x,t)\!\nonumber \vspace {2mm}\\
\nonumber \vspace {30mm}\\
&&+ \lambda _k(x,t)v_1(x,t)k(x,t)\nonumber \vspace {2mm}\\
\nonumber \vspace {30mm}\\
&&\displaystyle +\frac{2A(x,t)p(x,t)z_p(x,t)\lambda_k(x,t)}{(1+p^2(x,t))^2}f((1-\mathbb{I}_{\omega }(x)\tau(x,t) )k(x,t))\bigg] dx \,
dt.\label{cinci}
\end{eqnarray}

If we multiply  Equation (\ref{popt}) by $z_p$ and integrate over
$Q_T$ and use (\ref{eq-var}), then we get that
\begin{eqnarray}
&&-\int_{Q_T }\bigg\{ -d_2\nabla \lambda _p(x,t)\cdot
\nabla z_p(x,t)+\lambda _p(x,t)\theta \int_{\Omega }f'((1-\mathbb{I}_{\omega }(x')\tau
(x',t))k(x',t))\nonumber\vspace{2mm}\\
&&\times [(1-\mathbb{I}_{\omega }(x')\tau (x',t))z_k(x',t)-k(x',t)\mathbb{I}_{\omega }(x')v_2(x',t)]\varphi
(x',x)dx'-\delta _2 \lambda _p(x, t) z_p(x,t) \nonumber\vspace{2mm}\\
&& -\mathbb{I}_{\omega }(x)\xi (x,t)\lambda _p(x,t)z_p(x,t)-\mathbb{I}_{\omega }(x)v_3(x,t)p(x,t)\lambda _p(x,t)\bigg\} dx \, dt \nonumber \vspace{2mm}\\
&&=\int_{Q_T }\bigg\{ d_2\nabla \lambda_p (x,t)\cdot \nabla z_p(x,t)+\beta_0z_p(x,t)+\beta _2\mathbb{I}_{\omega }(x)\xi (x)z_p(x,t) \nonumber\vspace{2mm}\\
&&+\lambda _k (x,t)z_p(x,t)\frac{2A(x,t)p(x,t)}
{(1+p^2(x,t))^2}f((1-\mathbb{I}_{\omega }(x)\tau (x,t))k(x,t))+\delta _2\lambda_p
(x,t)z_p(x,t)  \nonumber\vspace{2mm}\\
&& +\mathbb{I}_{\omega }(x)\xi (x,t)\lambda _p(x,t)z_p(x,t) \bigg\} dx \, dt .  \nonumber
\end{eqnarray}

It follows that
\begin{eqnarray}\label{sase}
&&\displaystyle -\int_{Q_T}\beta_0z_p (x,t)dx \, dt \!-\beta _2\int_{Q_T}\mathbb{I}_{\omega }(x)\xi (x,t)z_p (x,t)dx \, dt \!  \nonumber\\
\nonumber~ \\
&&=\int_{Q_T} \lambda_k (x,t)z_p(x,t)\frac{2A(x,t)p(x,t)} {(1+p^2(x,t))^2}f((1-\mathbb{I}_{\omega }(x)\tau (x,t))k(x,t))dx \, dt \nonumber\\
\nonumber~ \\
&&\displaystyle +\theta \int_{Q_T}\int_{\Omega }\lambda_p(x,t)f'((1-\mathbb{I}_{\omega }(x')\tau (x',t))k(x',t)) \nonumber\\
\nonumber \\
&&\displaystyle  \times [(1-\mathbb{I}_{\omega }(x')\tau
(x',t))z_k(x',t)-k(x',t)\mathbb{I}_{\omega }(x')v_2(x',t)]\varphi (x',x)dx' \, dx
\,dt  \nonumber\\
\nonumber \\
&&\displaystyle - \int_{Q_T}\mathbb{I}_{\omega }(x)v_3(x,t)p(x,t)\lambda _p(x,t) dx \, dt .
%\label{sase}
\end{eqnarray}
By (\ref{cinci}) and (\ref{sase}) we get that
\begin{eqnarray}
&&\displaystyle \int_0^T\int_{\Omega }(cz_k (x,t)-\beta_0z_p(x,t)-\beta_1\mathbb{I}_{\omega }(x)\tau(x,t) z_k(x,t)-\beta_2\mathbb{I}_{\omega }(x)\xi (x,t) z_p(x,t))dx \, dt \nonumber \\
\nonumber \\
&&\displaystyle =\int_{Q_T} \bigg[-\frac{A(x,t)\lambda _k(x,t)}{1+p^2(x,t)}f'((1-\mathbb{I}_{\omega }(x)\tau(x,t) )k(x,t))k(x,t)\mathbb{I}_{\omega }(x)v_2(x,t) \nonumber \\
\nonumber \\
&&-v_1(x,t)k(x,t)\lambda_k(x,t)\!\bigg]dx \, dt \nonumber\\
\nonumber \\
&&\!\!\!\displaystyle -\theta \int_{Q_T}\int_{\Omega }\!\!\lambda_p(x,t)f'((1-\mathbb{I}_{\omega }(x')\tau
(x',t))k(x',t))k(x',t)\nonumber \\
\nonumber \\
&& \times \mathbb{I}_{\omega }(x')v_2(x',t)\varphi (x',x)dx'  dx \, dt -\int_{Q_T}\mathbb{I}_{\omega }(x)v_3(x,t)p(x,t)\lambda _p(x,t) dx \, dt .\label{sapte}
\end{eqnarray}
By (\ref{patru}) and (\ref{sapte}) we get that
$$dI(u)(v)=\int_{Q_T }k(x,t)[v_1(x,t)(1-\lambda _k(x,t))-v_2(x,t)\mathbb{I}_{\omega }(x)\beta^*(x,t)]dx \, dt $$
$$-\int_{Q_T}p(x,t)v_3(x,t)\mathbb{I}_{\omega }(x)(\beta _2+\lambda _p(x,t))dx \, dt,$$
and so we get the conclusion of the theorem. \vspace{3mm}

In the first paragraph, we describe the conceptual gradient-type algorithm for
the optimal control procedure.
In the second paragraph, we describe the numerical methods adopted for the solutions of System (\ref{dynamics_regional_pollution_treatment_new}) and (\ref{kopt})-(\ref{popt})
with the above mentioned boundary and initial/final conditions.
In the last paragraph, we report and discuss the computational results.

\subsection{Conceptual algorithm}

\begin{itemize}
	\item[STEP 0:] $\qquad $  Set $iter:=0$, Choose maxiter, $\delta
	>0$.
	\\
	\item[STEP 1:] $\qquad $  Initialize $(c^{(iter)},\tau^{(iter)},\xi ^{(iter)})\in \tilde{U}$.
	\\
	\item[STEP 2:]   $\qquad $  Compute
	$k^{(iter)},p^{(iter)},\lambda_k^{(iter)},\lambda
	_p^{(iter)},(\beta^*)^{(iter)}$ corresponding to \\
	$(c^{(iter)},\tau^{(iter)},\xi ^{(iter)})$, and
	$I^{(iter)}=I(c^{(iter)},\tau^{(iter)}, \xi ^{(iter)})$.
	\\
	\item[STEP 3:]  $\qquad $ Choose $(v_1,v_2,v_3)\in L^{\infty
	}(Q_T)\times L^{\infty }(\tilde{Q}_T)\times L^{\infty }(\tilde{Q}_T)$ such that 
	
	\noindent
	$(c^{(iter)},\tau
	^{(iter)},\xi ^{(iter)})+\varepsilon _0(v_1,v_2,v_3)\in \tilde{U}$ $(\varepsilon _0\geq 0)$
	and that $v_1(1-\lambda _k^{(iter)})\geq 0, \
	-v_2\mathbb{I}_{\omega }(\beta^*)^{(iter)}\geq 0$, $-v_3\mathbb{I}_{\omega }(\beta _2+\lambda _p^{(iter)})\geq 0$ a.e. $(x,t)\in Q_T$.
	\\
	\item[STEP 4:] $\qquad $   Compute
	$(c^{(new)},\tau^{(new)}, \xi ^{(new)}):=(c^{(iter)},\tau ^{(iter)}, \xi ^{(iter)})+\varepsilon
	_0(v_1,v_2,v_3)$.
	\\
	\item[STEP 5:] $\qquad $  Compute $\eta _0$ which maximizes
	
	$I(\eta (c^{(iter)},\tau^{(iter)}, \xi ^{(iter)})+(1-\eta )(c^{(new)},\tau^{(new)}, \xi ^{(new)}))$, for $\eta \in [0,1]$.
	
	Set $(c^{(iter+1)},\tau^{(iter+1)},\xi ^{(iter+1)})=\eta _0(c^{(iter)},\tau^{(iter)}, \xi ^{(iter)})$ \\
	$+(1-\eta _0)(c^{(new)},\tau^{(new)},\xi ^{(new)})$.
	\\
	
	\item[STEP 6:]  $\qquad $   If $I^{(iter+1)}-I^{(iter)}<\delta $,
	then
	
	\ \ \ \ \ $(c^*,\tau^*, \xi ^*):= (c^{(iter)},\tau^{(iter)},\xi ^{(iter)})$
	
	elseif $iter>maxiter$, then
	
	\ \ \ \ \ $(c^*,\tau^*,\xi ^*):= (c^{(iter)},\tau^{(iter)},\xi ^{(iter)})$
	
	else
	
	\ \ \ \ \ $iter:=iter+1$ and GO TO STEP 3.
\end{itemize}

In the numerical tests, we will use:
$$v_1(x,t)\in sgn (1-\lambda _k^{(iter)}(x,t)),$$
$$v_2(x,t)\in -sgn (\beta^*(x,t))^{(iter)}, \quad v_3(x,t)\in -sgn (\beta_2+\lambda _p^{(iter)}(x,t)),$$
a.e. $(x,t)\in Q_T$, and a.e. $(x,t)\in \tilde{Q}_T$, respectively.

\section{Numerical methods} \label{numerics}

Equations (\ref{dynamics_regional_pollution_treatment_new}) and (\ref{kopt})-(\ref{popt}) are discretized by finite elements in space and finite differences in time.

\subsection{Space discretization.} We first apply a standard Galerkin
procedure to the weak formulations of systems (\ref{dynamics_regional_pollution_treatment_new}) and (\ref{kopt})-(\ref{popt}), with the associated boundary and initial/ final conditions, respectively.
The two-dimensional set $\overline{\Omega }=[-1,1]^2$ is discretized by a
uniform grid of $64\times 64$ bilinear finite elements ($\uQ^1$), yielding a
total amount of $N=4225$ discretization nodes. The stiffness matrix
is computed exactly, whereas the mass matrix is obtained by
applying the mass lumping technique.

In the following, we describe in detail this procedure applied to System (\ref{dynamics_regional_pollution_treatment_new}),
being the treatment of equations (\ref{kopt})-(\ref{popt}) analogous.

Let us first denote by $V$ the Sobolev space $H^1(\Omega)$, by
$(\cdot,\cdot)$ the $L^2 -$ inner product and by
$a_1(\cdot,\cdot),\,a_2(\cdot,\cdot):V\times V\rightarrow
\mathbb{R}$ the bilinear forms
\[
\begin{array}{ll}
\displaystyle a_1(u,v):=\int_\Omega d_1 \nabla u \cdot \nabla v dx, & \quad \textnormal{for all} \,\, u,v\in V, \vspace{0.2cm}\\
\displaystyle a_2(u,v):=\int_\Omega d_2 \nabla u \cdot \nabla v
dx, & \quad \textnormal{for all}\,\, u,v\in V.
\end{array}
\]
Then, using  Green's formula, the variational formulation of
problem (\ref{dynamics_regional_pollution_treatment_new}) reads:

given $(c,\,\tau)\in U,$ find $k,\,p$ such that
\begin{equation}\label{sys_var_tax}
\left\{
\begin{array}{ll}
\displaystyle \frac{\partial}{\partial t} (k(t),\phi) + a_1(k(t),\phi) = (G_1(k(t),p(t),c(t),\tau(t)),\phi), & \  \textnormal{for all} \ \phi\in V, \vspace{0.2cm}\\
\displaystyle \frac{\partial}{\partial t} (p(t),\phi) +
a_2(p(t),\phi) = (G_2(k(t),p(t),c(t),\tau(t)),\phi), & \
\textnormal{for all} \ \phi\in V,
\end{array}
\right.
\end{equation}
subject to the initial conditions on $k$ and $p$, and where $G_1$ and $G_2$ are the reaction terms in (\ref{dynamics_regional_pollution_treatment_new}).

Let us now introduce the $\uP^1$ finite element space $V^h$,
associated to a partition of the one-dimensional domain $\Omega$
into $N-1$ subintervals of size $h$. 

The space semi-discretization of problem (\ref{sys_var_tax}) then reads:

find $k_h,\,p_h:[0,T]\rightarrow V^h$ such that
\[
\left\{
\begin{array}{ll}
\displaystyle \frac{\partial}{\partial t} (k_h(t),\phi_h) + a_1(k_h(t),\phi_h) = (G_1(k_h(t),p_h(t),c_h(t),\tau_h(t)),\phi_h), & \  \textnormal{for all} \ \phi_h\in V^h, \vspace{0.2cm}\\
\displaystyle \frac{\partial}{\partial t} (p_h(t),\phi_h) +
a_2(p_h(t),\phi_h) = (G_2(k_h(t),p_h(t),c_h(t),\tau_h(t)),\phi_h),
& \  \textnormal{for all} \ \phi_h\in V^h,
\end{array}
\right.
\]
where $c_h,\,\tau_h:[0,T]\rightarrow V^h$ are convenient finite element approximations of $c$ and $\tau$.
The previous variational equation is equivalent to
\begin{equation}\label{sys_semid_tax}
\left\{
\begin{array}{ll}
\displaystyle \frac{\partial}{\partial t} (k_h(t),\phi_i) + a_1(k_h(t),\phi_i) = (G_1(k_h(t),p_h(t),c_h(t),\tau_h(t)),\phi_i), \vspace{0.2cm}\\
\displaystyle \frac{\partial}{\partial t} (p_h(t),\phi_i) +
a_2(p_h(t),\phi_i) =
(G_2(k_h(t),p_h(t),c_h(t),\tau_h(t)),\phi_i),\,
\end{array}
\right.
\end{equation}
for all $i=1,...,N$, where $\phi_i$, $i=1,...,N$, are the {\em hat
	basis functions} of the $\uP^1$ finite element space $V^h$.

Expanding the unknowns $k_h(x,t),\,p_h(x,t)$ and $c_h(x,t),\,\tau_h(x,t)$ with respect to the basis $\{\phi_i\}_{i=1}^N$,
we have
\begin{equation}\label{sol_expand_hat_tax}
\begin{array}{ll}
\displaystyle k_h(x,t)=\sum_{j=1}^N K_j(t) \phi_j(x),& \displaystyle \quad p_h(x,t)=\sum_{j=1}^N P_j(t), \vspace{0.2cm}\\
\displaystyle c_h(x,t)=\sum_{j=1}^N C_j(t) \phi_j(x),& \displaystyle \quad \tau_h(x,t)=\sum_{j=1}^N T_j(t)
\end{array}
\end{equation}
where $K_j,\,P_j:[0,T]\rightarrow\mathbb{R}$ are unknown scalar functions to be determined, and
we assume that $C_j,\,T_j:[0,T]\rightarrow\mathbb{R}$ are given.

Denoting by
$\mathrm{A}_1,\,\mathrm{A}_2$ the stiffness matrices associated
with the bilinear forms $a_1(\cdot,\cdot),\,a_2(\cdot,\cdot)$, by
$\mathrm{M}$ the lumped mass matrix,
and substituting the expressions
(\ref{sol_expand_hat_tax}) into system (\ref{sys_semid_tax}), the
semi-discrete problem then reads:

find $\uK,\uP:[0,T]\rightarrow \mathbb{R}^N$ such that
\begin{equation}\label{sys_semid_alg_tax}
\left\{
\begin{array}{l}
\displaystyle \mathrm{M}\frac{d\uK}{dt}(t) + \mathrm{A}_1 \uK(t) = \mathrm{M}\uG_1(\uK(t),\uP(t),\uC(t),\uTau(t)), \vspace{0.2cm}\\
\displaystyle \mathrm{M}\frac{d\uP}{dt}(t) + \mathrm{A}_2 \uP(t) =
\mathrm{M}\uG_2(\uK(t),\uP(t),\uC(t),\uTau(t)),
\end{array}
\right.
\end{equation}
where
\[
\begin{array}{lll}
\uK(t)                                & =     & [...,K_j(t),...]^T,  \vspace{0.2cm}\\
\uP(t)                                & =     & [...,P_j(t),...]^T, \vspace{0.2cm}\\
\uC(t)                                & =     & [...,C_j(t),...]^T,  \vspace{0.2cm}\\
\uTau(t)                              & =     & [...,T_j(t),...]^T,  \vspace{0.2cm}\\
\uG_1(\uK(t),\uP(t),\uC(t),\uTau(t))  & =     & [...,G_1(K_j(t),P_j(t),C_j(t),T_j(t)),...]^T,  \vspace{0.2cm}\\
\uG_2(\uK(t),\uP(t),\uC(t),\uTau(t))  & =     & [...,G_2(K_j(t),P_j(t),C_j(t),T_j(t)),...]^T.
\end{array}
\]

\subsection{Time discretization.} After space discretization by finite elements, we obtain the semi-discrete problem (\ref{sys_semid_alg_tax}) that
consists of a system of ODEs. We solve this system by employing a
first order semi-implicit finite difference scheme, 
where the linear
diffusion and reaction terms are approximated by Backward Euler,
whereas the non-linear reaction terms are approximated by Forward
Euler.

In formulae, denoting by $\Delta t$ the time step size, at the generic time step $n$,

given $\uK^n,\,\uP^n$, we have to find $\uK^{n+1},\,\uP^{n+1}$ such that
\begin{equation*}\label{sys_d_alg_tax}
\left\{
\begin{array}{l}
\displaystyle \mathrm{M}\frac{\uK^{n+1}-\uK^n}{\Delta t} + \mathrm{A}_1 \uK^{n+1} = \mathrm{M}\uG_1(\uK^n,\uP^n,\uC^n,\uTau^n), \vspace{0.2cm}\\
\displaystyle \mathrm{M}\frac{\uP^{n+1}-\uP^n}{\Delta t} + \mathrm{A}_2 \uP^{n+1} = \mathrm{M}\uG_2(\uK^n,\uP^n,\uC^n,\uTau^n).
\end{array}
\right.
\end{equation*}
As a result, at each time step it is required the solution of the linear system of algebraic equations
\[
\left[
\begin{array}{cc}
\displaystyle \frac{1}{\Delta t} \mathrm{M}+\mathrm{A}_1 & 0 \vspace{0.2cm} \\
0 & \displaystyle \frac{1}{\Delta t} \mathrm{M}+\mathrm{A}_2
\end{array}
\right]
\left[
\begin{array}{cc}
\uK^{n+1} \vspace{0.2cm} \\
\uP^{n+1}
\end{array}
\right]
=
\left[
\begin{array}{cc}
\displaystyle \frac{1}{\Delta t} \mathrm{M} \uK^n + \mathrm{M}\uG_1(\uK^n,\uP^n,\uC^n,\uTau^n) \vspace{0.2cm} \\
\displaystyle \frac{1}{\Delta t} \mathrm{M} \uP^n + \mathrm{M}\uG_2(\uK^n,\uP^n,\uC^n,\uTau^n)
\end{array}
\right],
\]
which in our specific case has the dimension of $8450$ degrees of
freedom. The linear system is solved by the Gaussian elimination
with the built-in function of MATLAB. The time step size is
$\Delta t=0.05$.

\section{Numerical simulations}\label{sec:3}

\begin{figure}[!htb]
	\begin{center}
		\includegraphics[width=6cm]{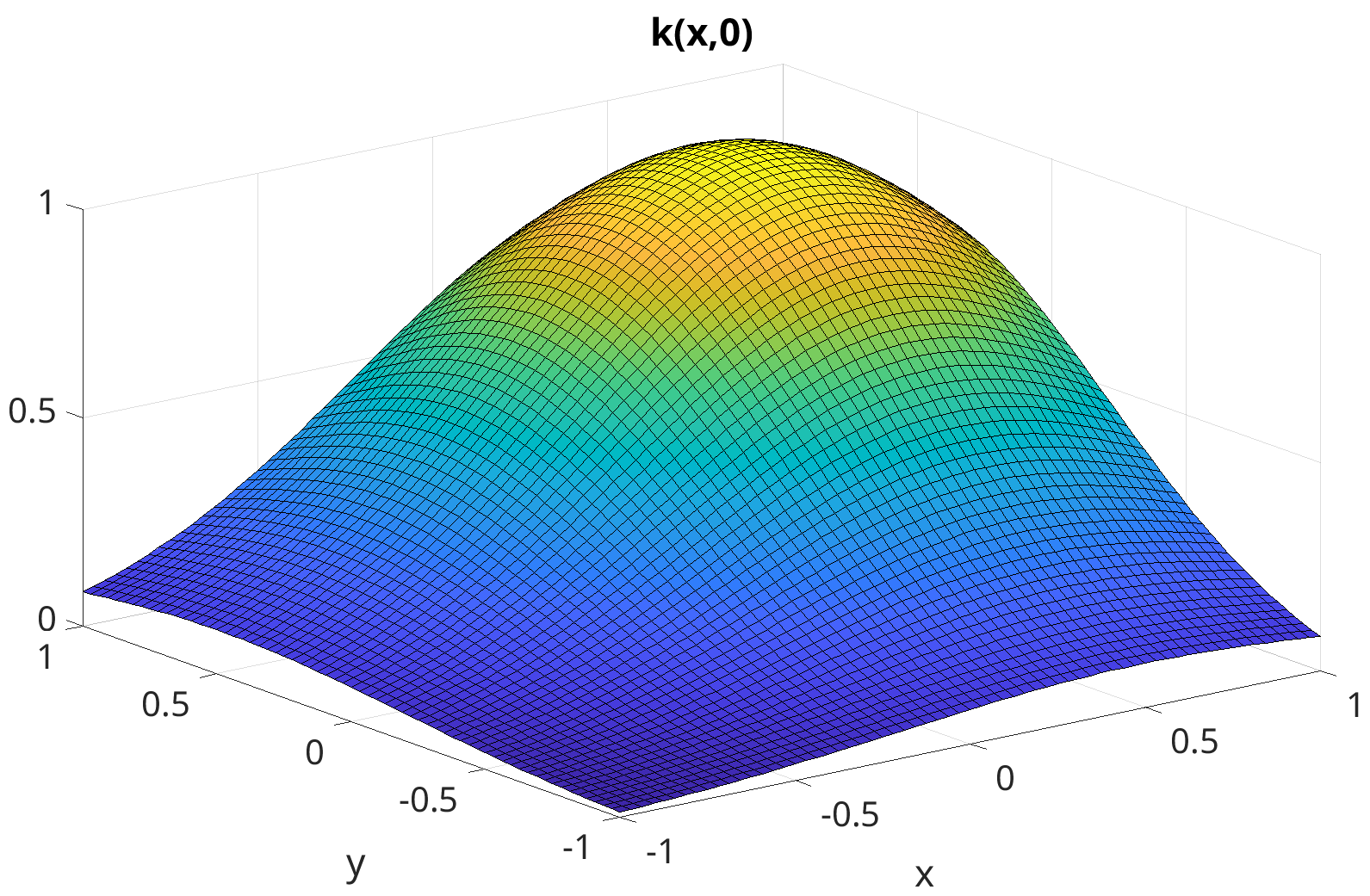}
		\includegraphics[width=6cm]{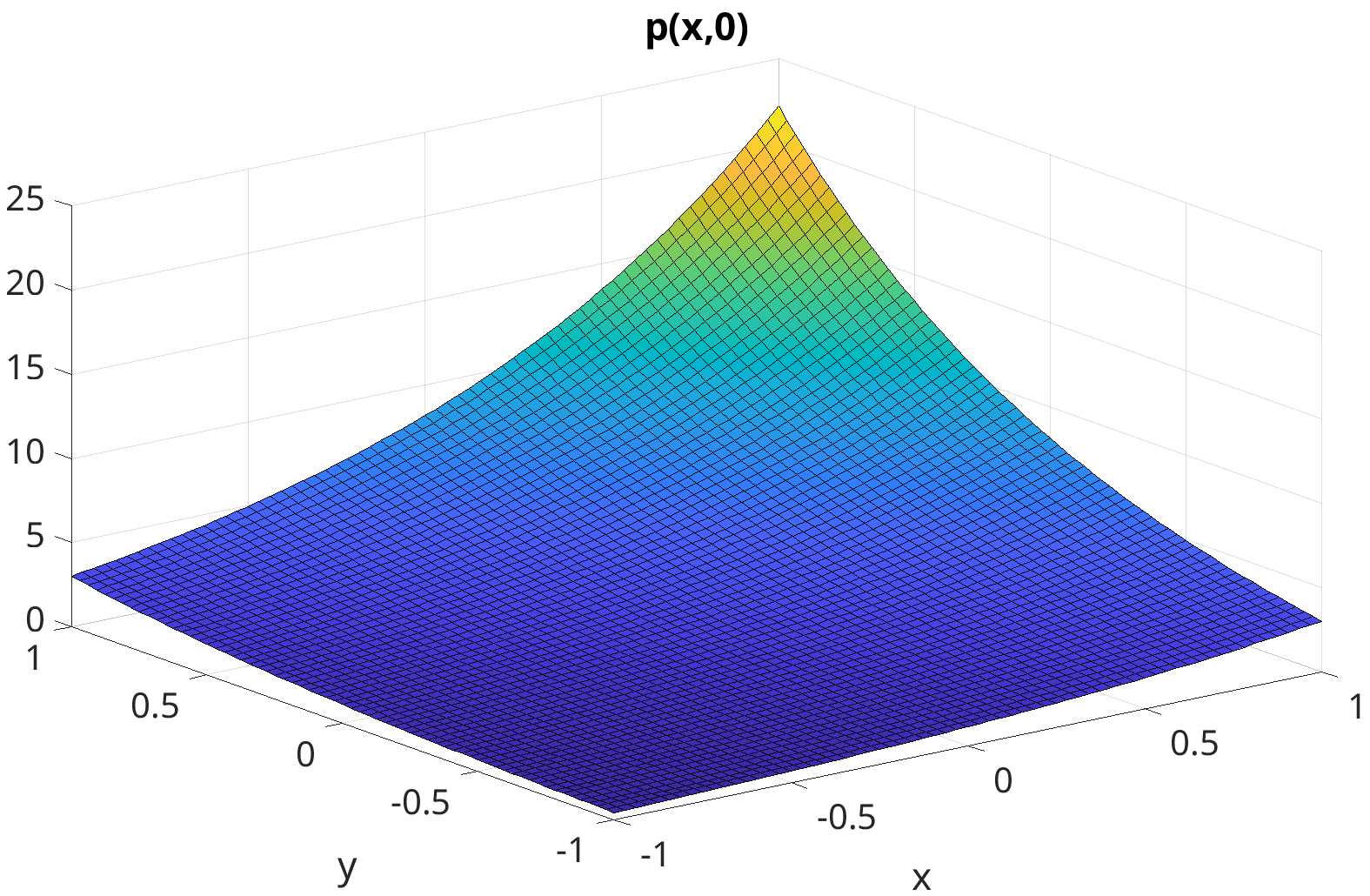}
		\caption{Initial distributions (at $t=0$) of state variables $k(x,t)$ (left) and $p(x,t)$ (right).}
		\label{fig_init_k_p}
	\end{center}
\end{figure}

Let us use the above algorithm with the following parameters'
values $d_1 = 1$,  $d_2 = 1$,  $\delta_1 = 0.05$, $\delta_2 =
0.03$, $A=1$, $s = 0.6$, $\theta = 2$, $\alpha_1 = 0.7$, $\alpha_2
= 1$, $\gamma = 4$, $\varphi = 0.3$, and initial conditions
\[
k_0(x_1, x_2) = 0.1 \exp\left\{ \dfrac{-(x_1-0.5)^2-(x_2-0.5)^2}{0.1}\right\} 
\]
and
\[
p_0(x_1, x_2) = \exp \{ x_1+x_2\}; 
\]
see Fig. \ref{fig_init_k_p}. We  have taken  $\beta_0=\beta_1=\beta_2=1$ and $L=0.5$.

\begin{figure}[!htb]
	\begin{center}
		\includegraphics[width=6cm]{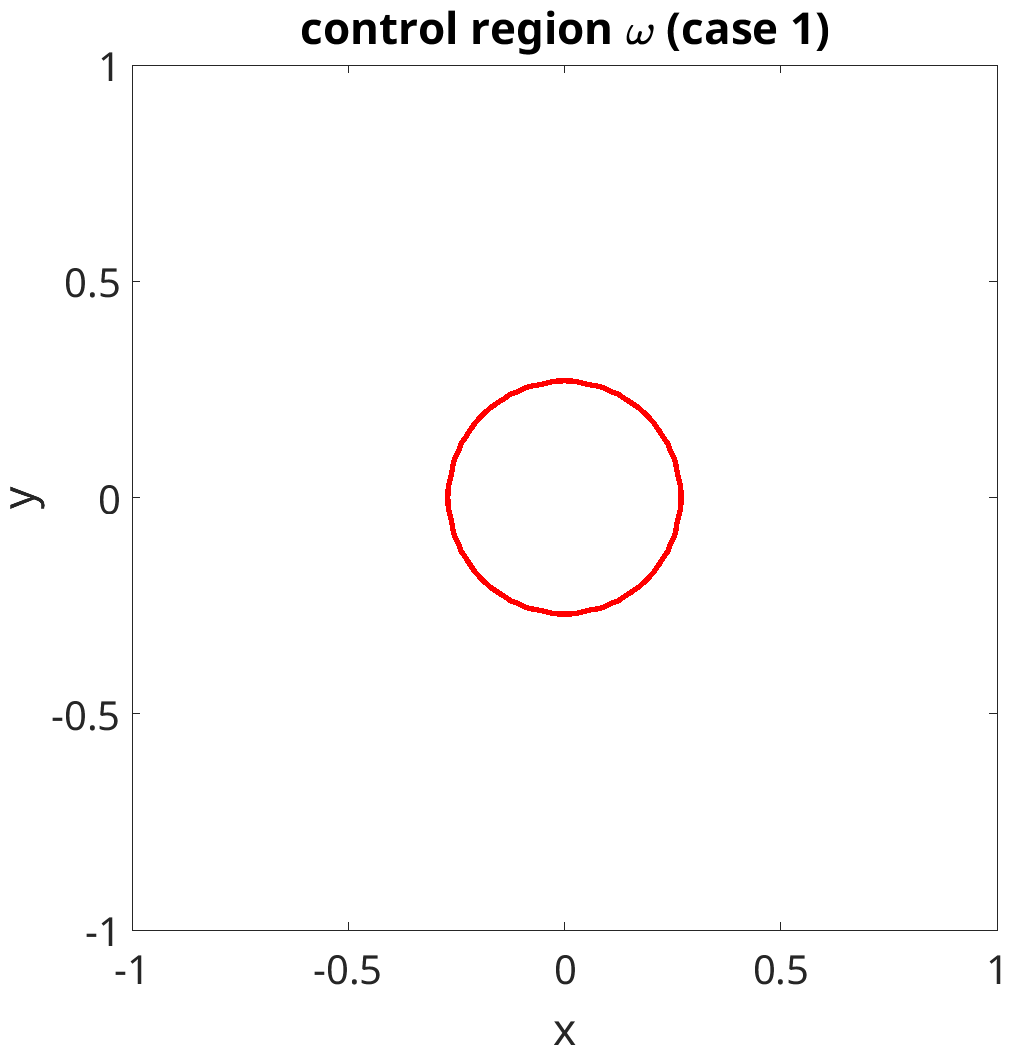}
		\includegraphics[width=6cm]{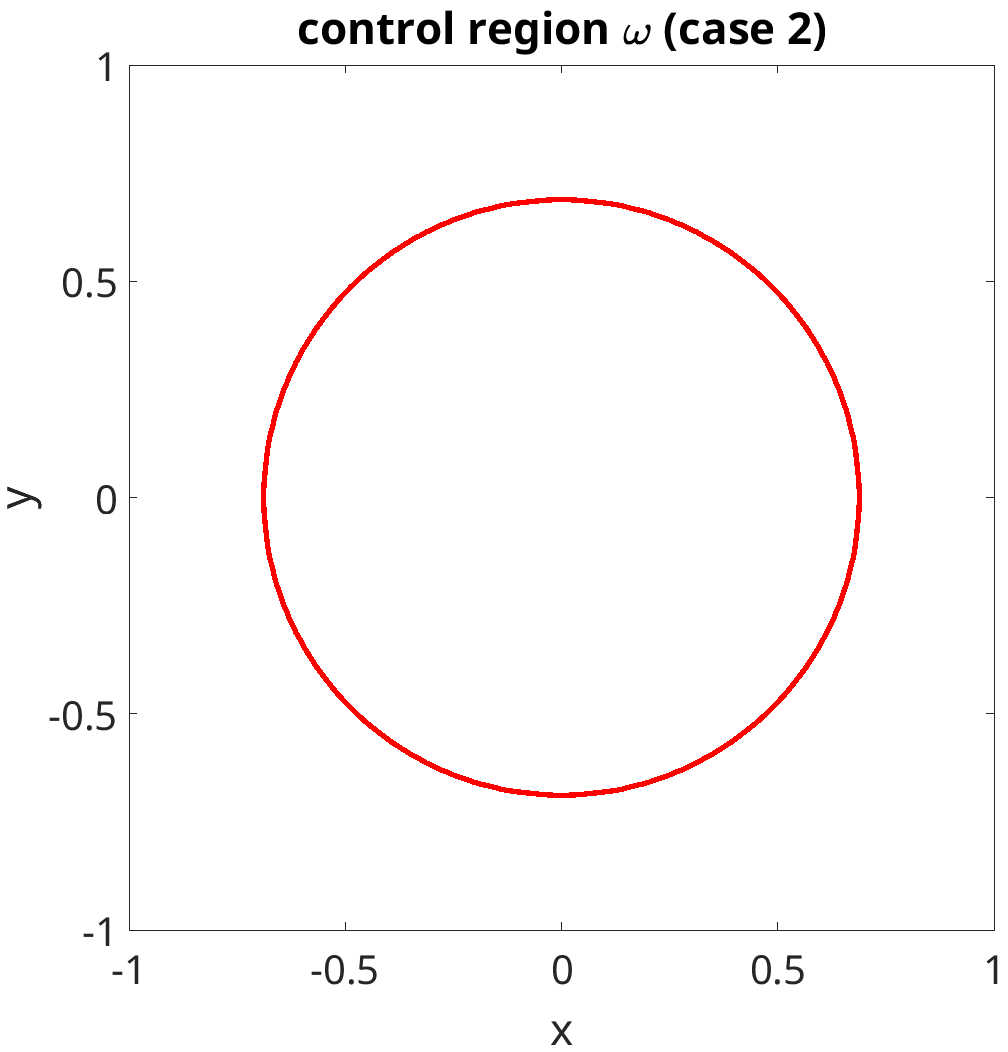}
		\caption{Control regions $\omega$: small ({\bf case 1}, left) and big ({\bf case 2}, right).}
		\label{fig_reg_contr}
	\end{center}
\end{figure}

We consider the following four scenarios; as  time horizon for the control we have taken T=5:
\begin{itemize}
	\item {\bf case 0}: no control is applied;
	\item {\bf case 1}: the control region $\omega$ is the small circle displayed in Fig. \ref{fig_reg_contr} (left);
	\item {\bf case 2}: the control region $\omega$ is the big circle displayed in Fig. \ref{fig_reg_contr} (right);
	\item {\bf case 3}: the control region $\omega$ is the whole domain $\Omega=[-1,1]^2$.
\end{itemize}

%\newpage

\section{Concluding remarks} \label{concluding_remarks}

The numerical results/tests agree with what intuition tells us. However, besides of what agrees with the intuition, we get additional information.
Note that the numerical results give also an evaluation of the optimal control magnitude ($c^*, \tau ^*, \xi ^*)$ and of the corresponding states $k^*(x, t), p^*(x, t).$

\begin{figure}[!htb]
	\begin{center}
		\includegraphics[width=4.8cm]{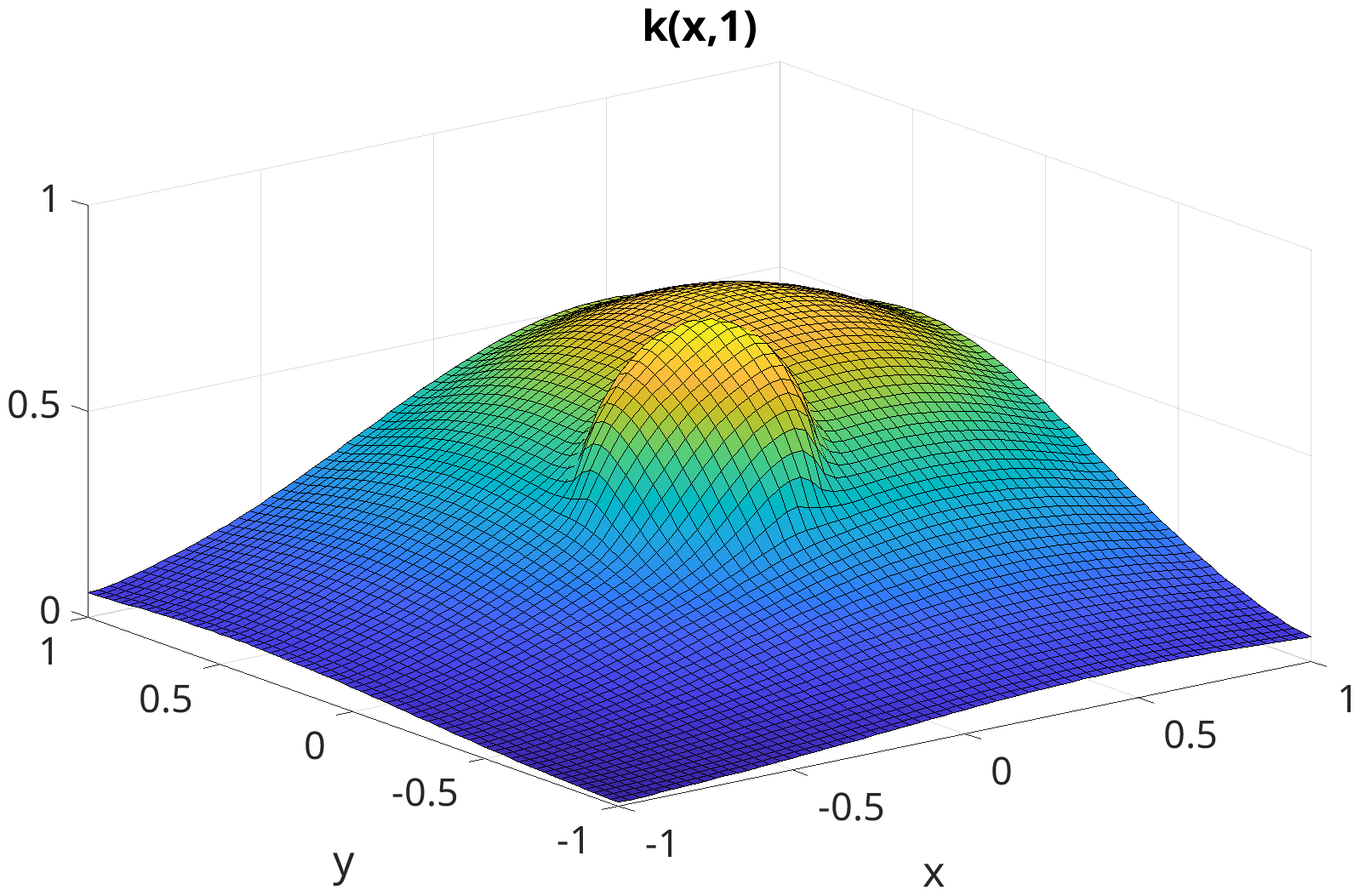}
		\includegraphics[width=4.8cm]{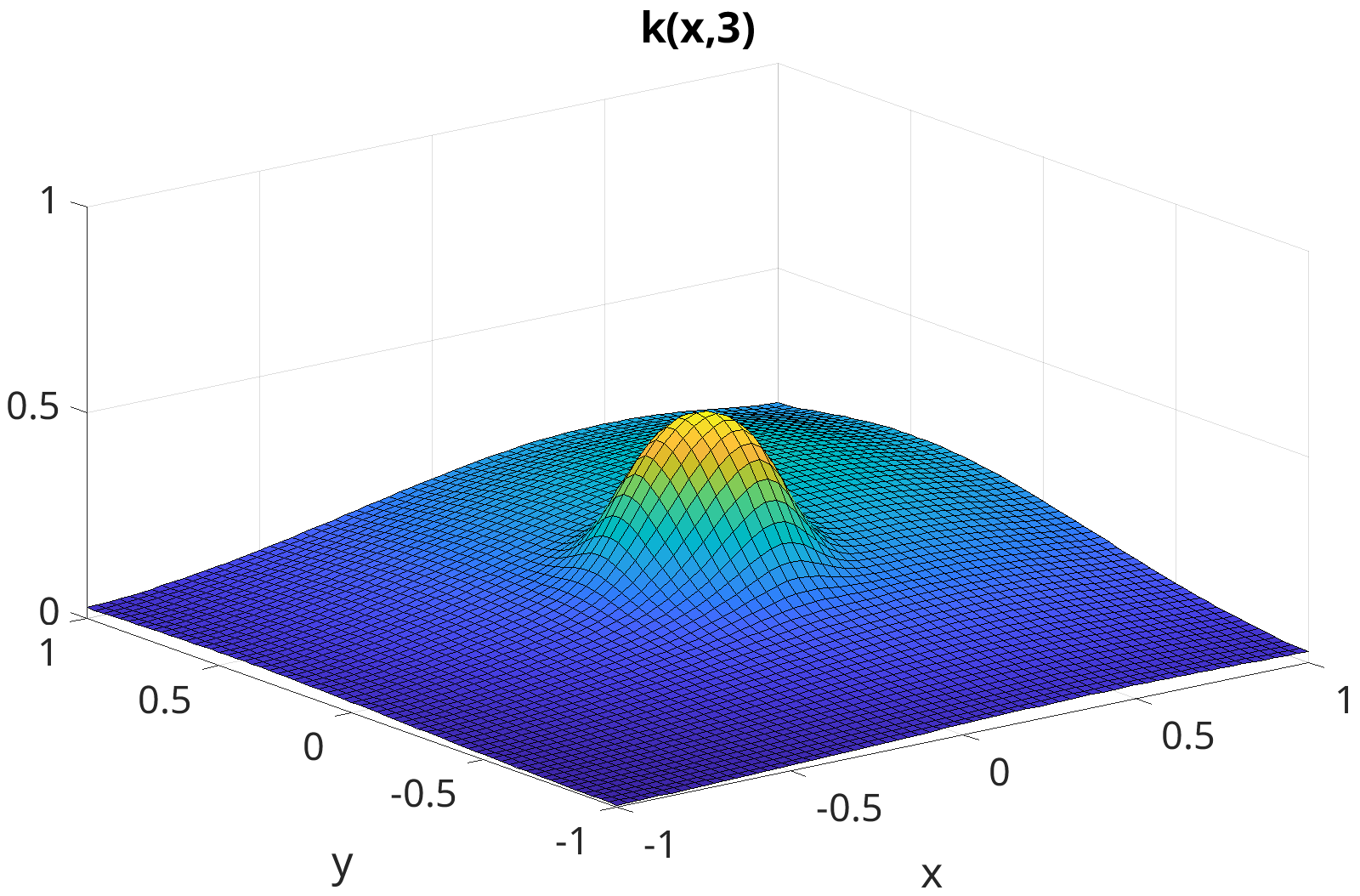}
		\includegraphics[width=4.8cm]{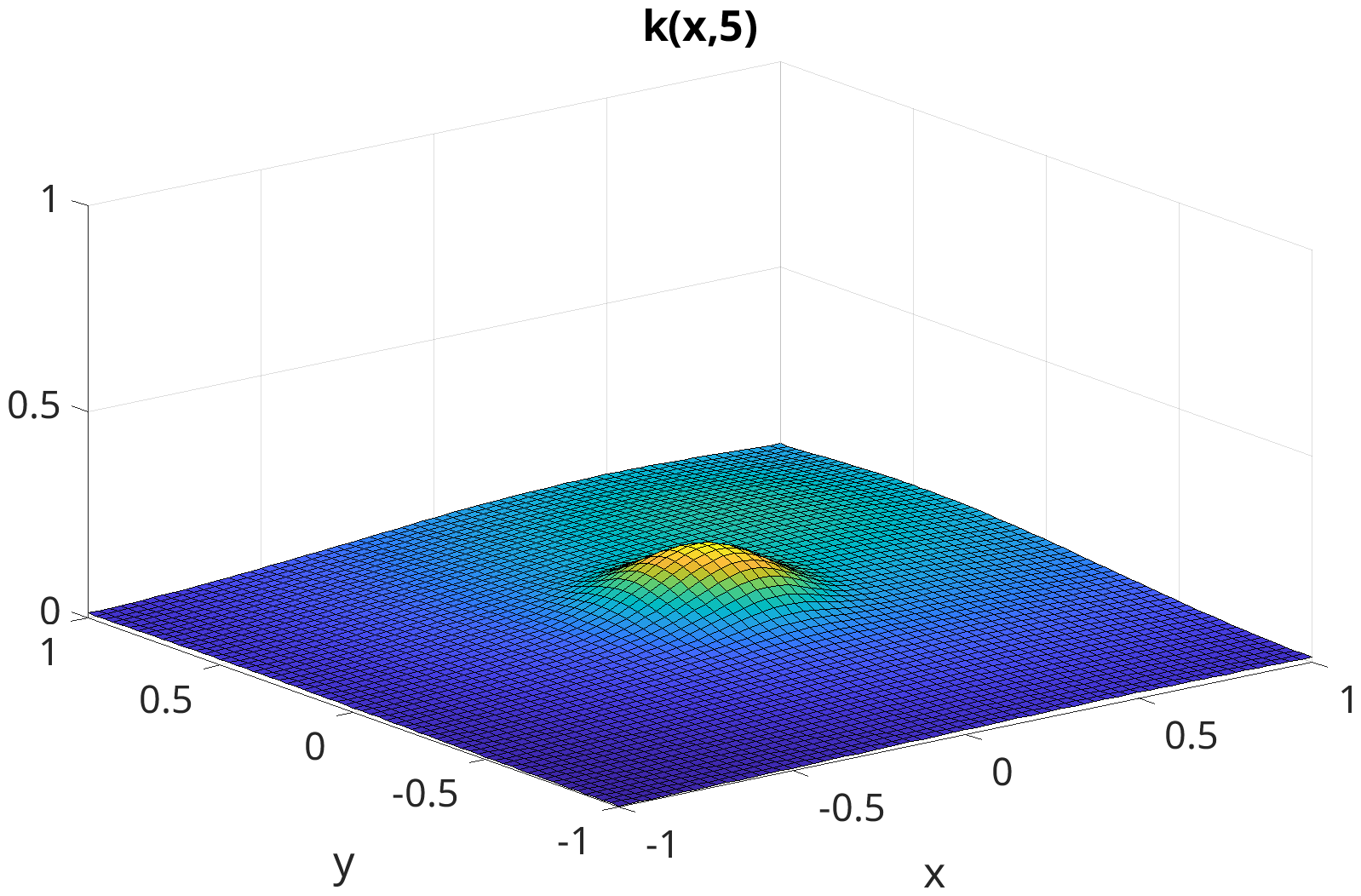}\\
		\includegraphics[width=4.8cm]{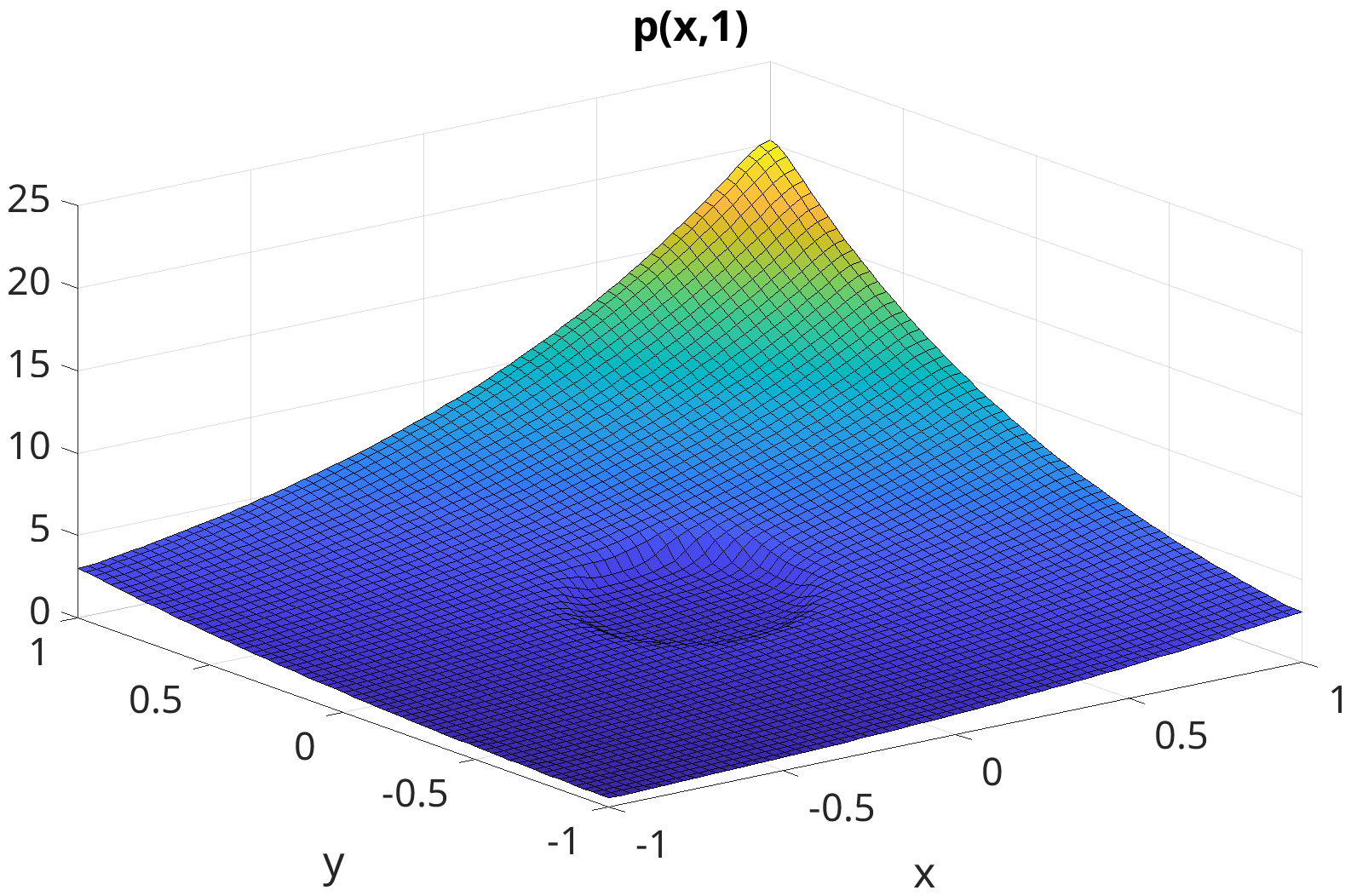}
		\includegraphics[width=4.8cm]{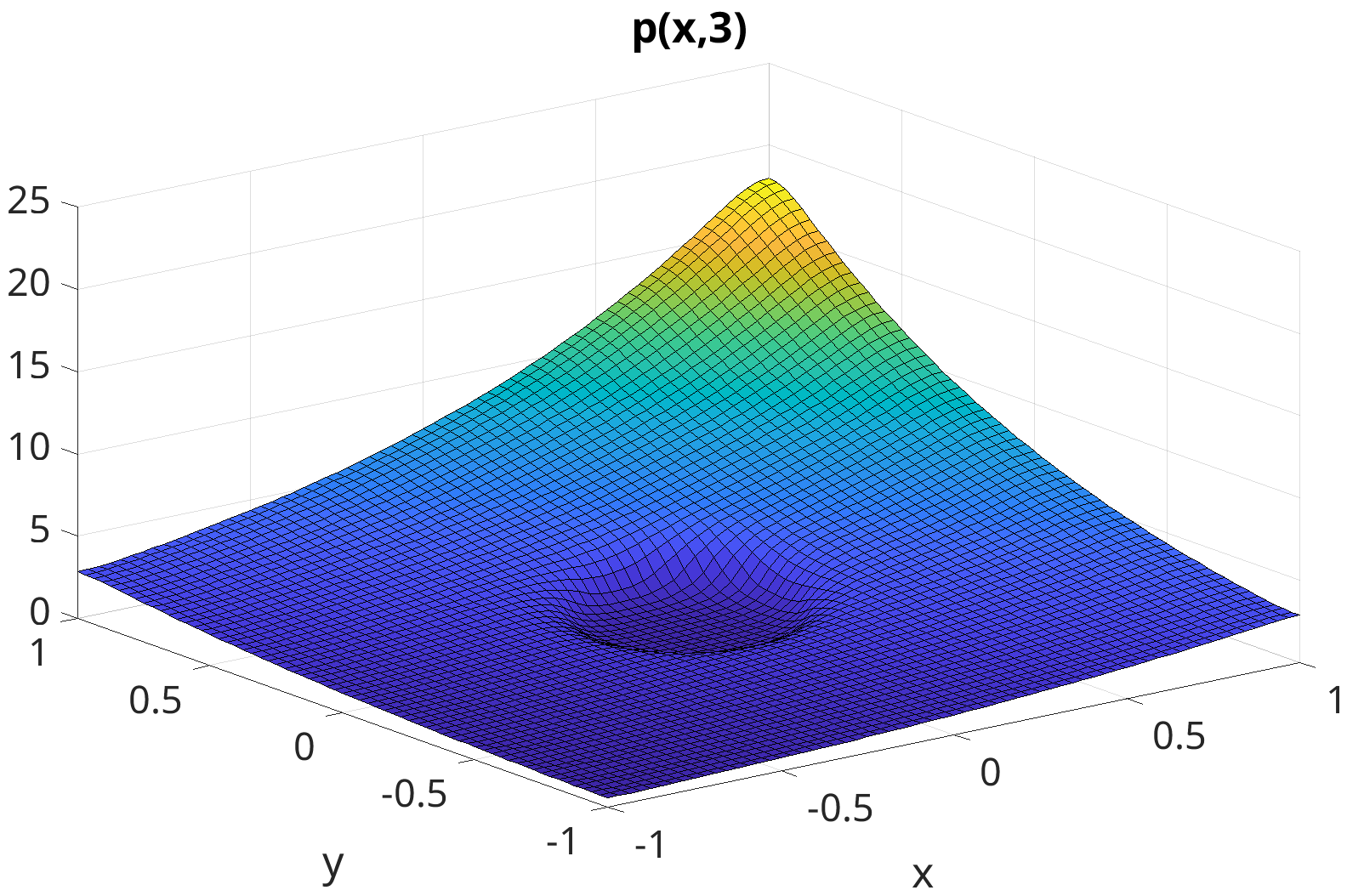}
		\includegraphics[width=4.8cm]{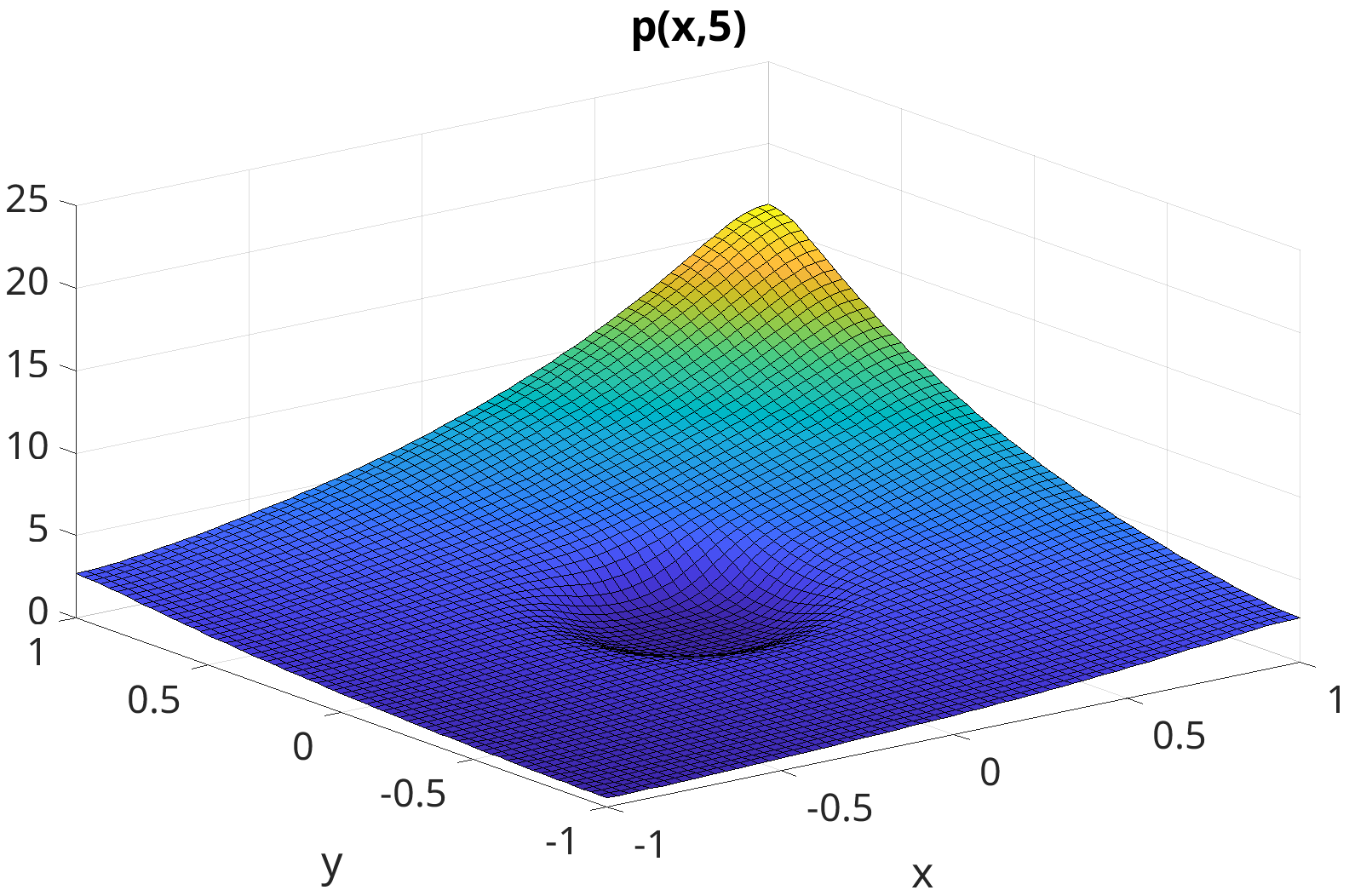}
		\caption{{\bf Case 1}, small control region. Space distributions of the optimal state variables $k(x, t)$ (first row) and $p(x, t)$ (second row) at three time instants.}
		\label{fig_reg_contr_small_k_p}
	\end{center}
\end{figure}

Figure \ref{fig_reg_contr_small_k_p}. Since our main goal is to get a maximal value of the cost functional it is natural noticing a decrease of the capital when
time increases (because there is no interest in keeping a high reserve of $k$ for time larger than the time horizon $T$).

\begin{figure}[!htb]
	\begin{center}
		\includegraphics[width=4.8cm]{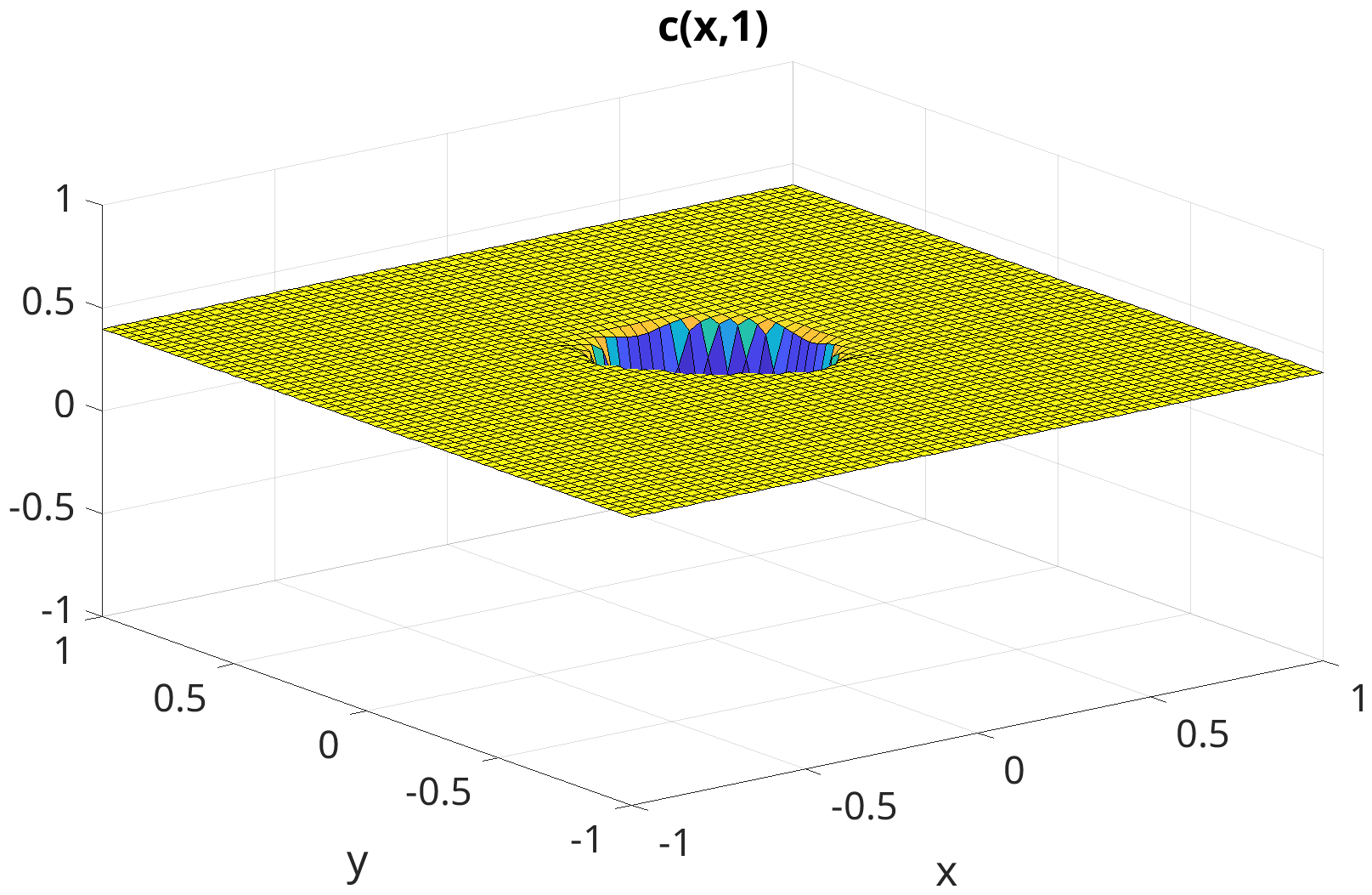}
		\includegraphics[width=4.8cm]{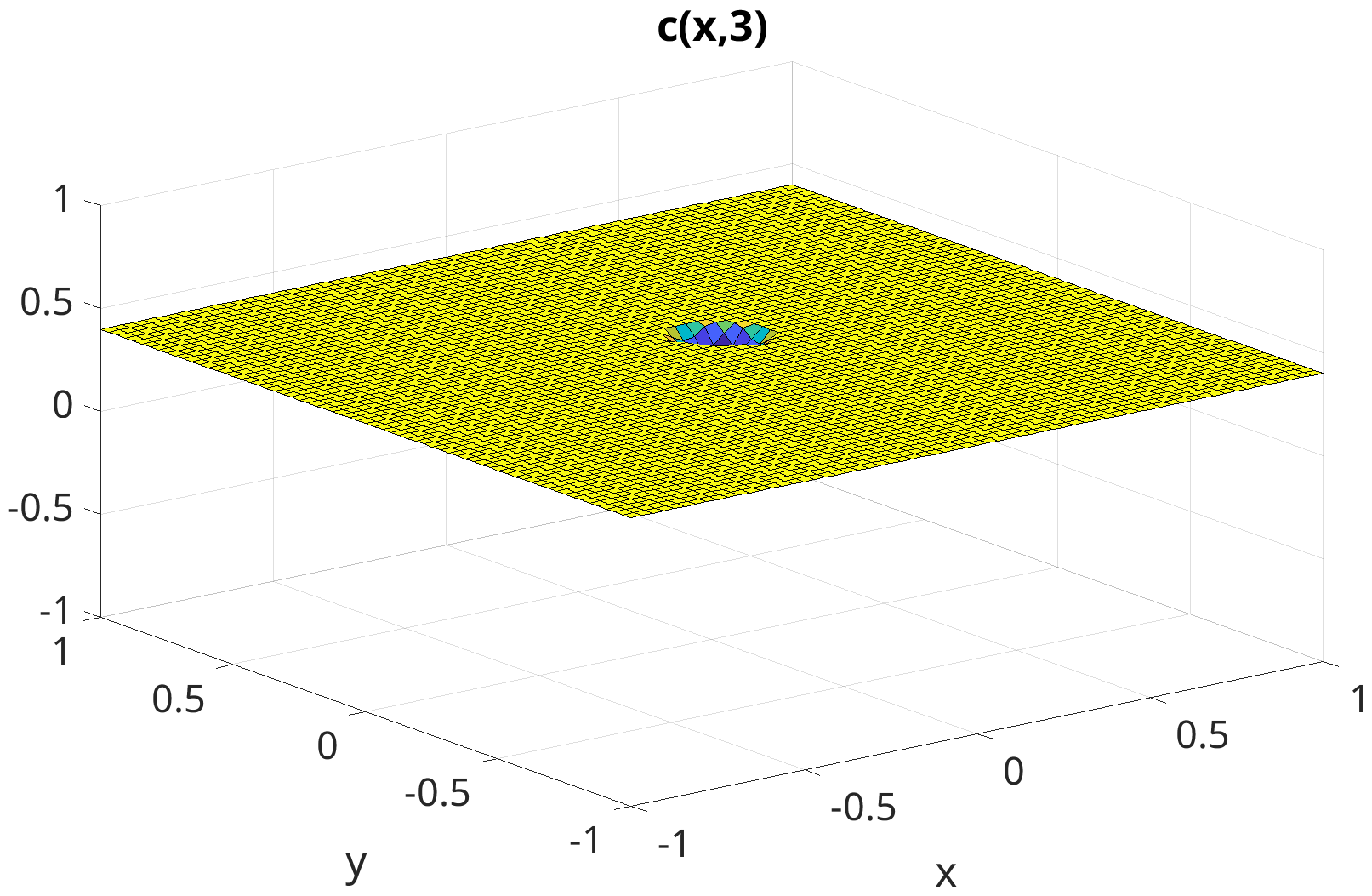}
		\includegraphics[width=4.8cm]{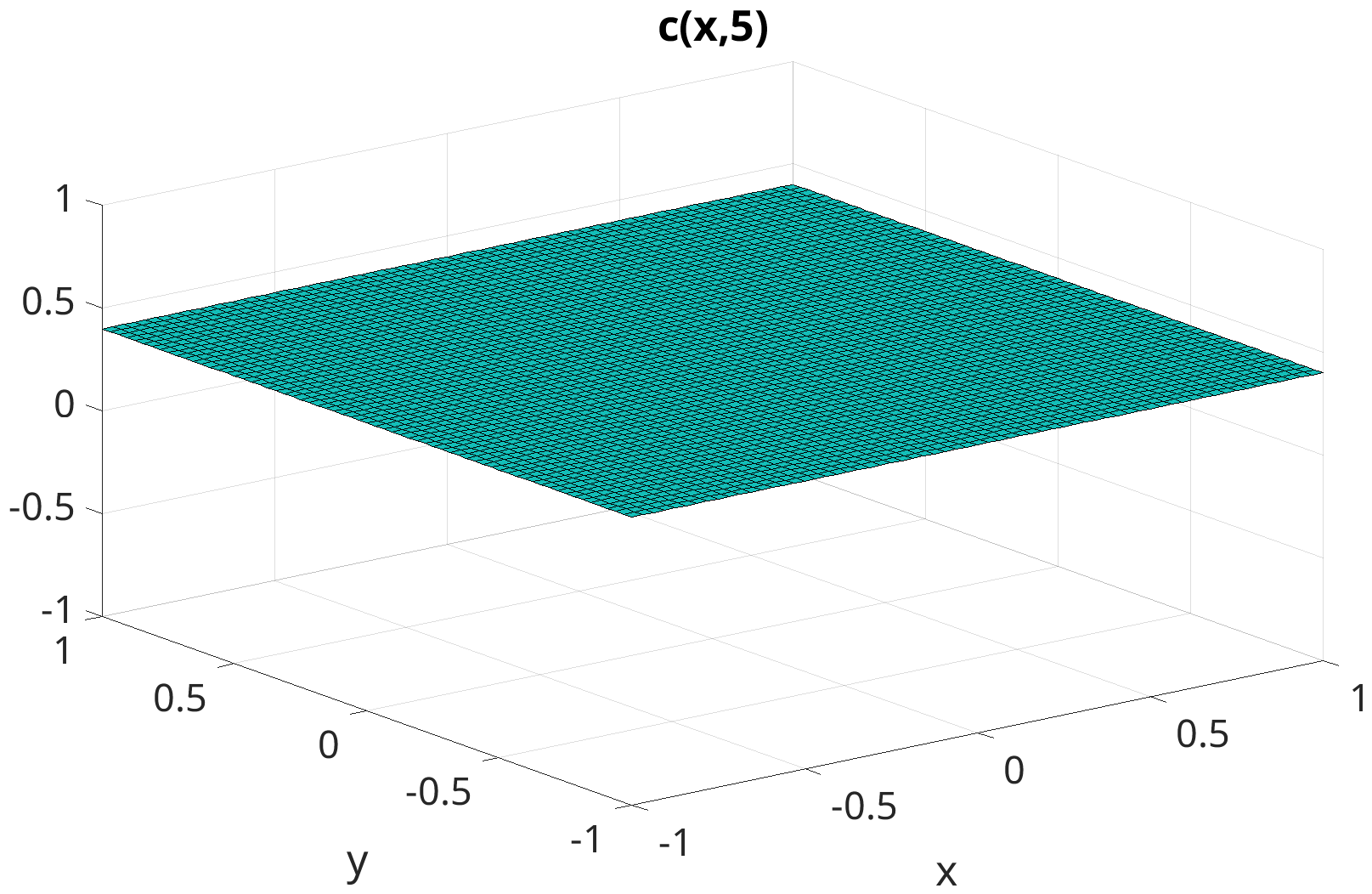}\\
		\includegraphics[width=4.8cm]{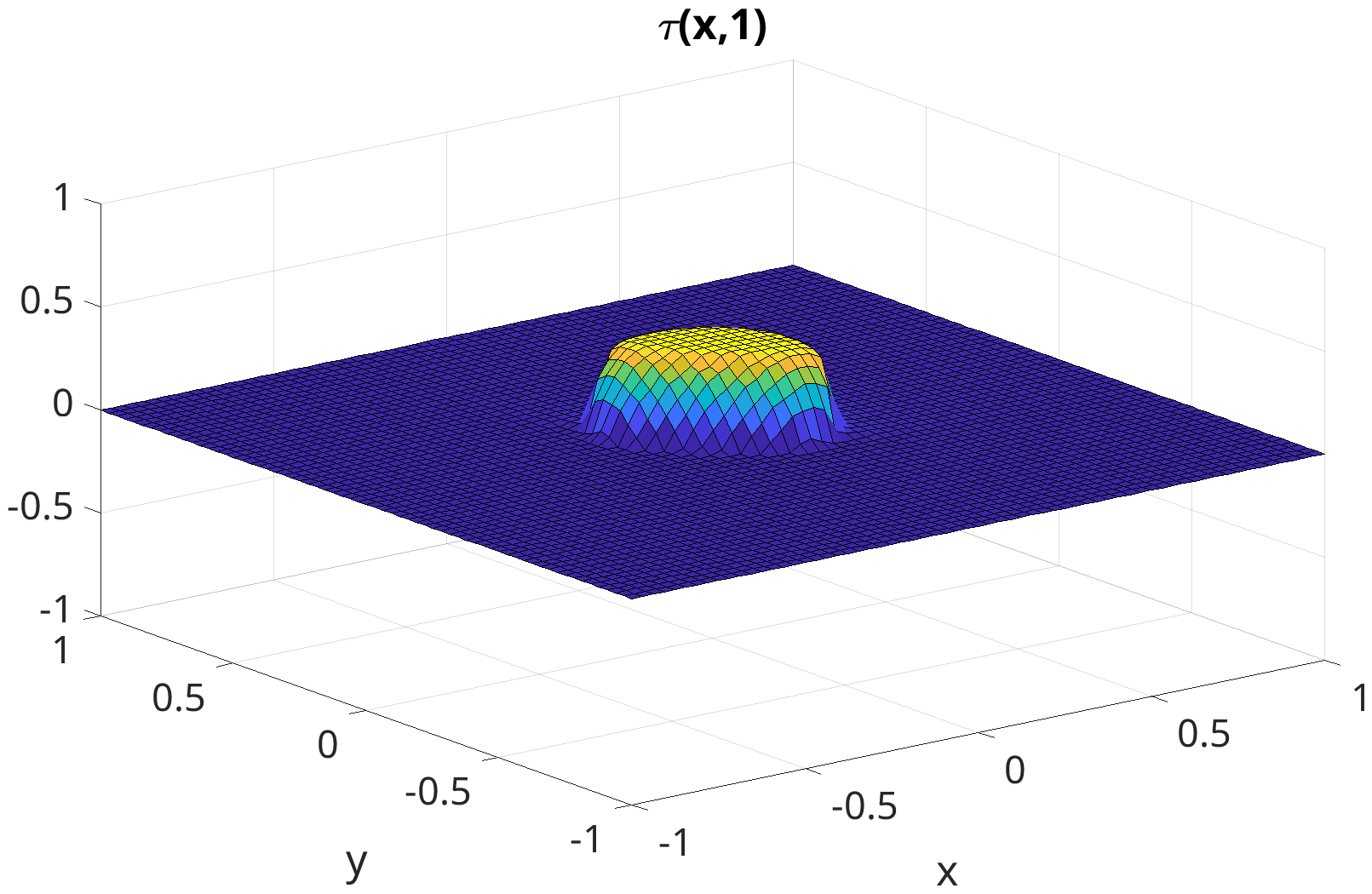}
		\includegraphics[width=4.8cm]{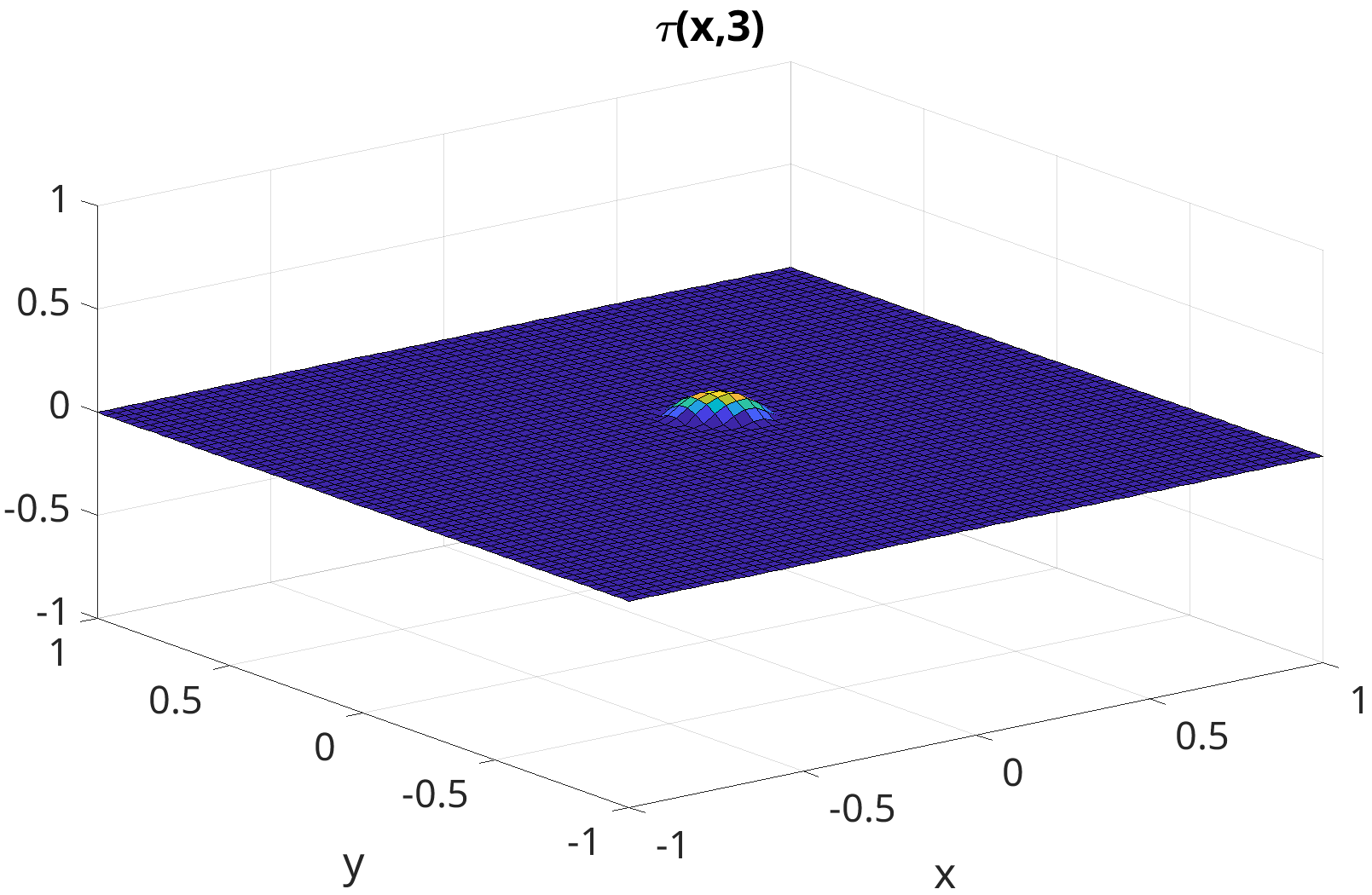}
		\includegraphics[width=4.8cm]{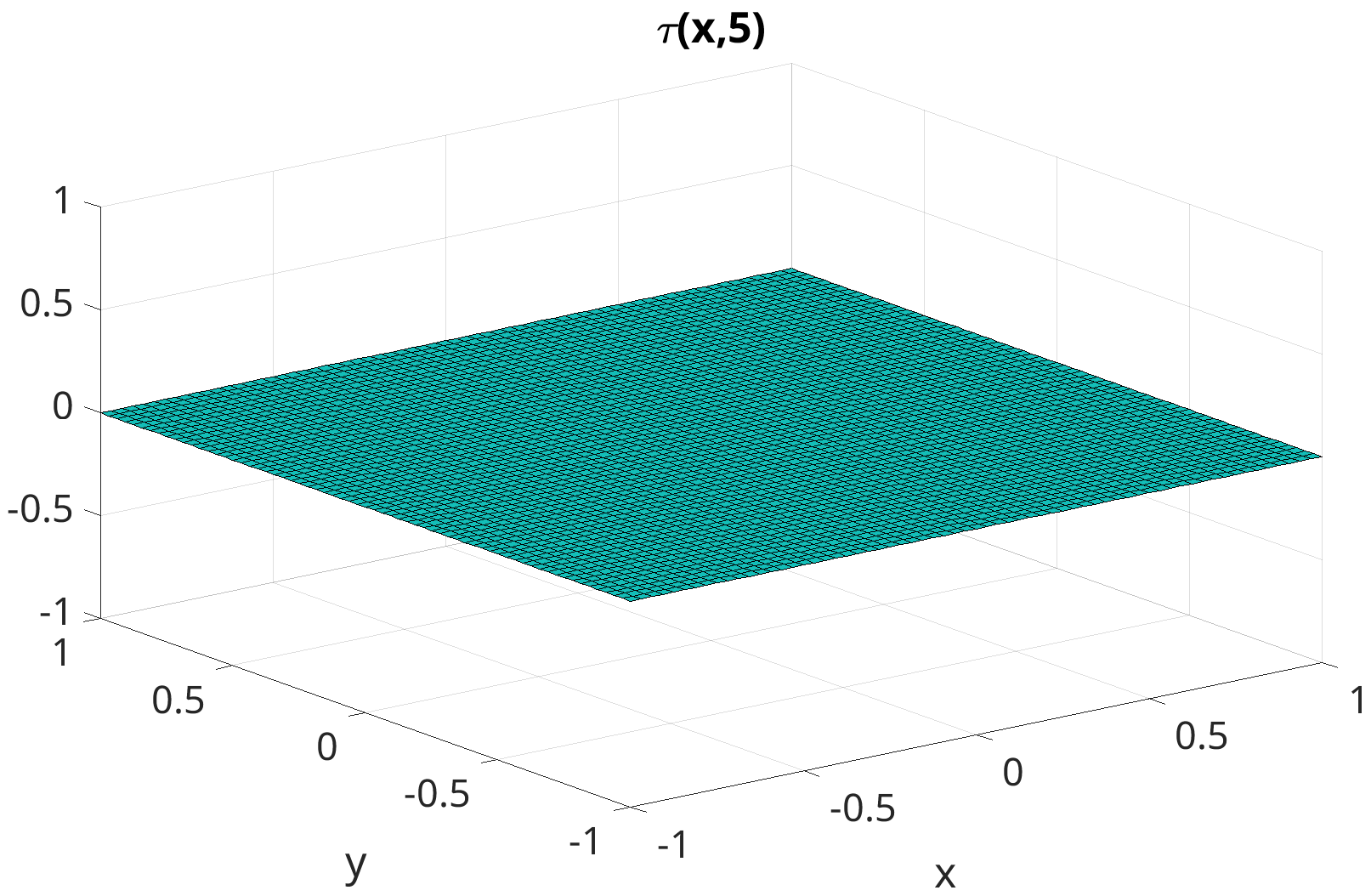}
		\includegraphics[width=4.8cm]{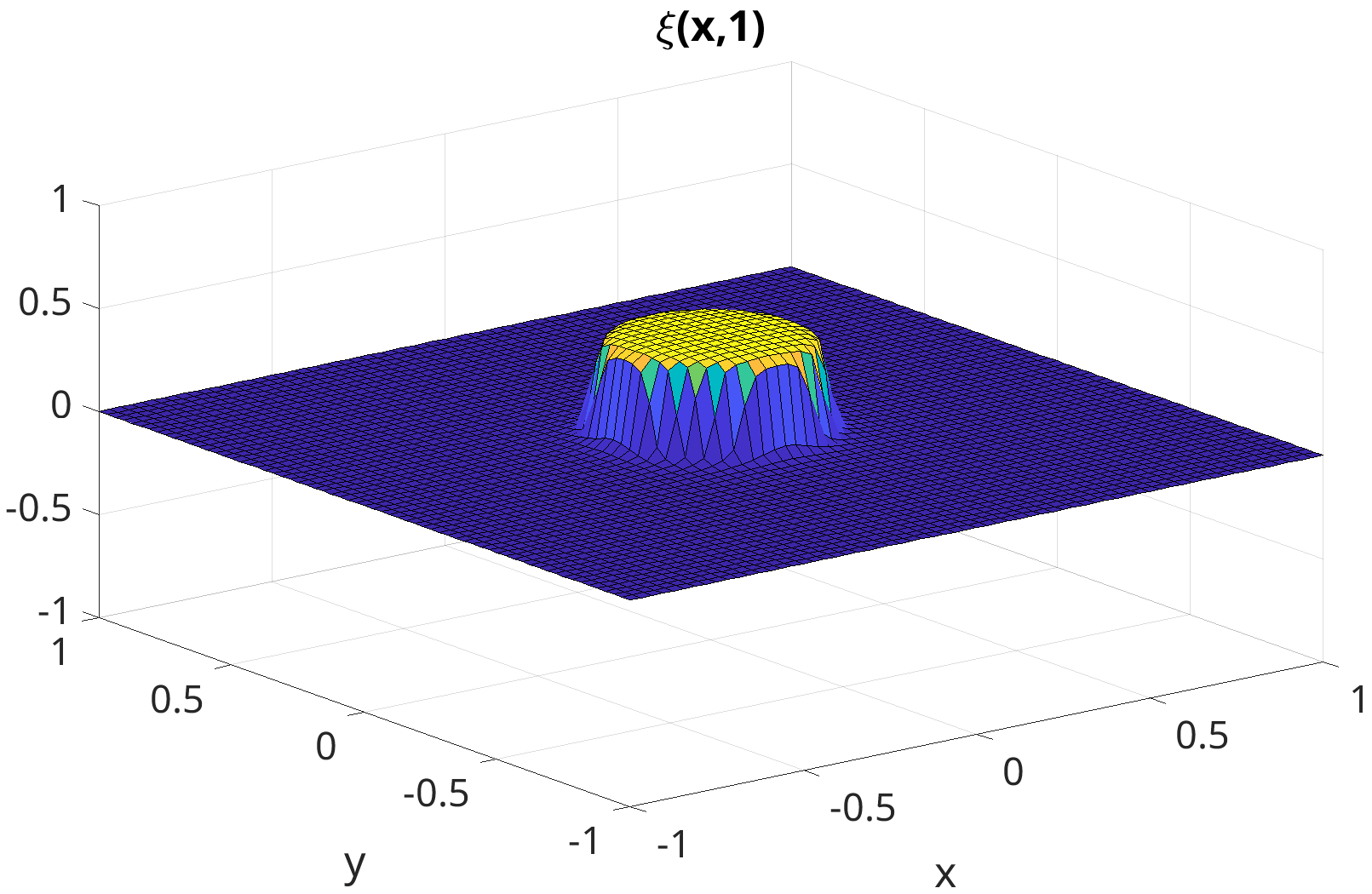}
		\includegraphics[width=4.8cm]{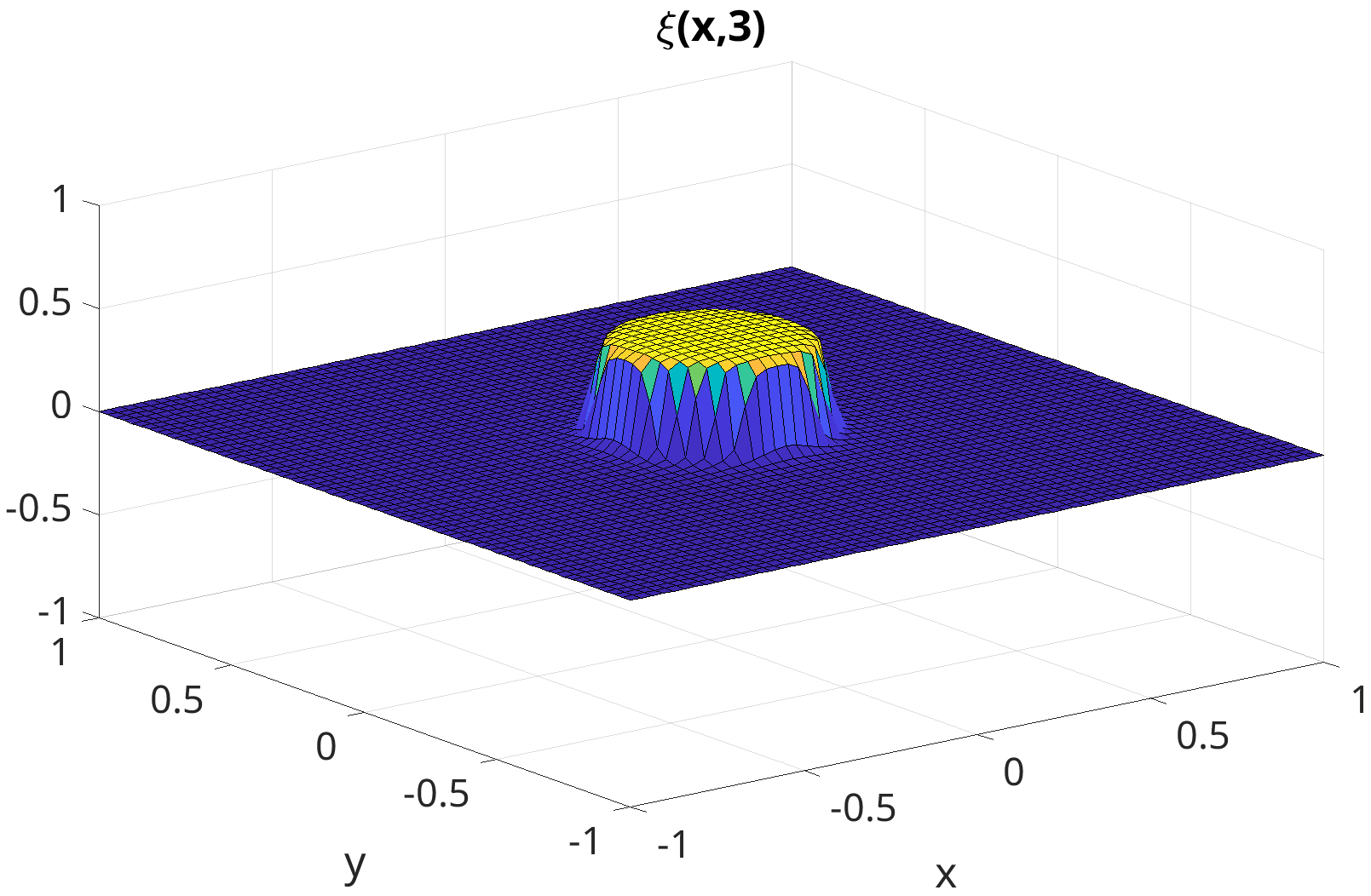}
		\includegraphics[width=4.8cm]{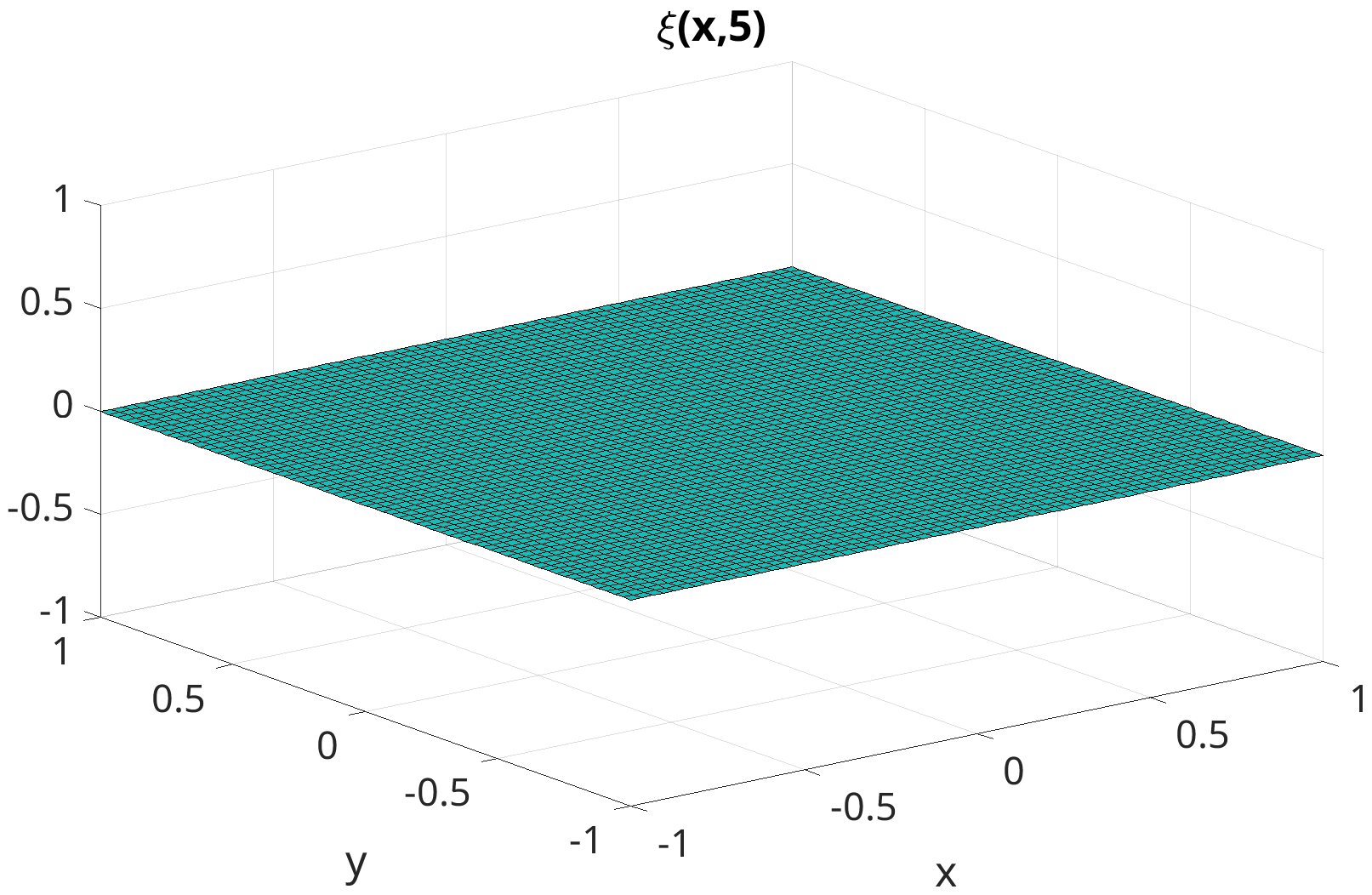}
		\caption{{\bf Case 1}, small control region. Space distributions of the optimal control variables $c(x, t)$ (first row), $\tau (x, t)$ (second row) and $\xi(x, t)$ (third row) at three time instants.}
		\label{fig_reg_contr_small_c_tau_xi}
	\end{center}
\end{figure}

Figure \ref{fig_reg_contr_small_c_tau_xi}. The level of the optimal consumption per unit of capital, $c*,$ becomes maximal when time is close to the horizon $ T.$ The optimal environmental taxation rate and the optimal pollution abatement approach $ 0$ as time approaches the time horizon $T.$ Note that $c^* + \tau^* = 1 - s$ for any $x \in \omega$ and $t \in (0, T).$

The same remarks are true for the other  cases.

Moreover, in Case 2 (the large control region) we get a higher level of capital and a lower level of pollution.

\begin{figure}[!htb]
	\begin{center}
		\includegraphics[width=4.8cm]{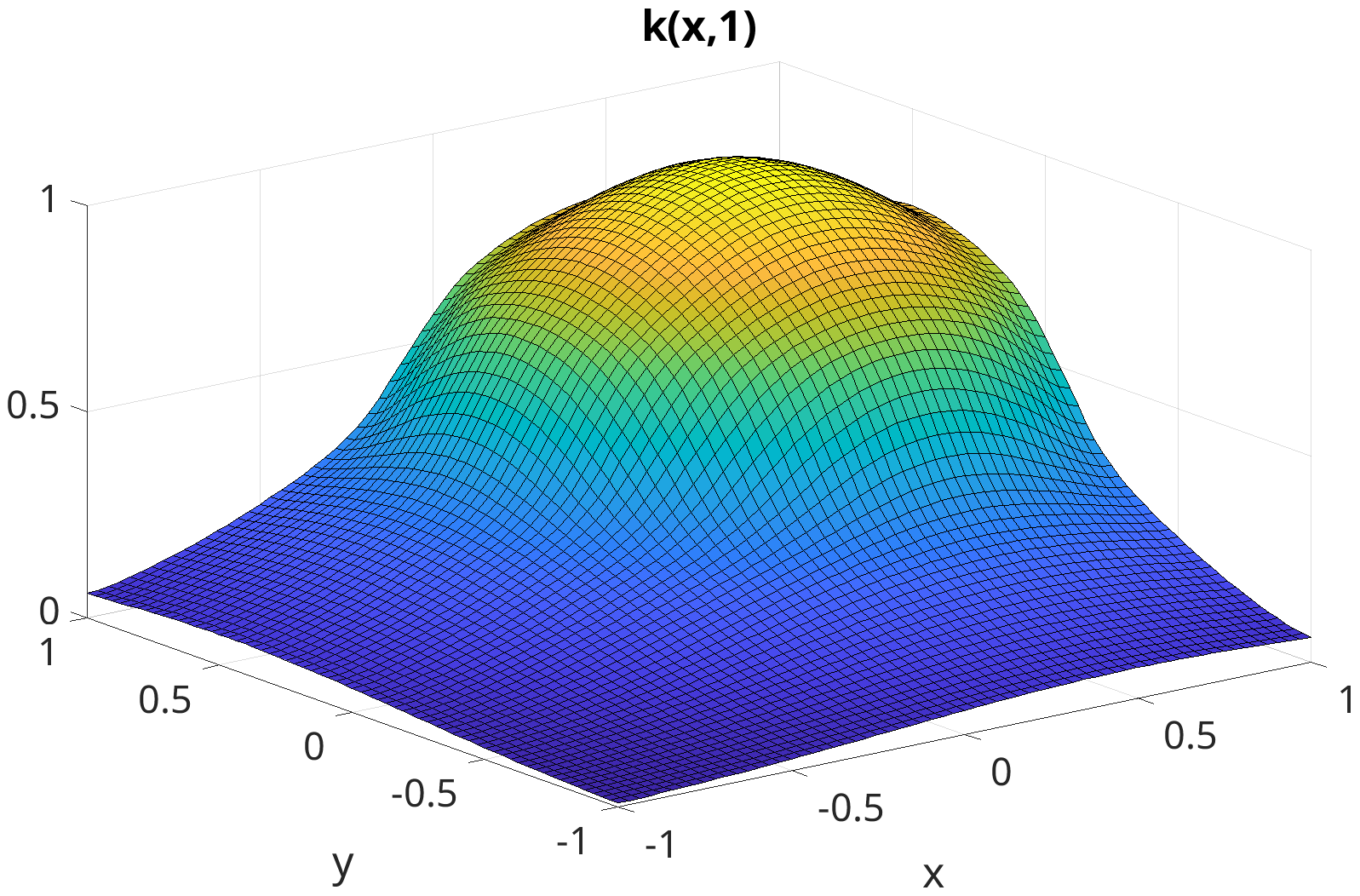}
		\includegraphics[width=4.8cm]{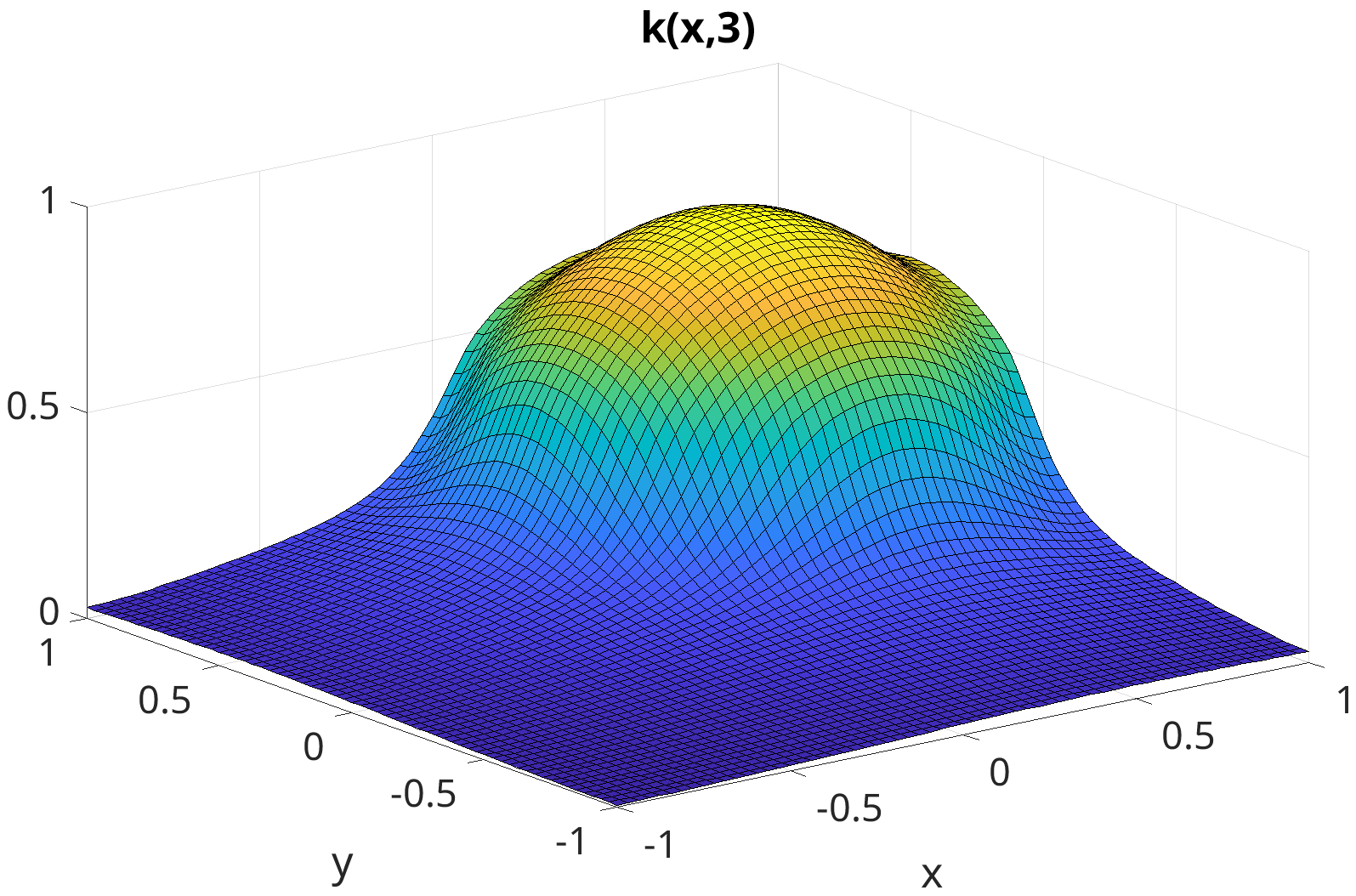}
		\includegraphics[width=4.8cm]{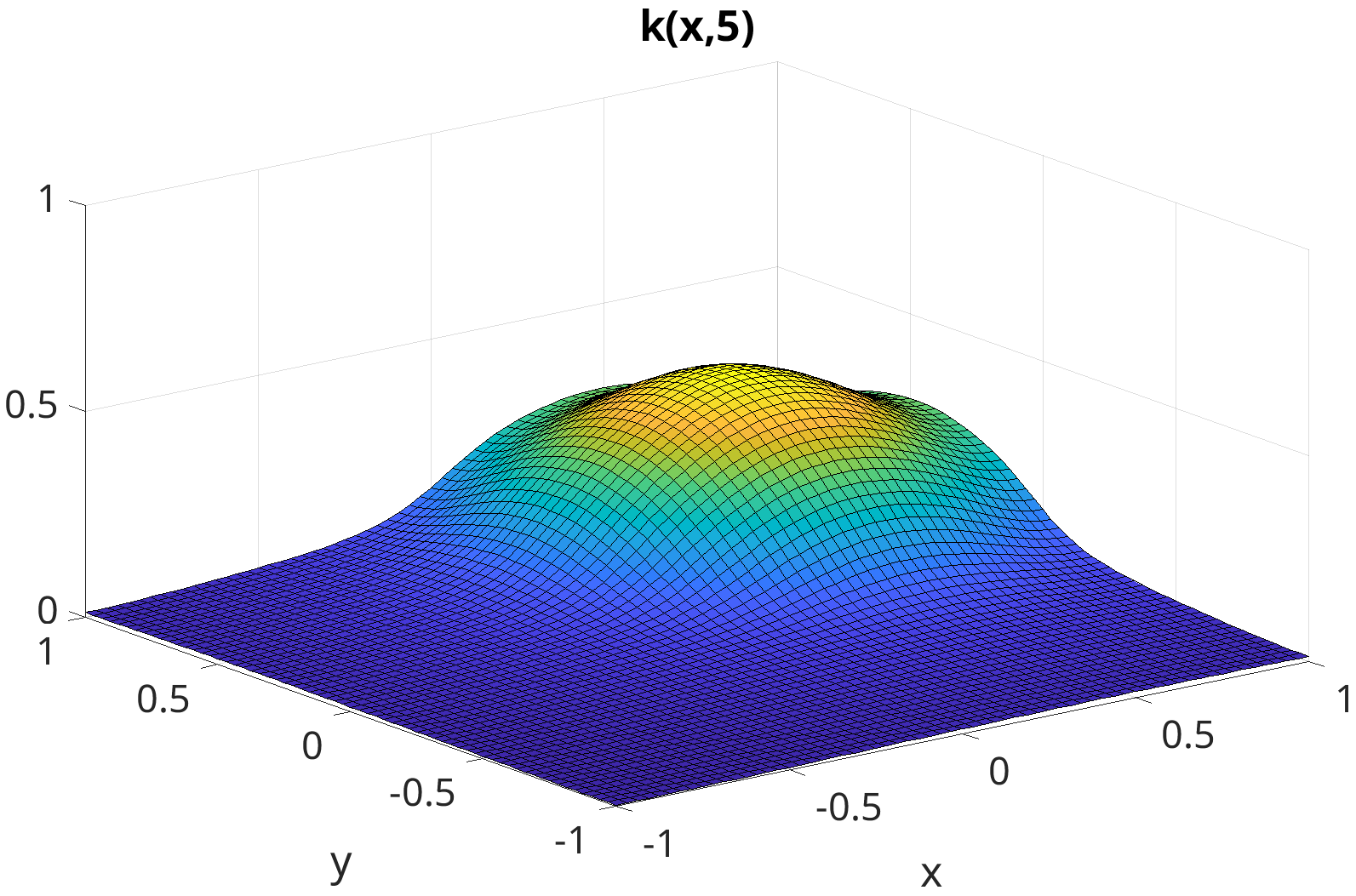}\\
		\includegraphics[width=4.8cm]{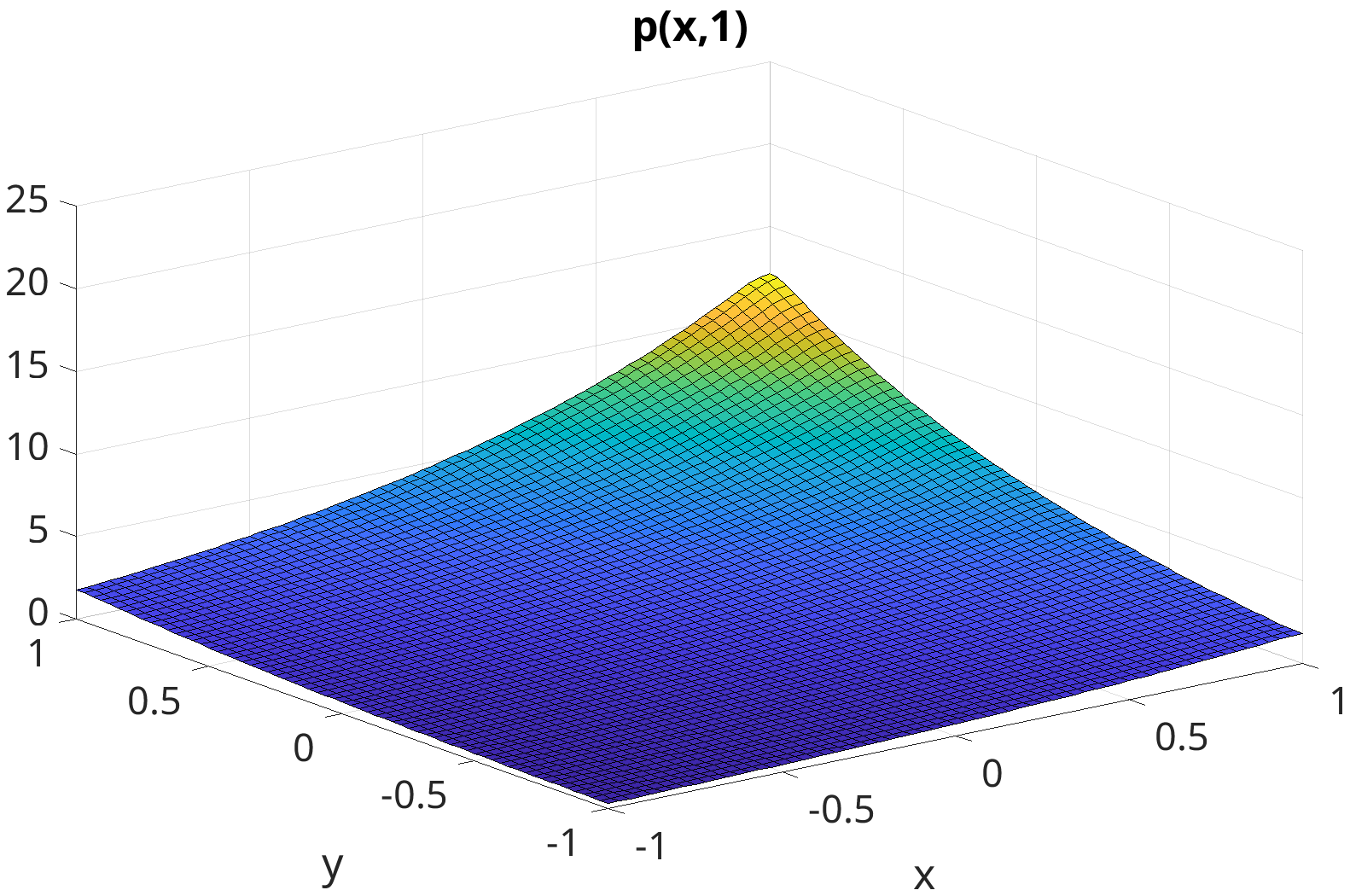}
		\includegraphics[width=4.8cm]{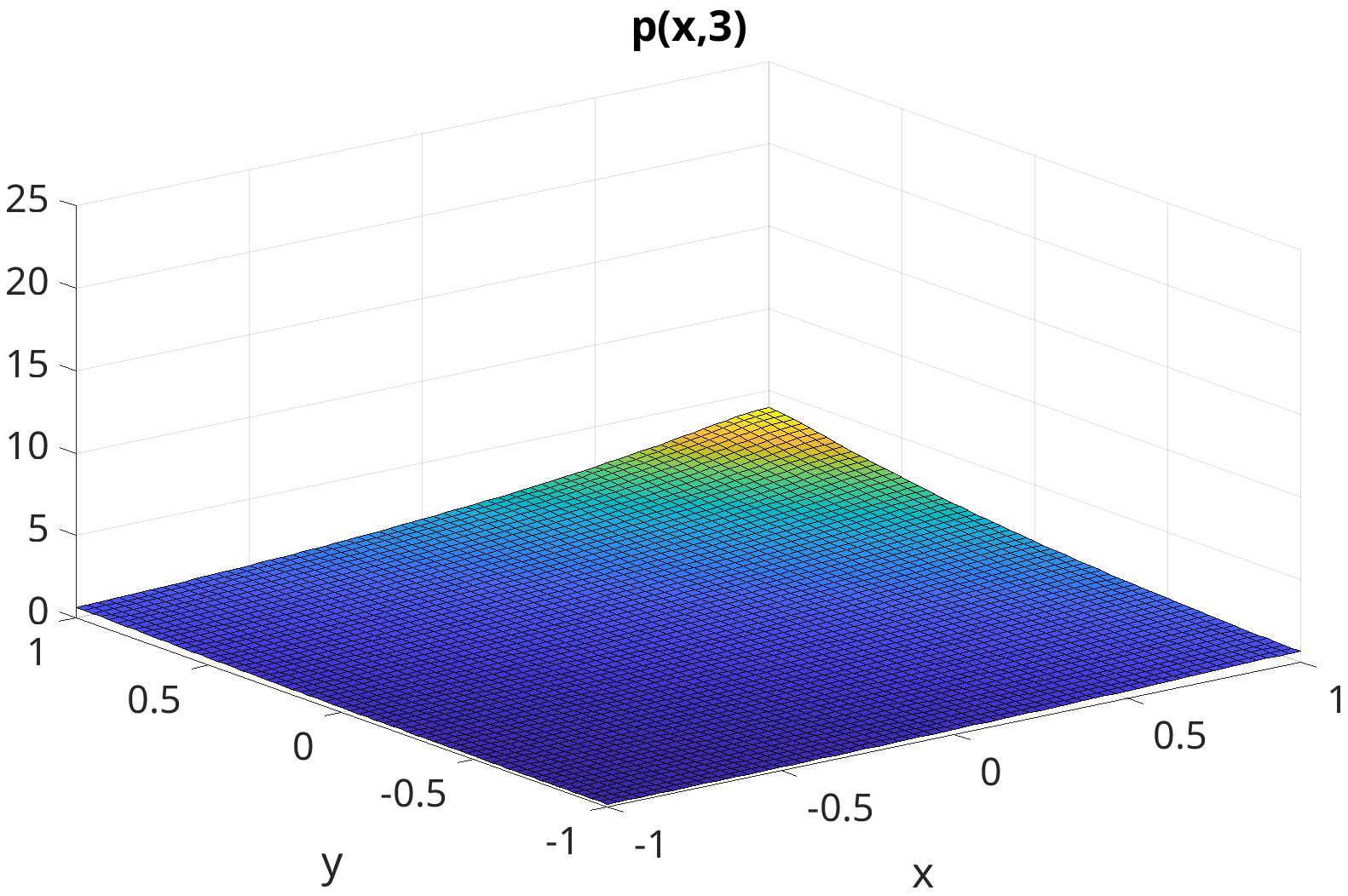}
		\includegraphics[width=4.8cm]{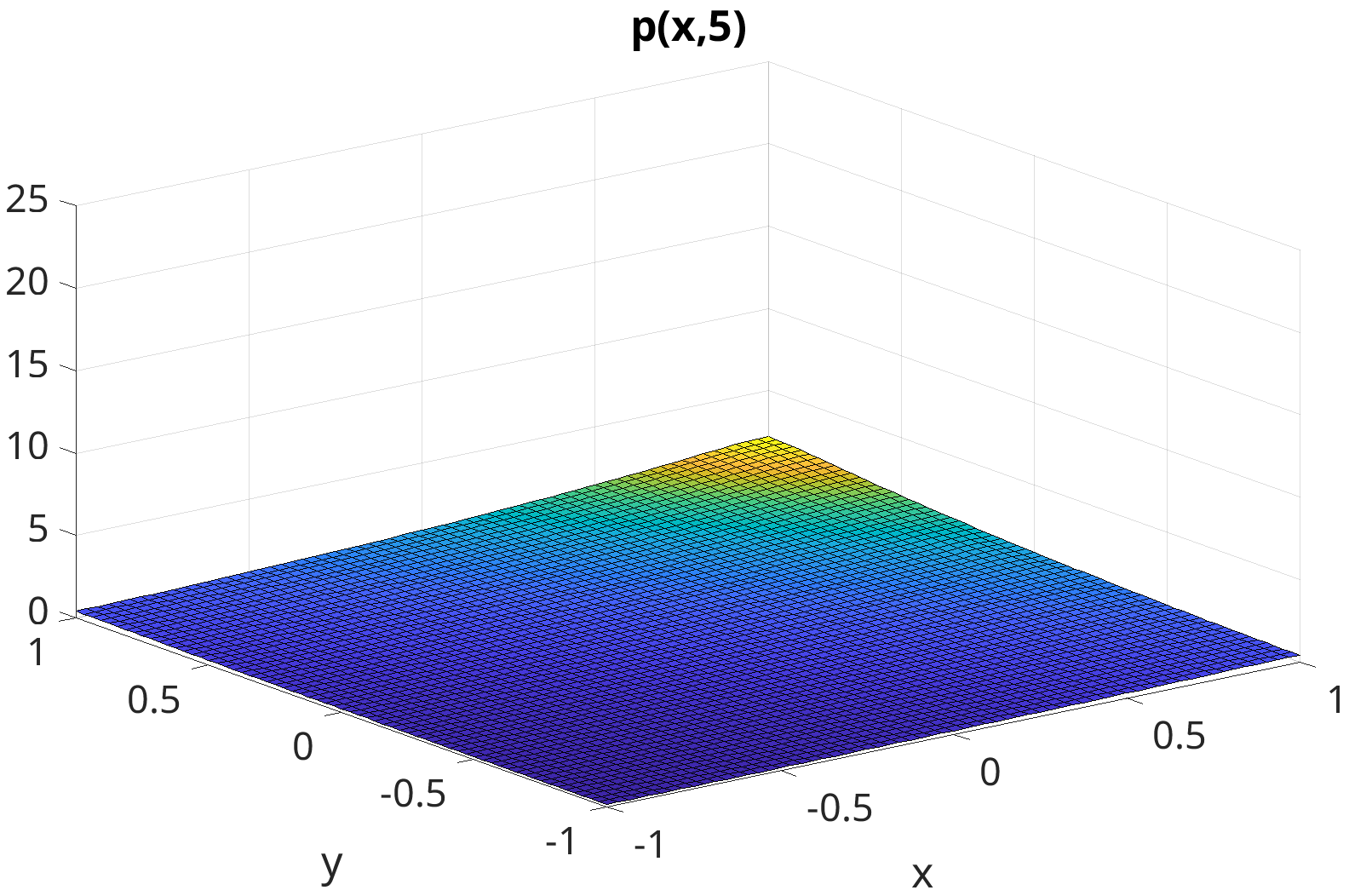}
		\caption{{\bf Case 3}, the control region is the whole domain. Space distributions of the optimal state variables $k(x,t)$ (first row) and $p(x,t)$ (second row) at three time instants.}
		\label{fig_reg_tot_contr_k_p}
	\end{center}
\end{figure}

\begin{figure}[!htb]
	\begin{center}
		\includegraphics[width=4.8cm]{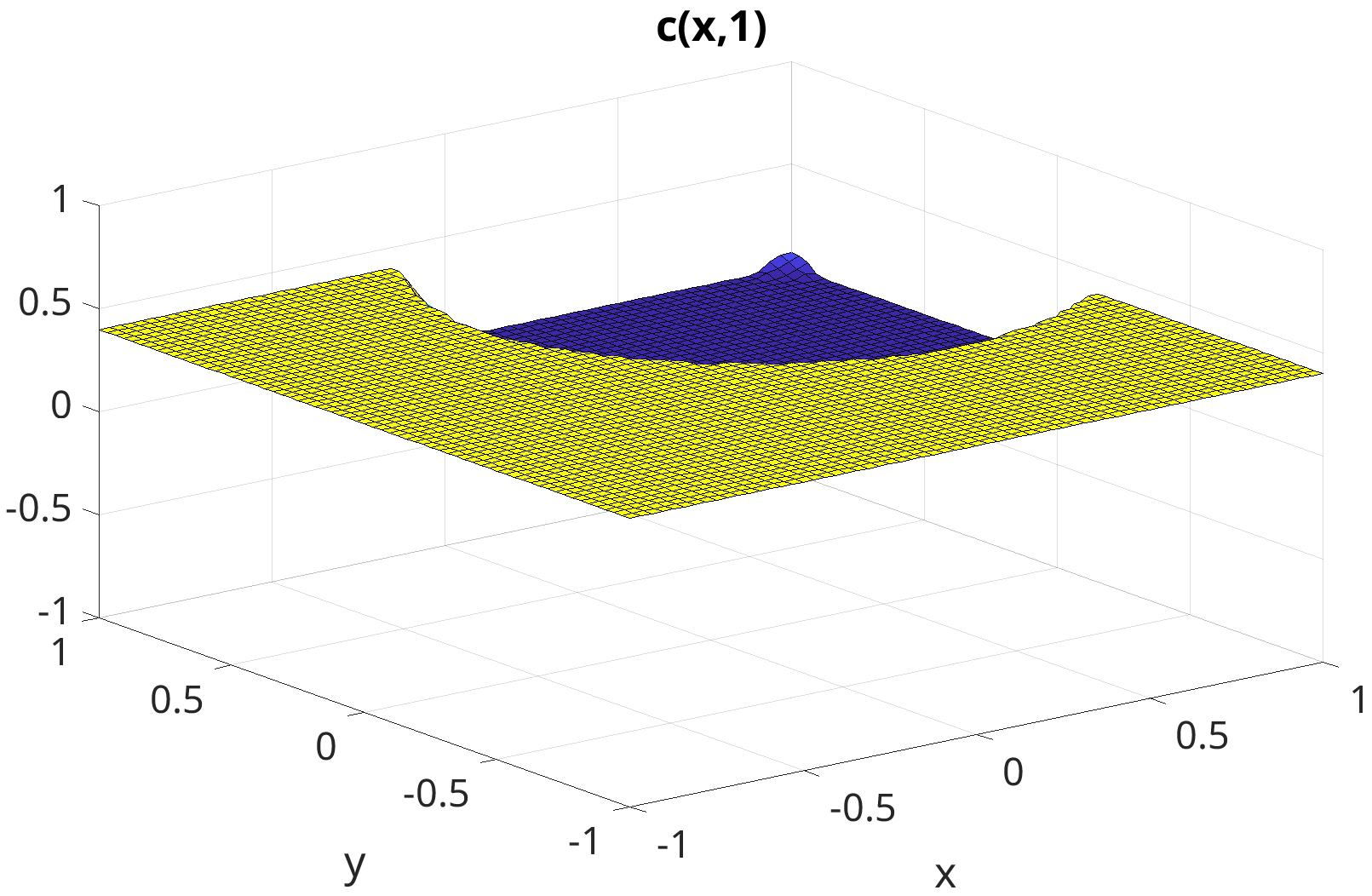}
		\includegraphics[width=4.8cm]{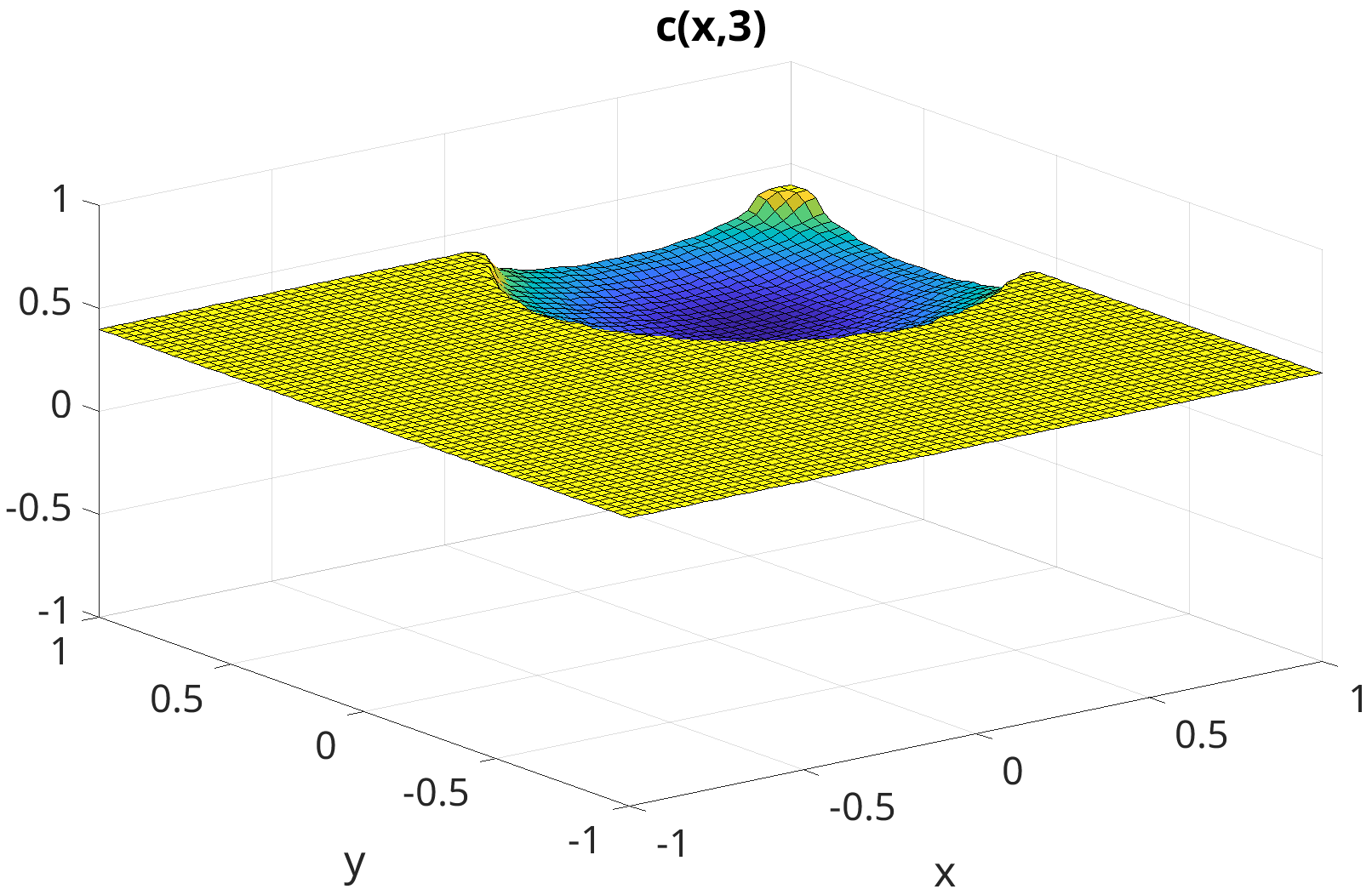}
		\includegraphics[width=4.8cm]{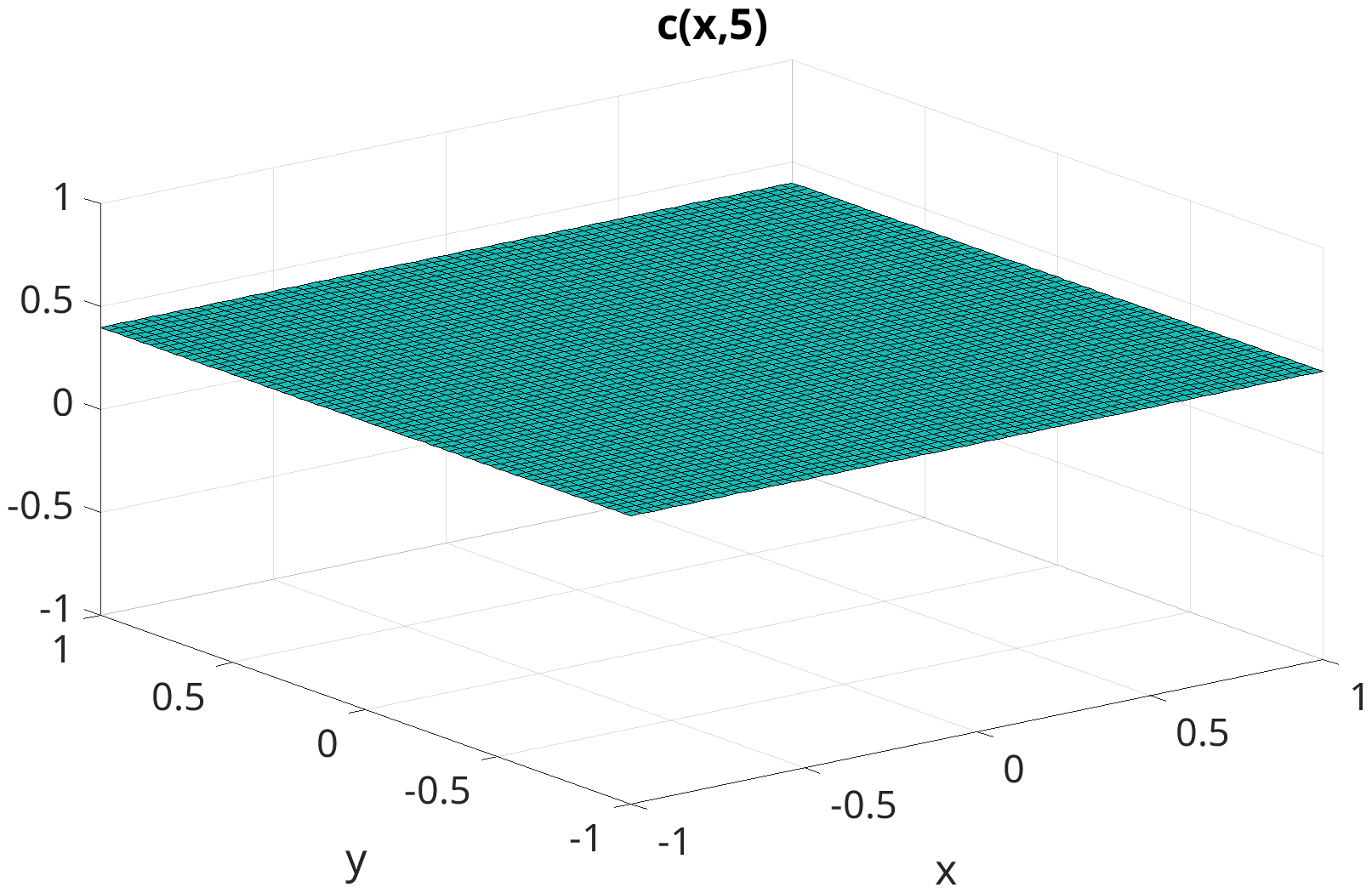}\\
		\includegraphics[width=4.8cm]{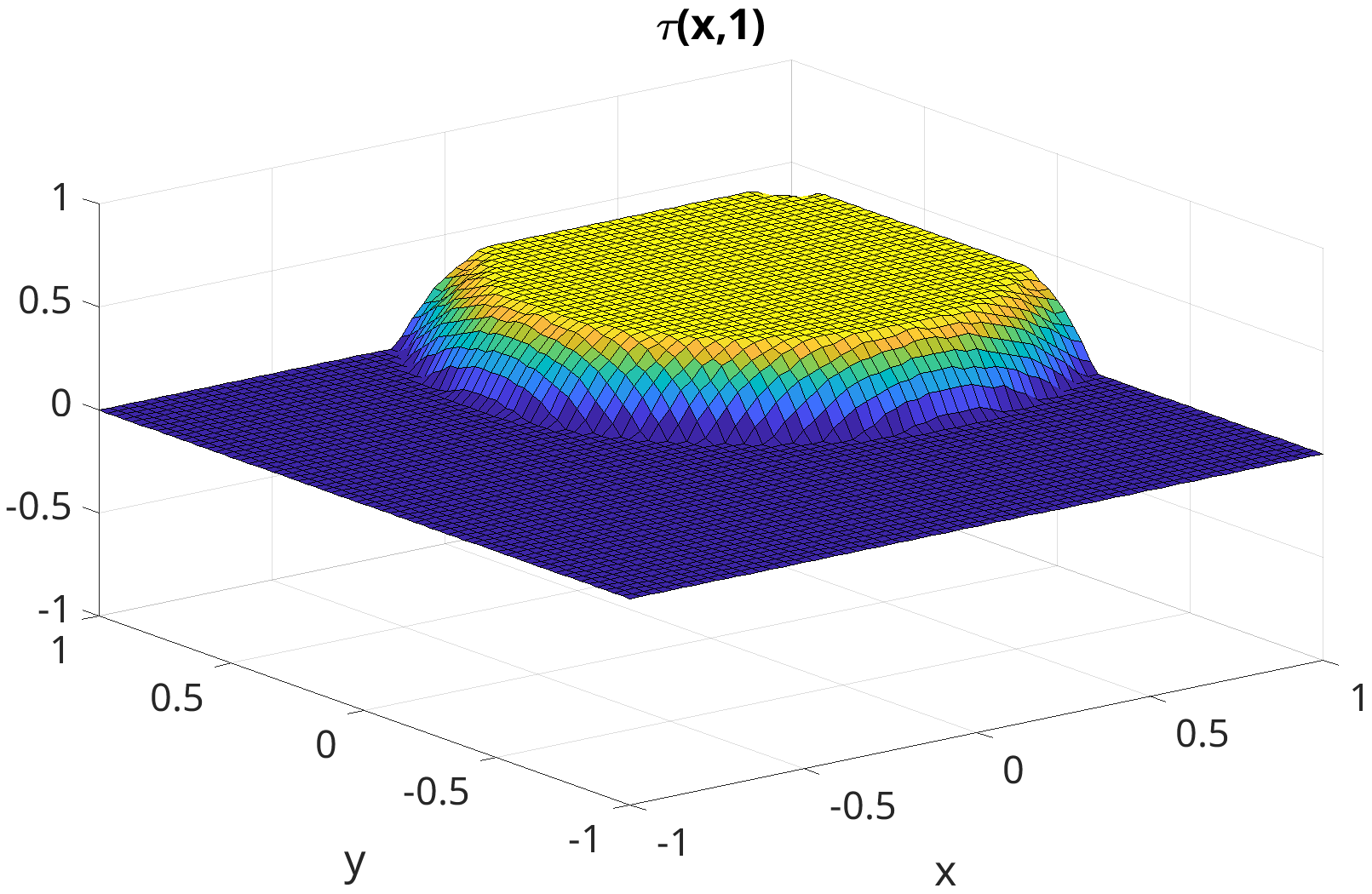}
		\includegraphics[width=4.8cm]{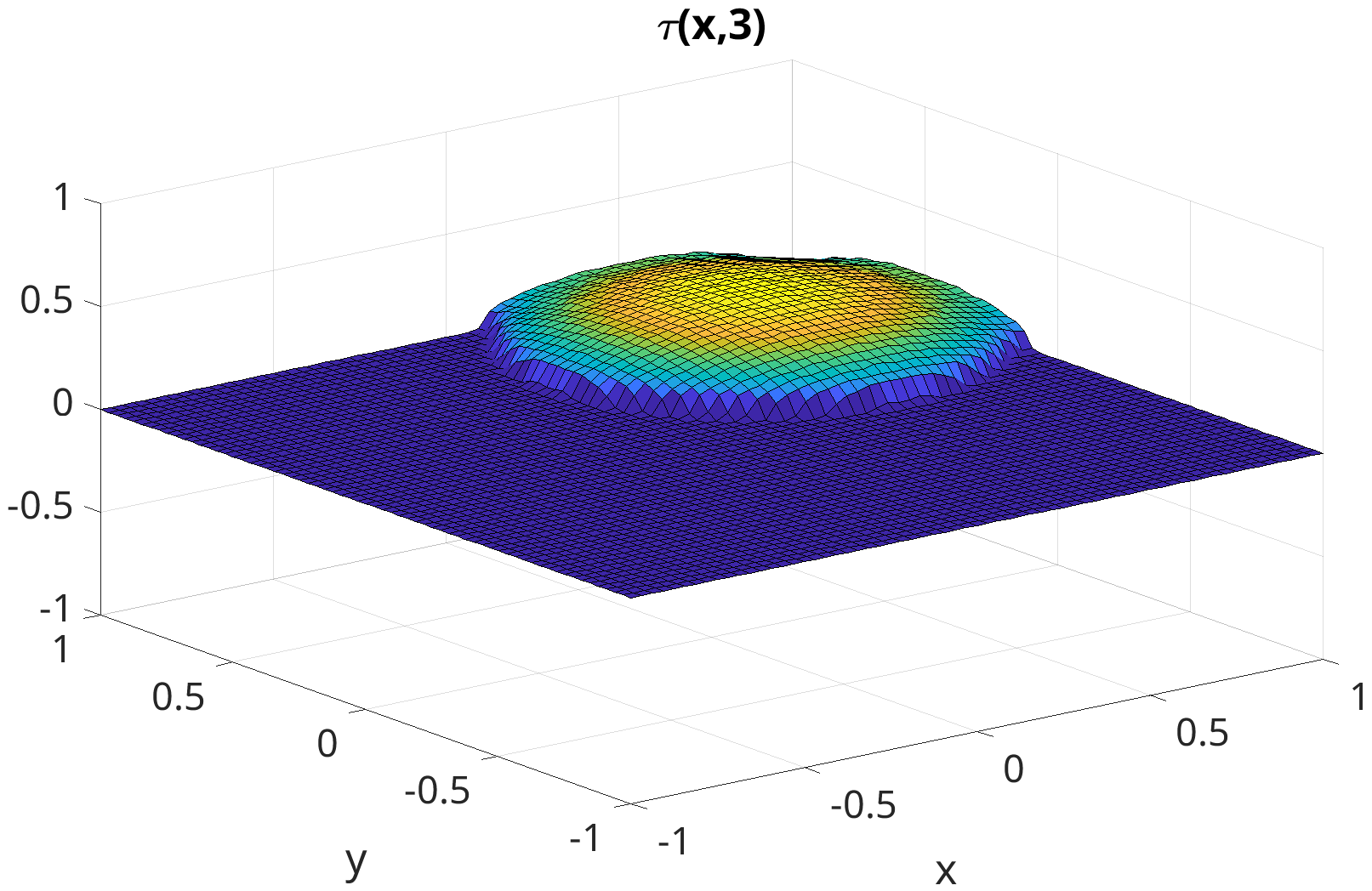}
		\includegraphics[width=4.8cm]{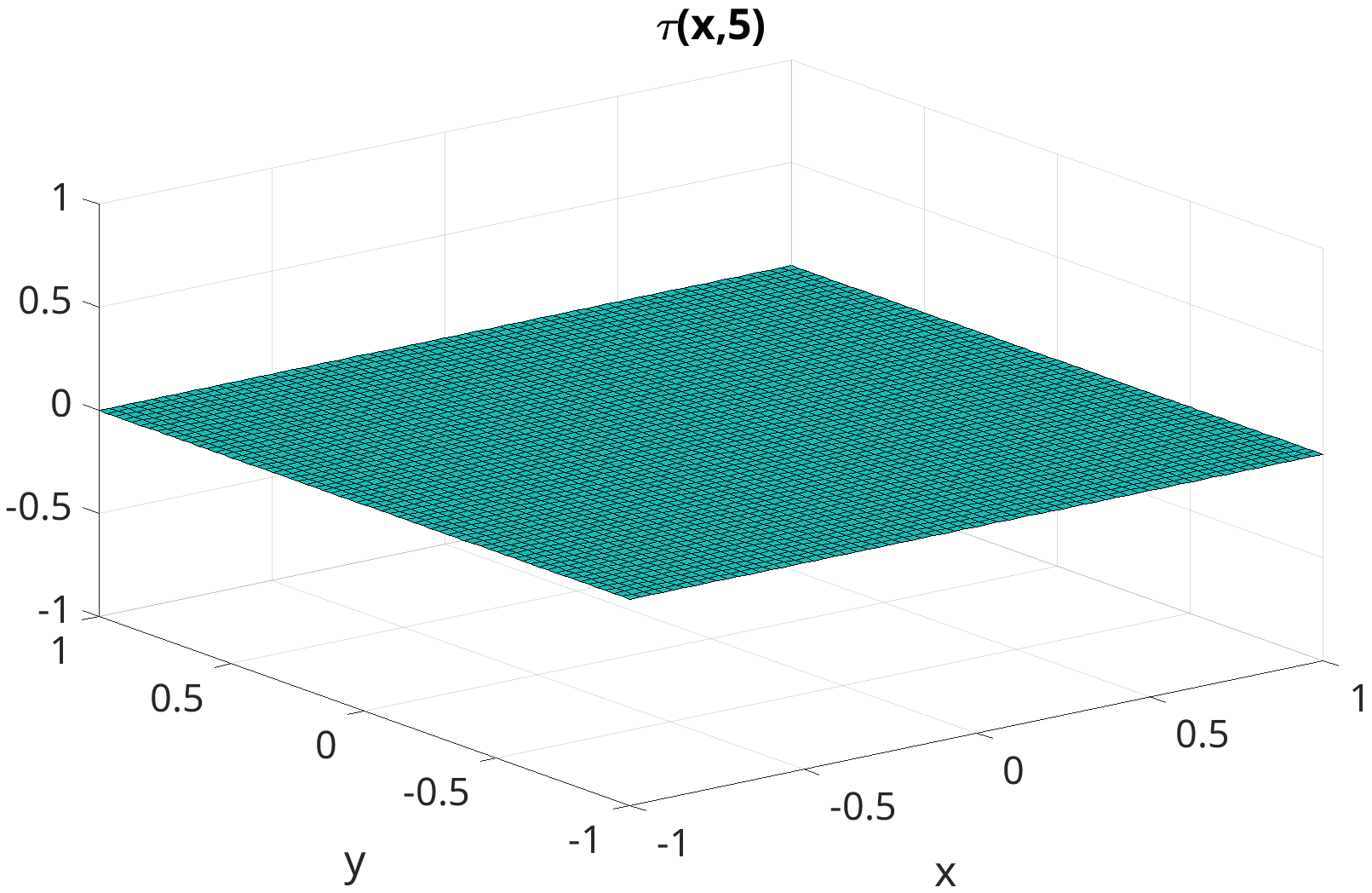}
		\includegraphics[width=4.8cm]{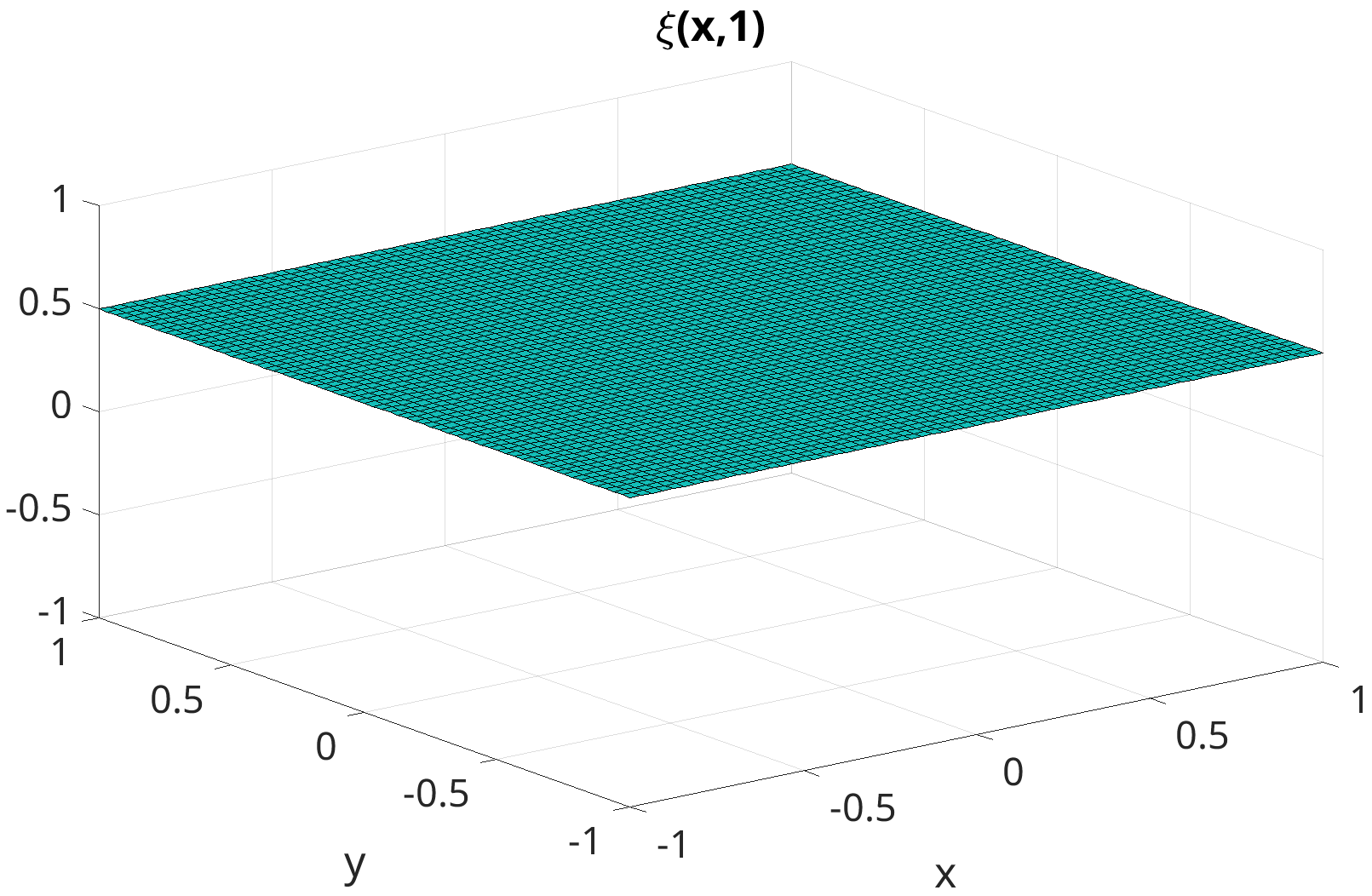}
		\includegraphics[width=4.8cm]{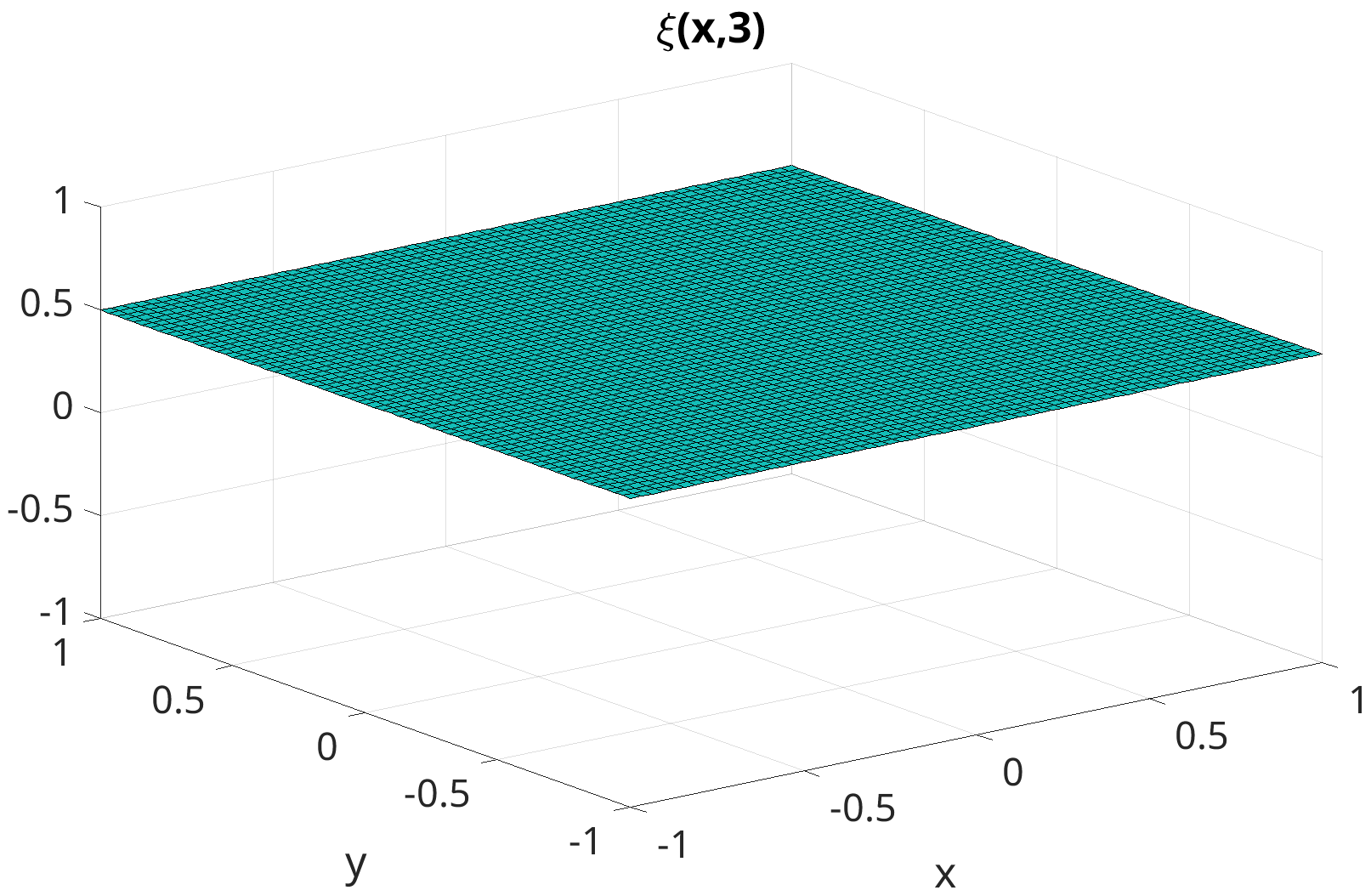}
		\includegraphics[width=4.8cm]{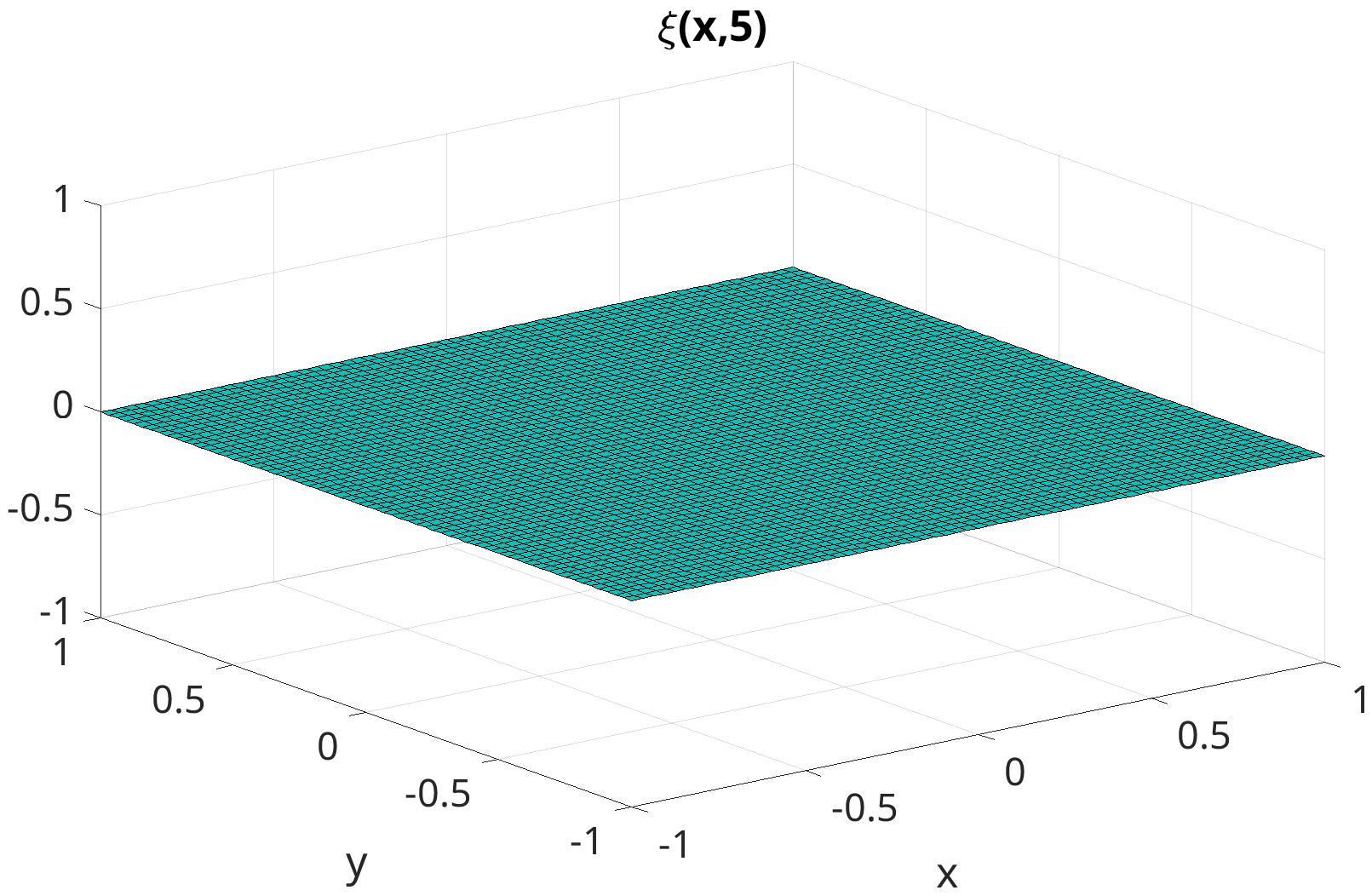}
		\caption{{\bf Case 3}, the control region is the whole domain. Space distributions of the optimalcontrol variables $c(x,t)$ (first row), $\tau(x,t)$ (second row) and $\xi(x,t)$ (third row) at three time instants.}
		\label{fig_reg_tot_contr_c_tau_xi}
	\end{center}
\end{figure}

Figure \ref{fig_reg_tot_contr_k_p} (corresponding to Case 3) shows an important diminishment of $k(x,t)$ as $t\rightarrow T$. On the other hand, Figure 	\ref{fig_reg_tot_contr_c_tau_xi} shows that in this case
$\tau ^*(x,t)$ is big for $x$ close to the region with initial high level of pollution if $t$ is small and it is small for large $t$ or for $x$ far from this region.

%==========================

%\newpage

\begin{figure}[!htb]
	\begin{center}
		\includegraphics[width=4.8cm]{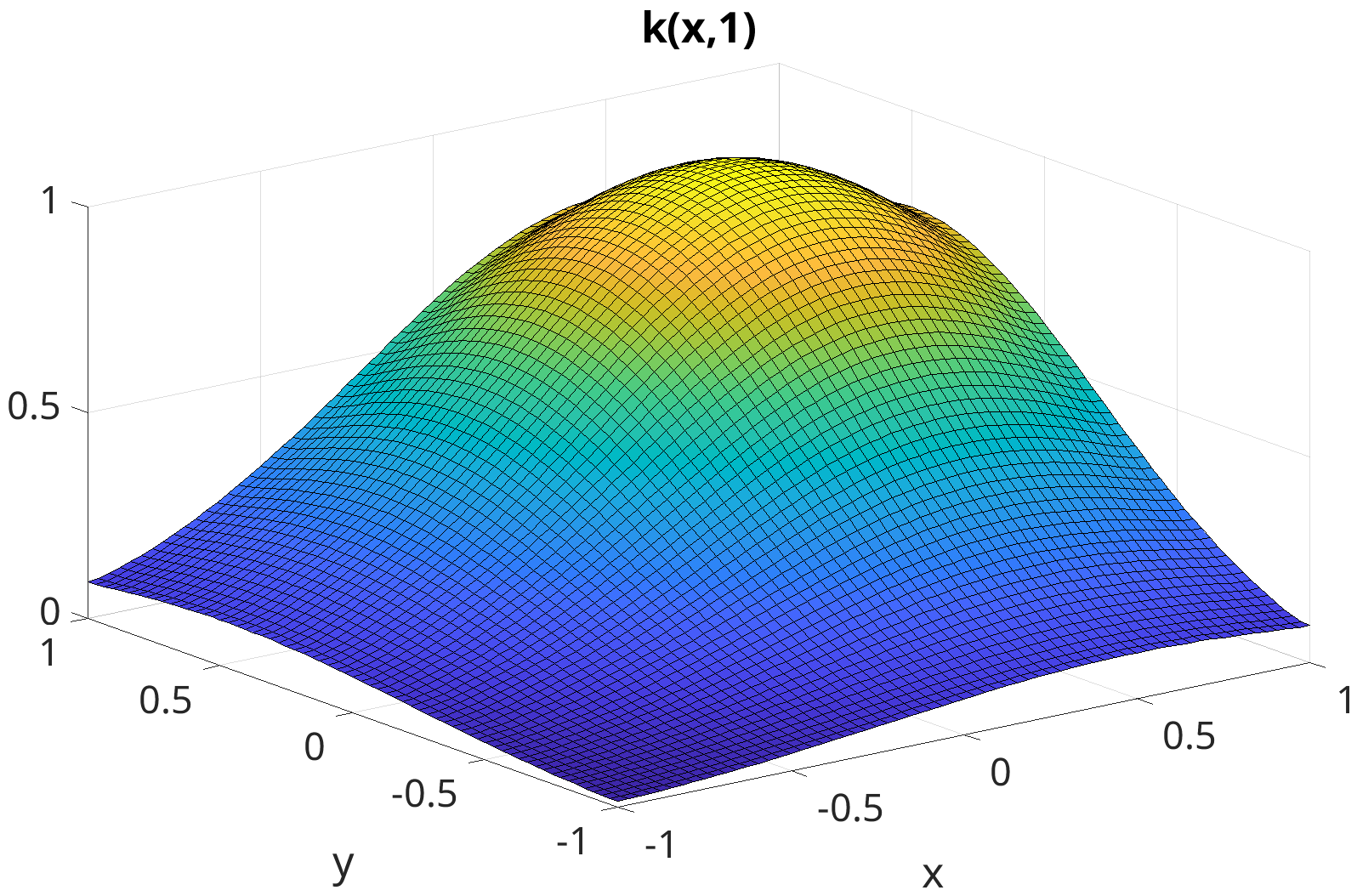}
		\includegraphics[width=4.8cm]{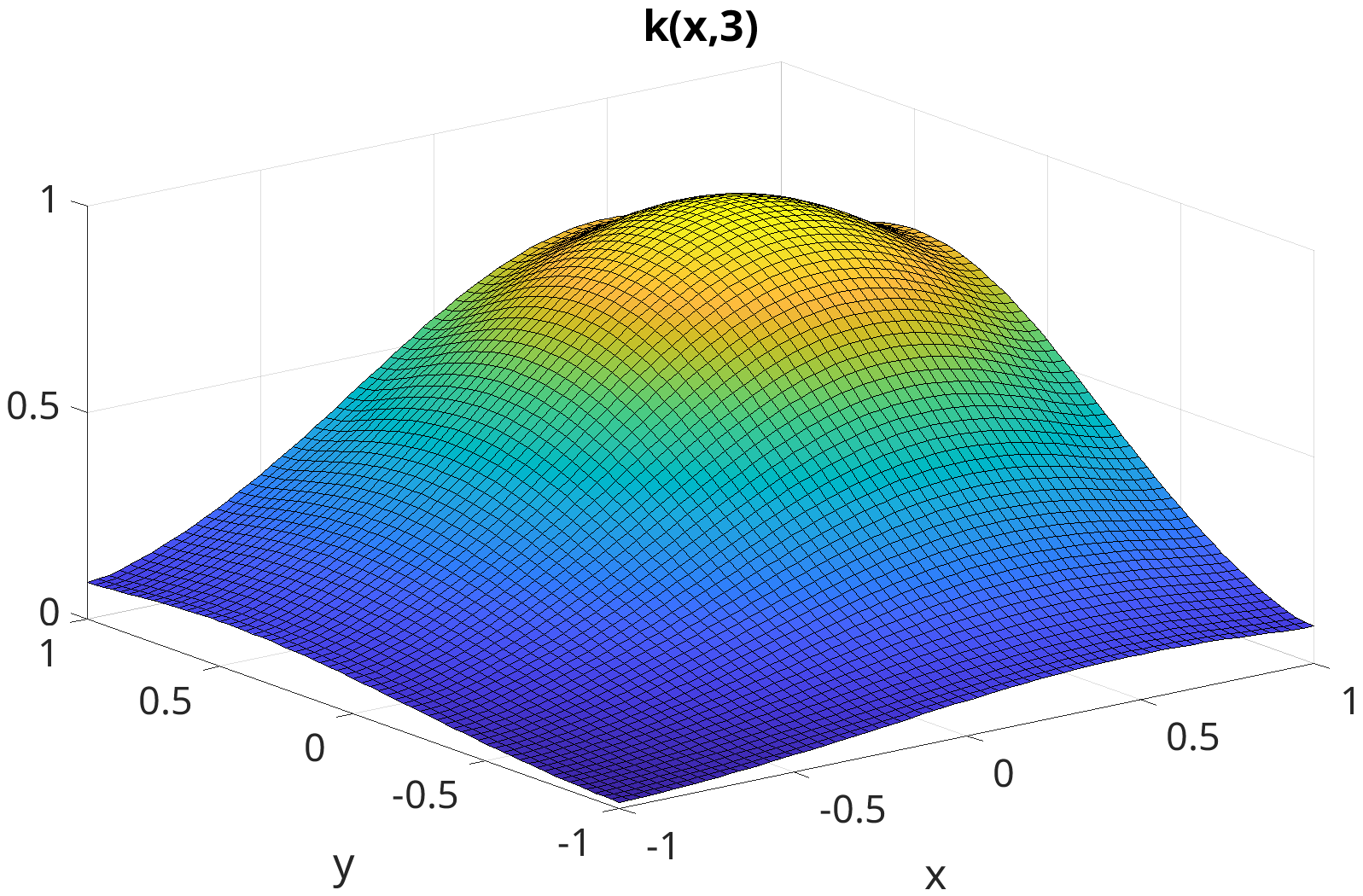}
		\includegraphics[width=4.8cm]{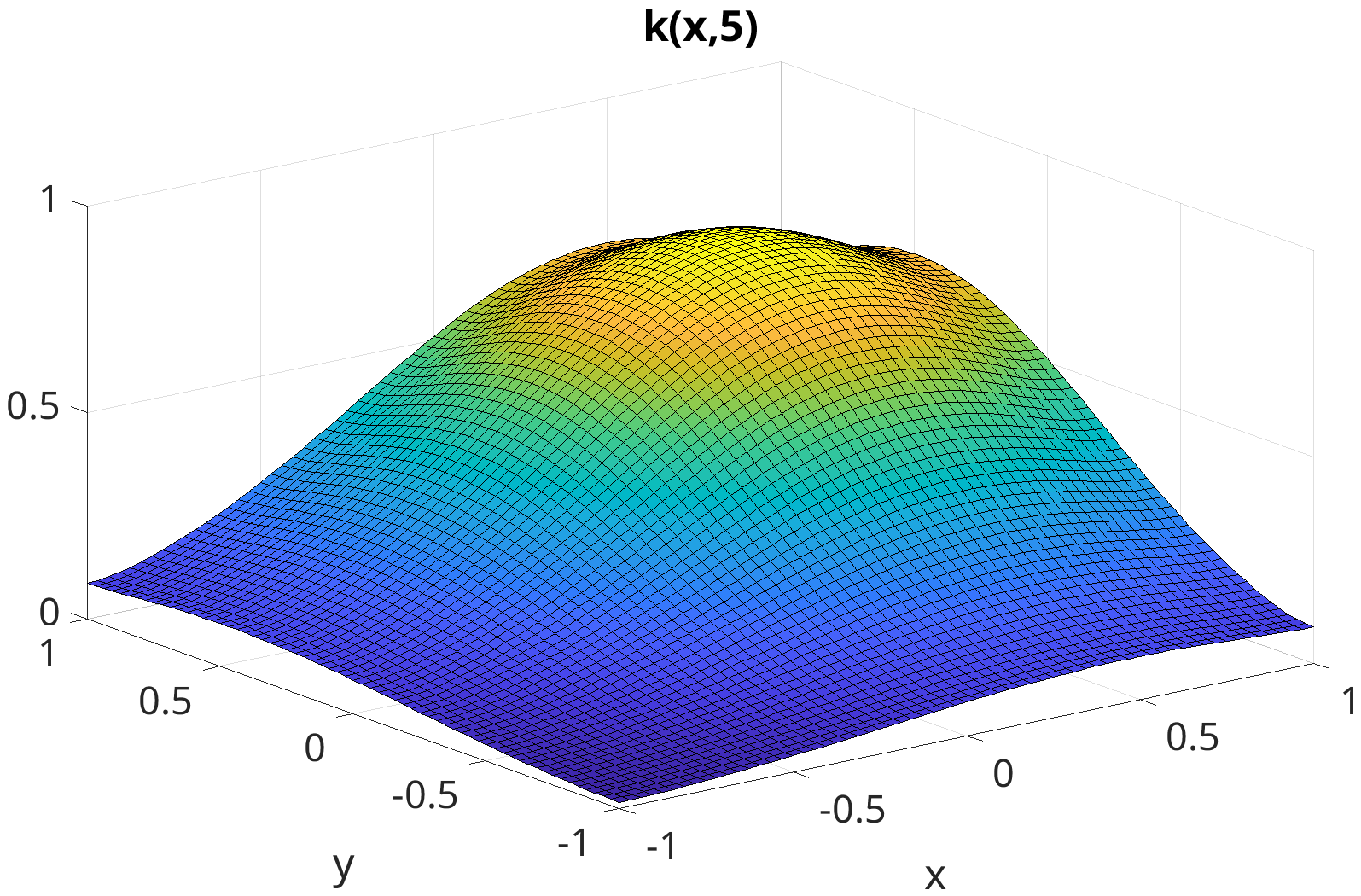}\\
		\includegraphics[width=4.8cm]{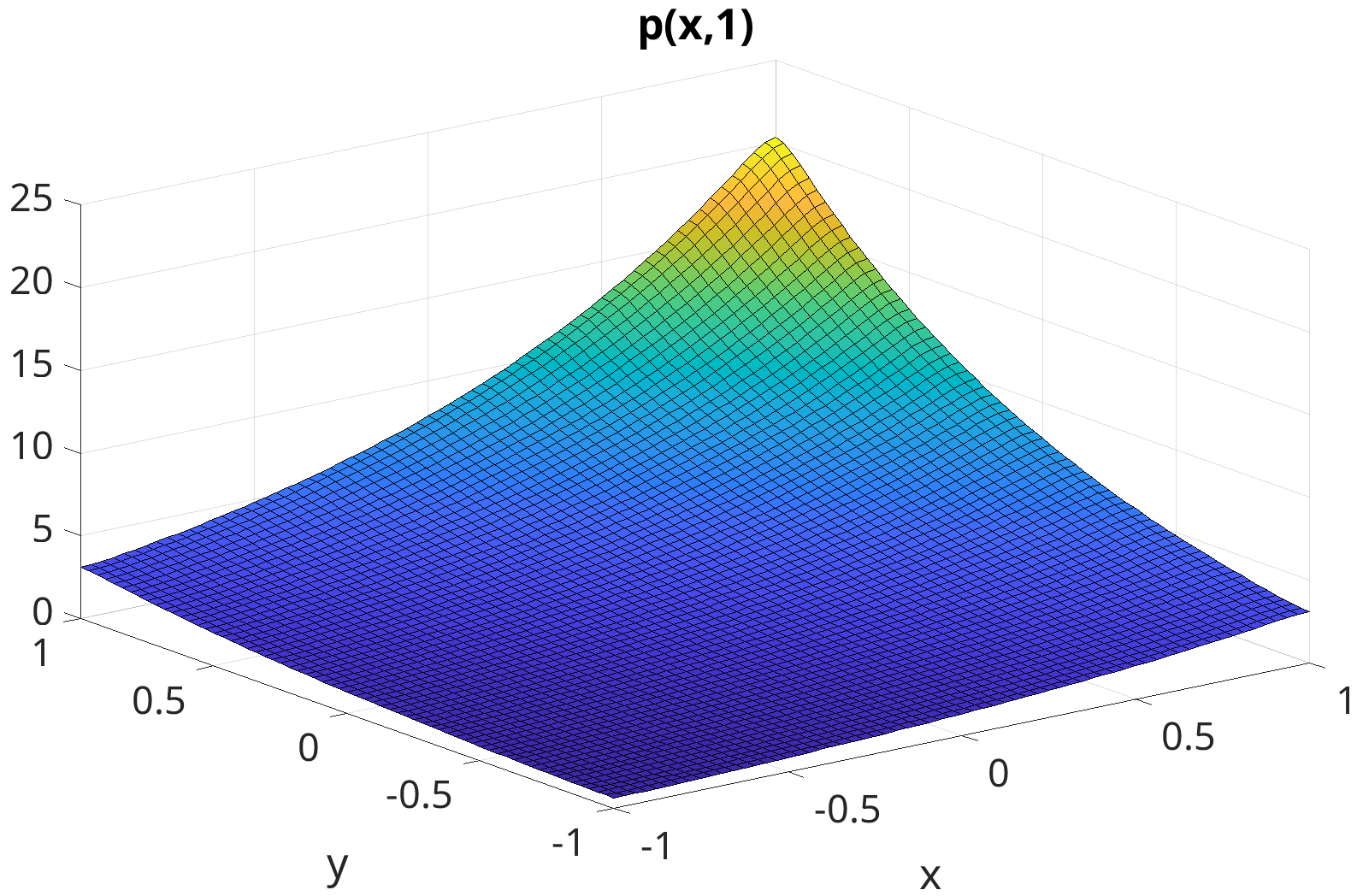}
		\includegraphics[width=4.8cm]{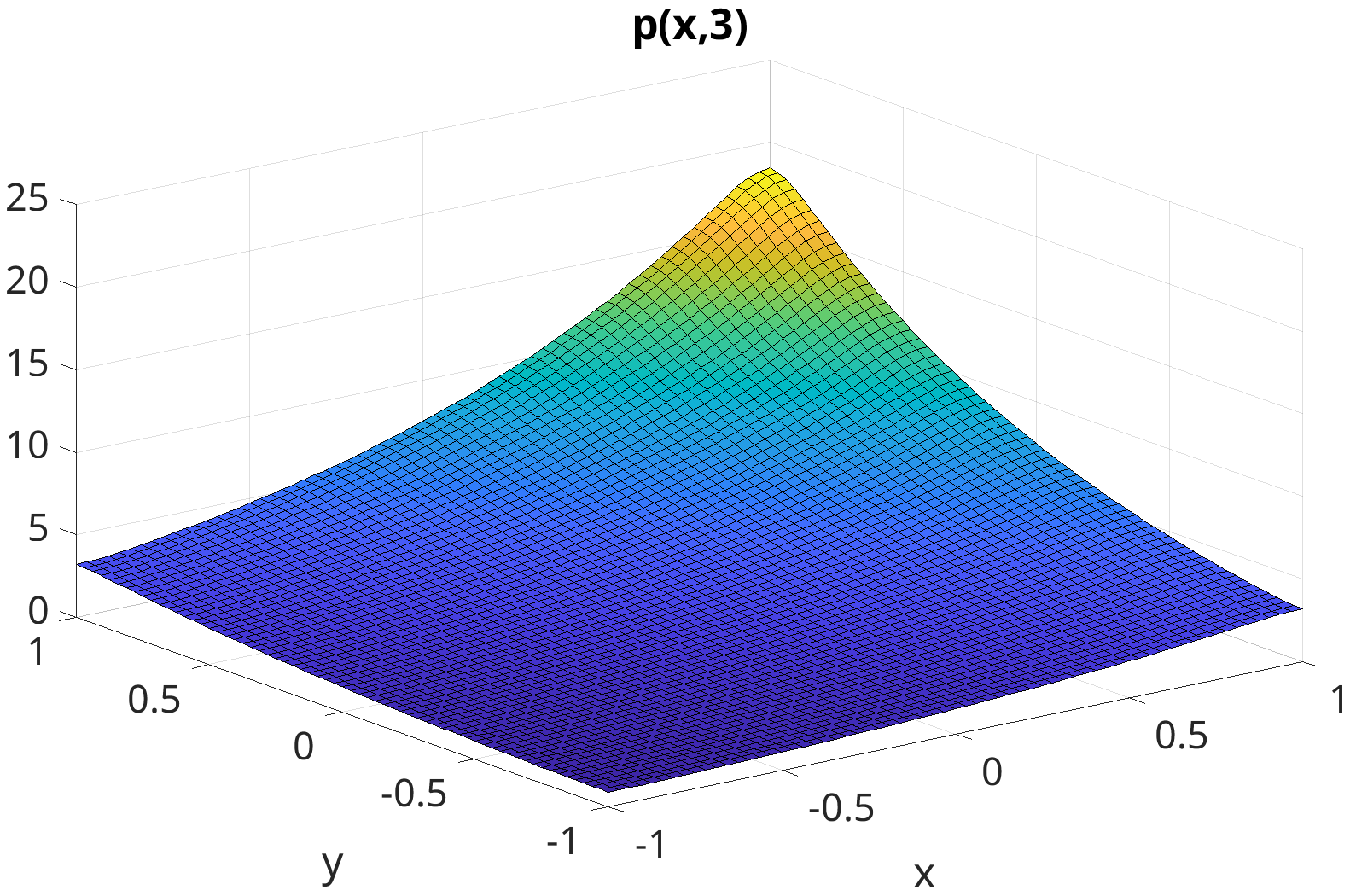}
		\includegraphics[width=4.8cm]{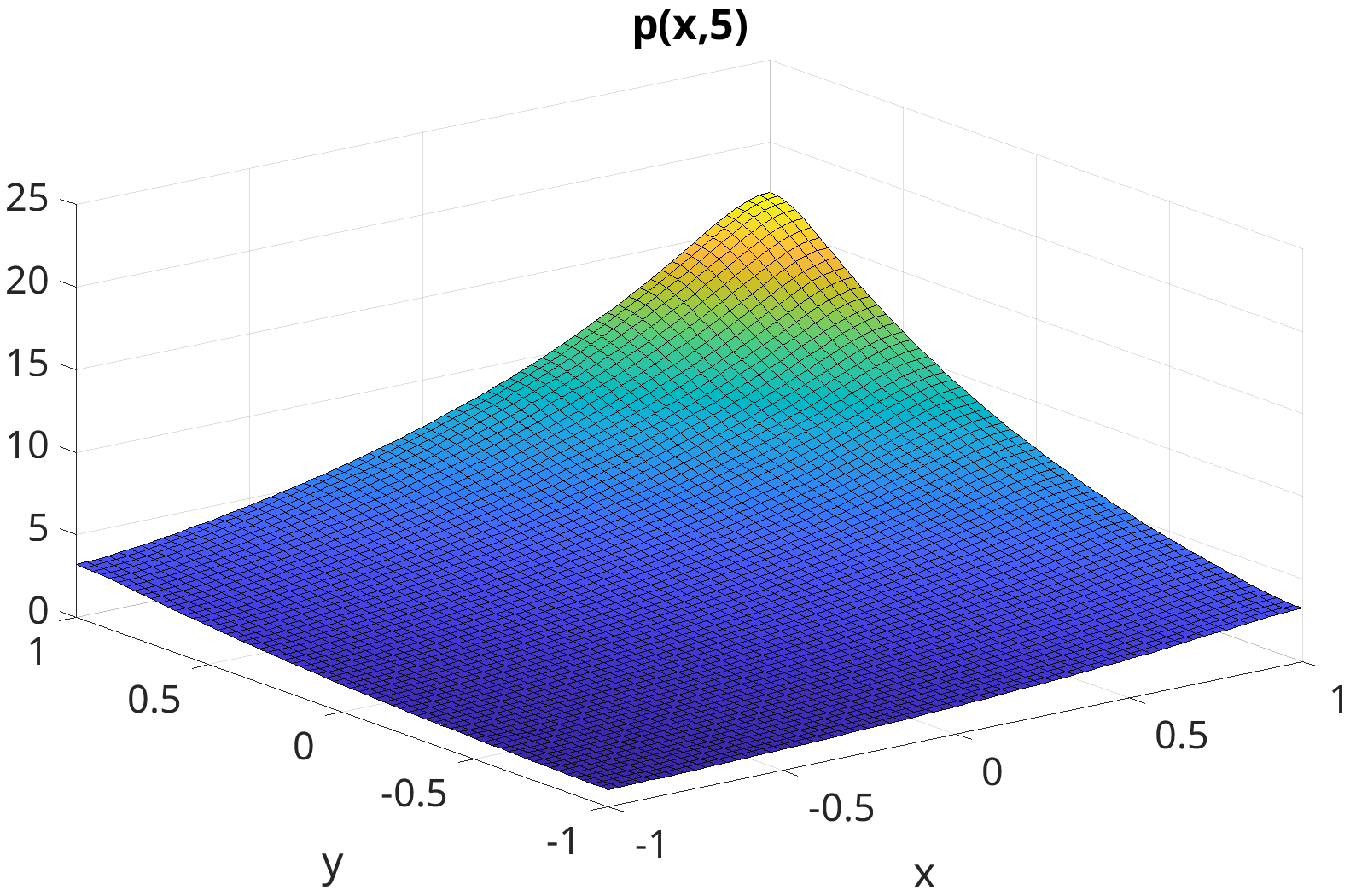}
		\caption{{\bf Case 0}, no control. Space distributions of state variables $k(x,t)$ (first row) and $p(x,t)$ (second row) at three time instants.}
		\label{fig_no_contr_k_p}
	\end{center}
\end{figure}

\newpage

\begin{figure}[!htb]
	\begin{center}
		\includegraphics[width=4.8cm]{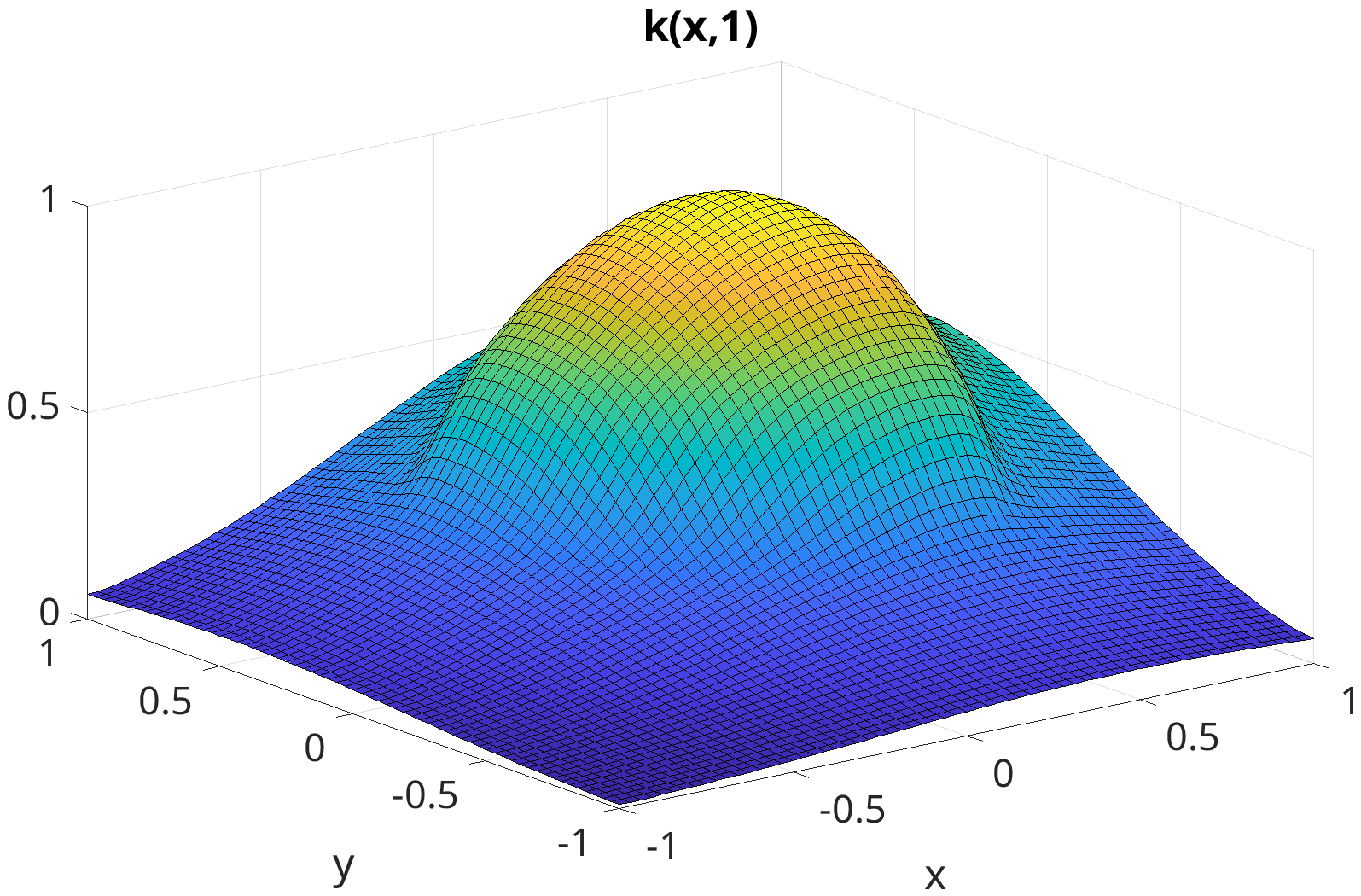}
		\includegraphics[width=4.8cm]{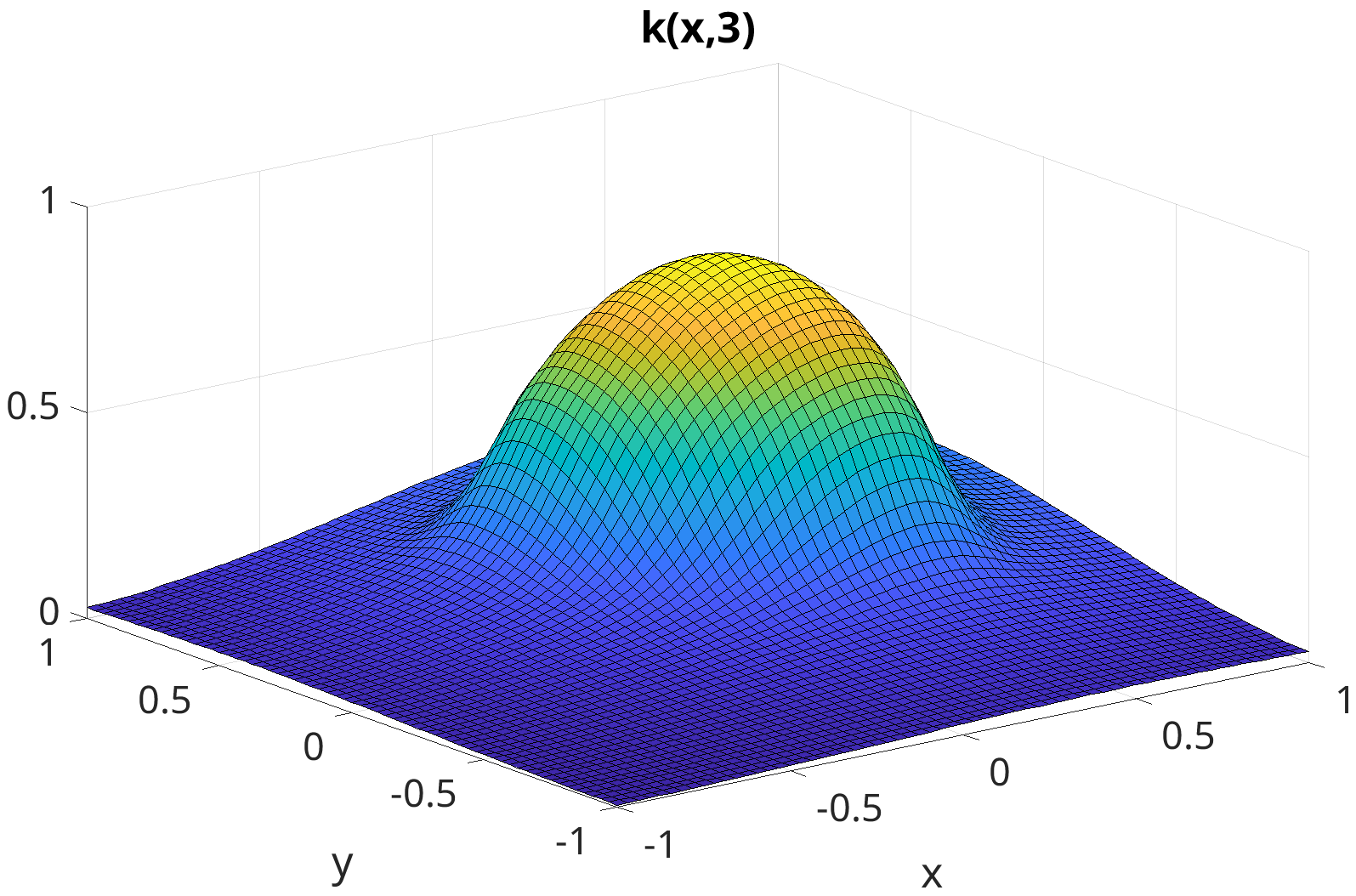}
		\includegraphics[width=4.8cm]{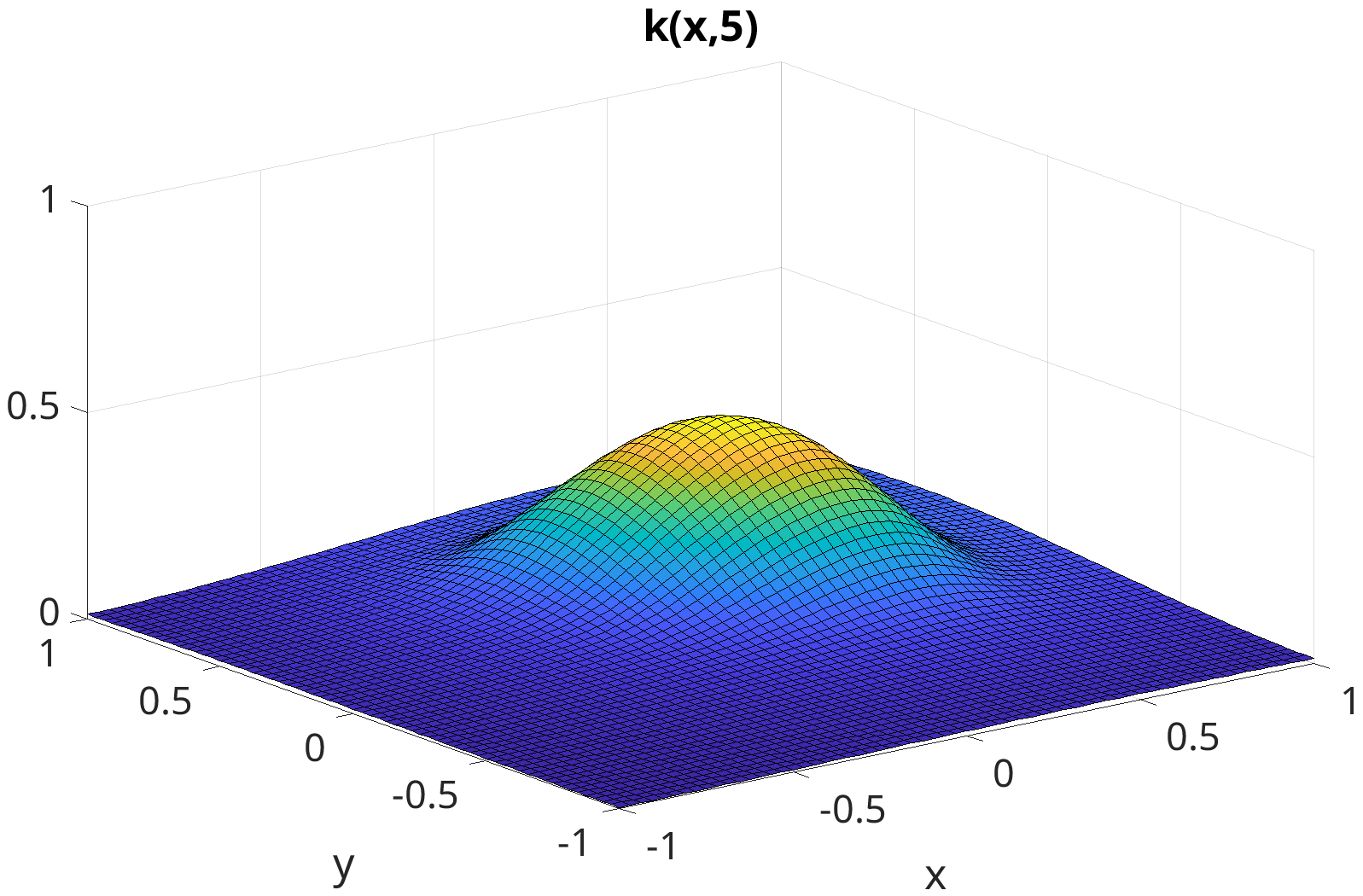}\\
		\includegraphics[width=4.8cm]{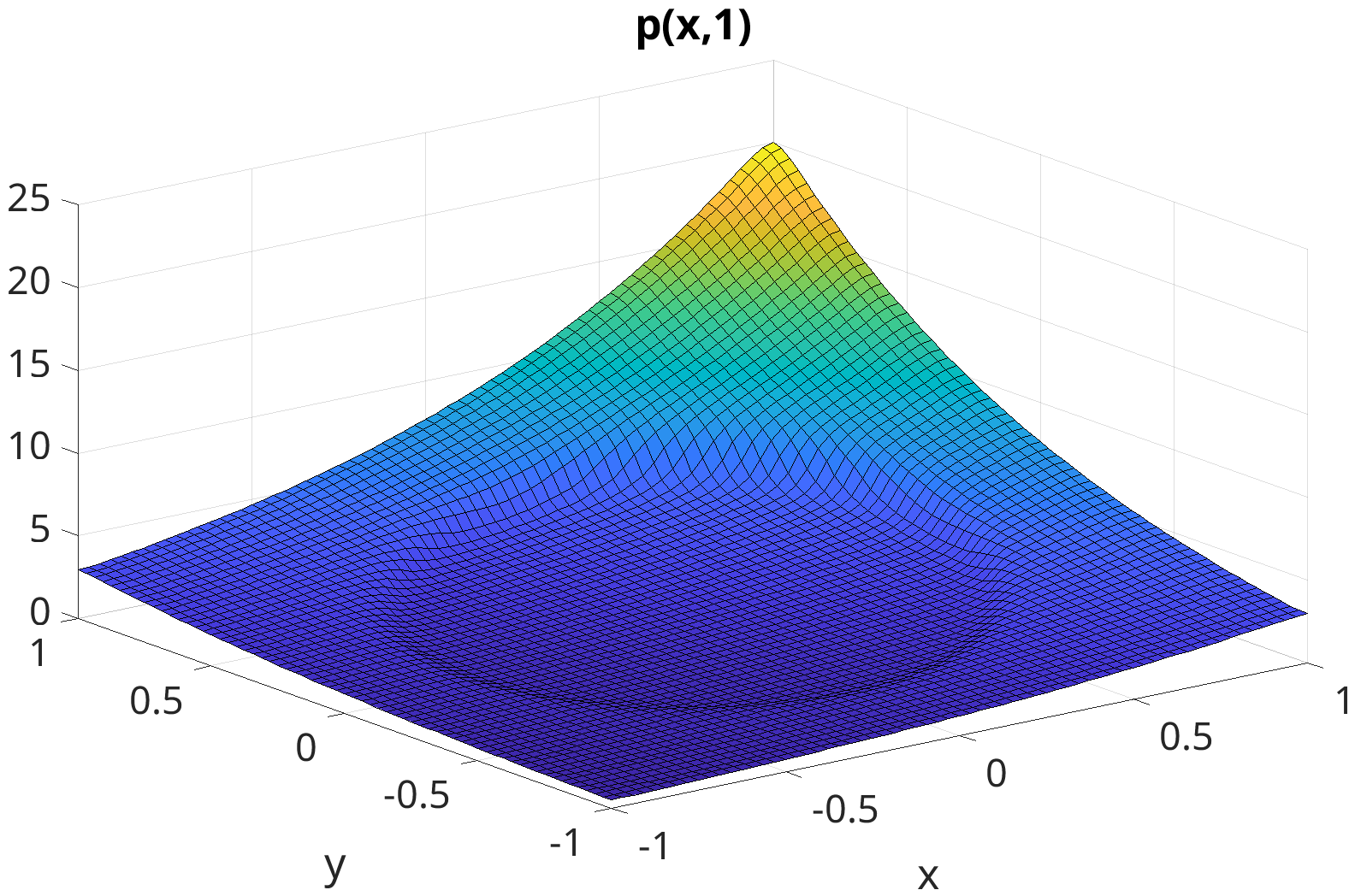}
		\includegraphics[width=4.8cm]{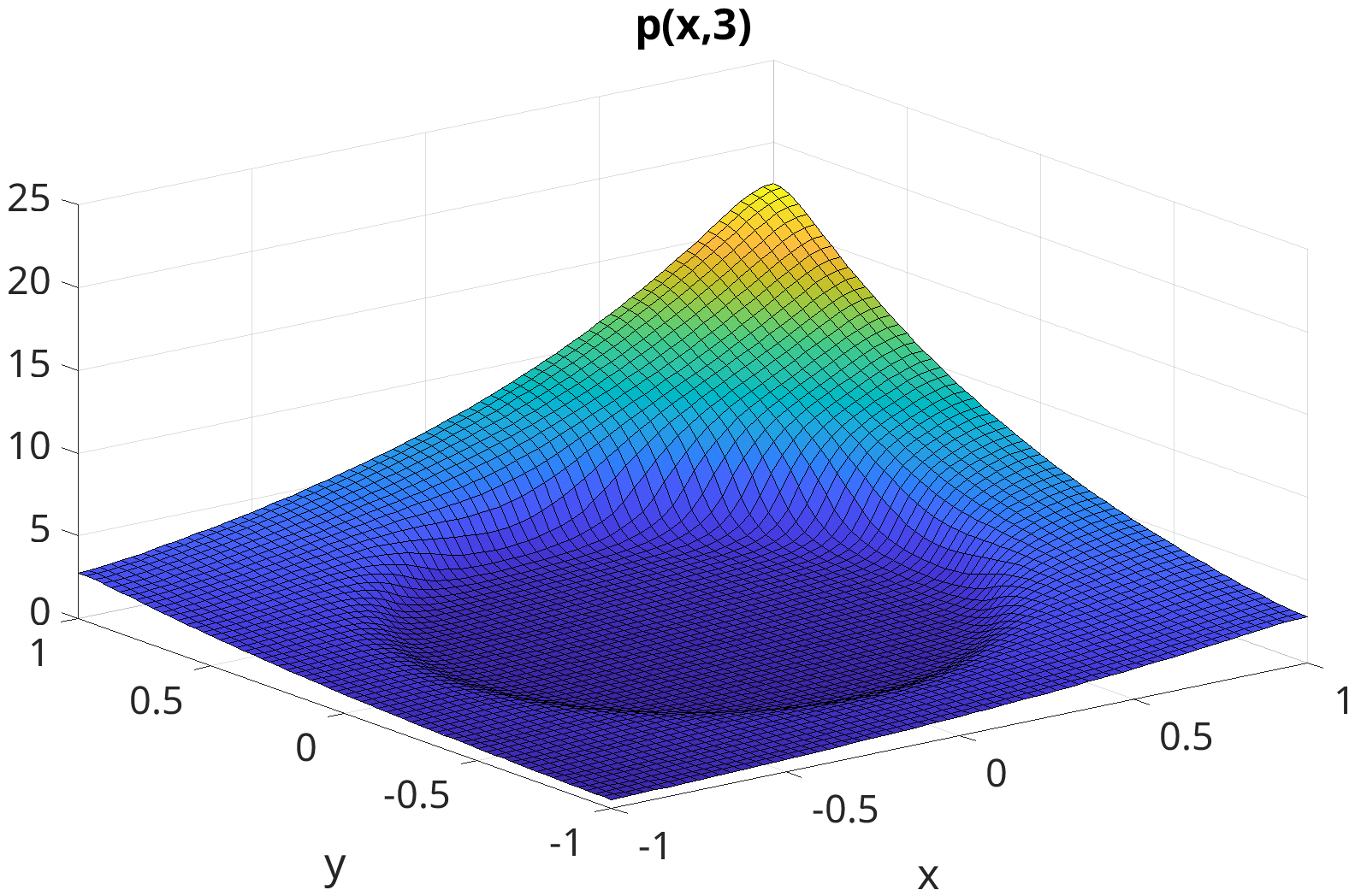}
		\includegraphics[width=4.8cm]{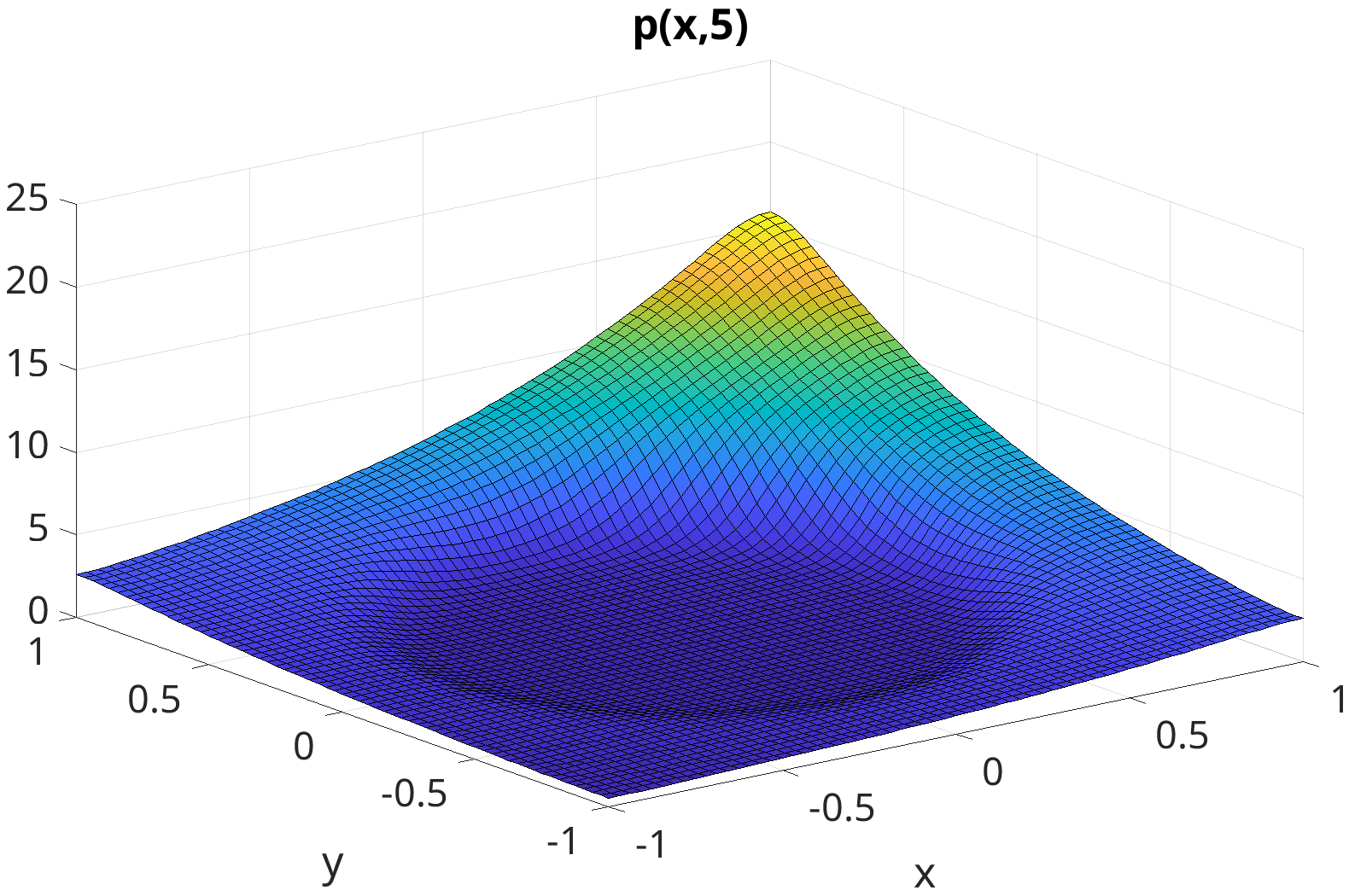}
		\caption{{\bf Case 2}, big control region. Space distributions of the optimal state variables $k(x, t)$ (first row) and $p(x, t)$ (second row) at three time instants.}
		\label{fig_reg_contr_big_k_p}
	\end{center}
\end{figure}

\begin{figure}[!htb]
	\begin{center}
		\includegraphics[width=4.8cm]{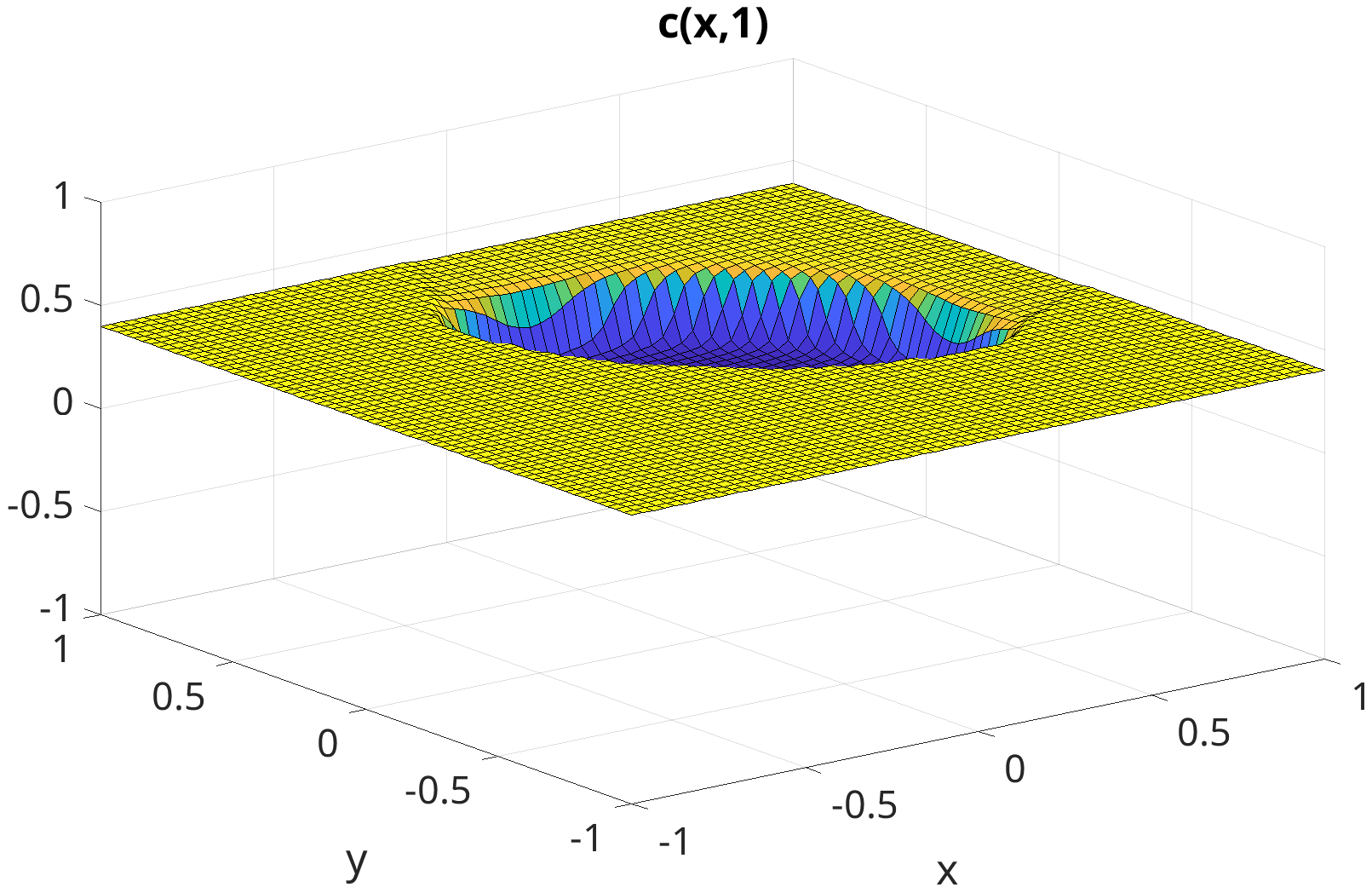}
		\includegraphics[width=4.8cm]{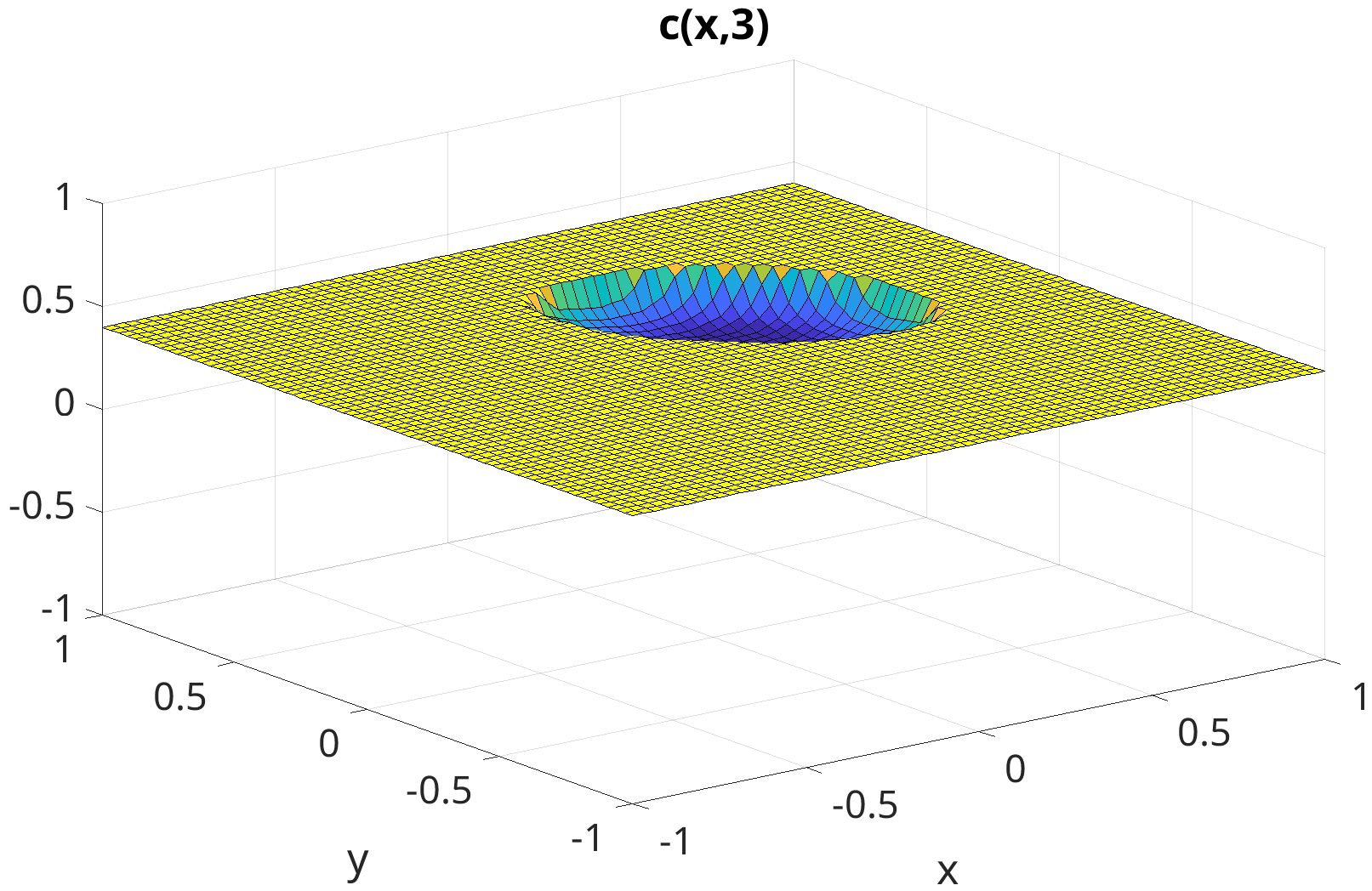}
		\includegraphics[width=4.8cm]{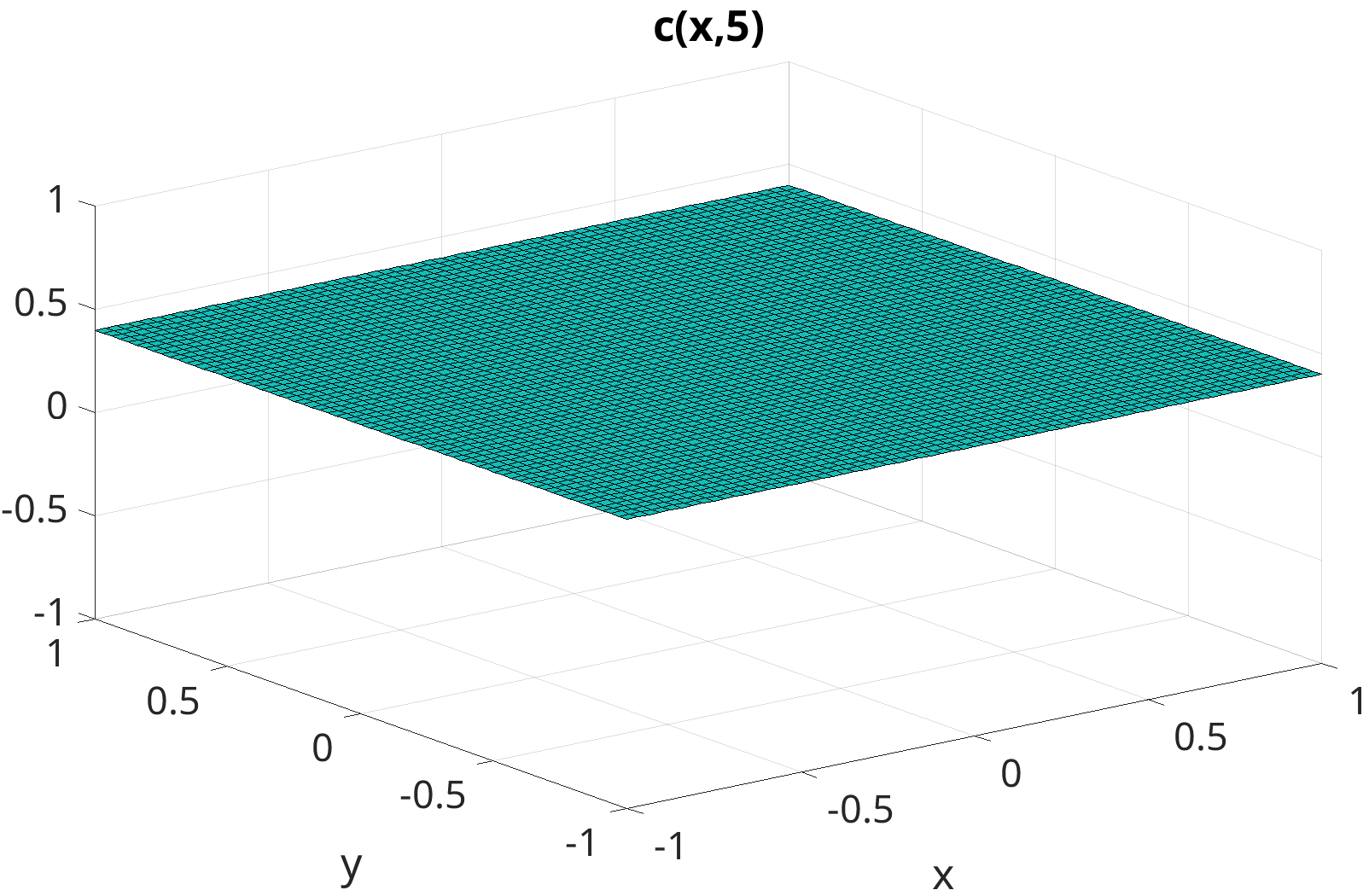}\\
		\includegraphics[width=4.8cm]{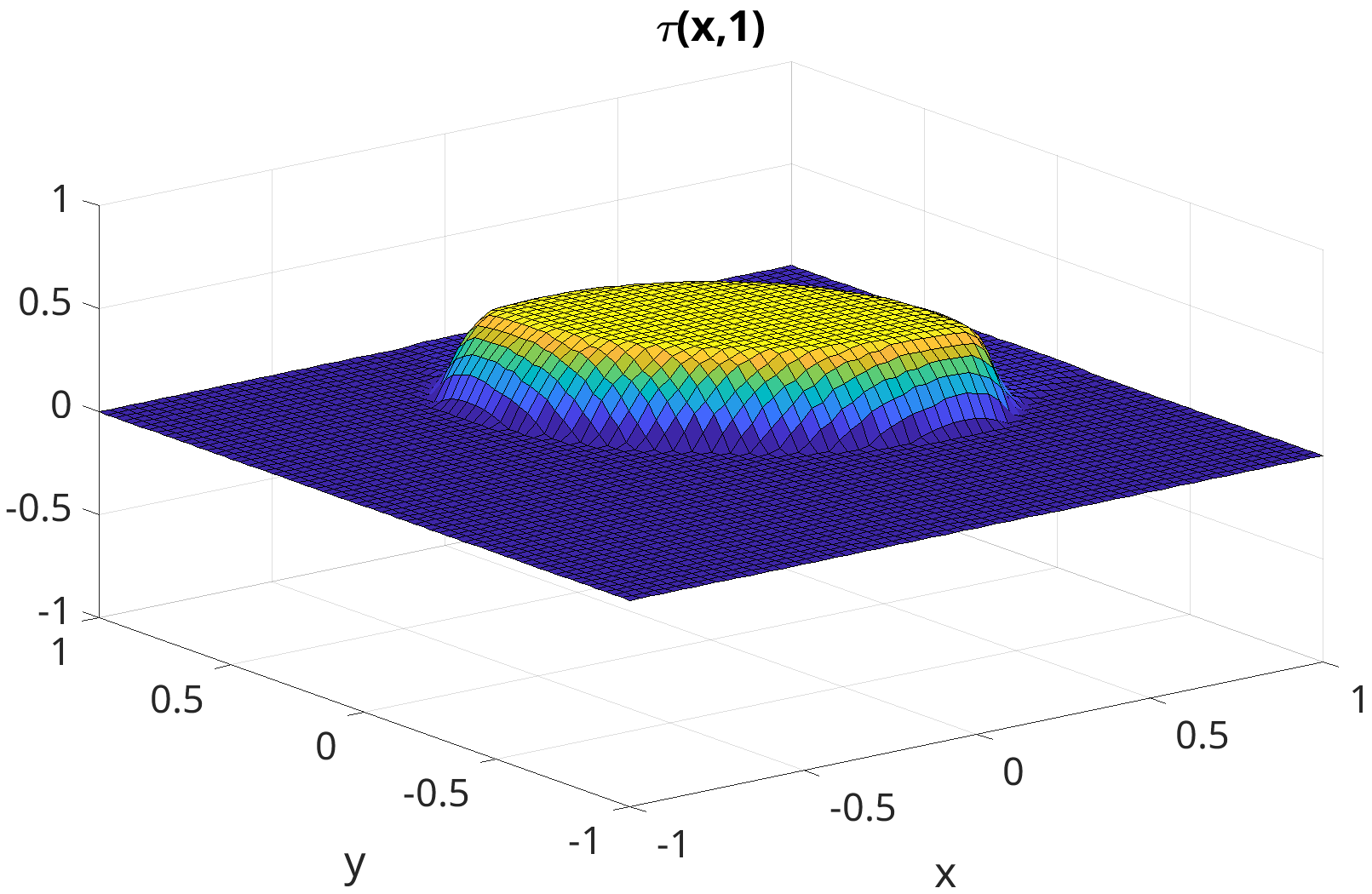}
		\includegraphics[width=4.8cm]{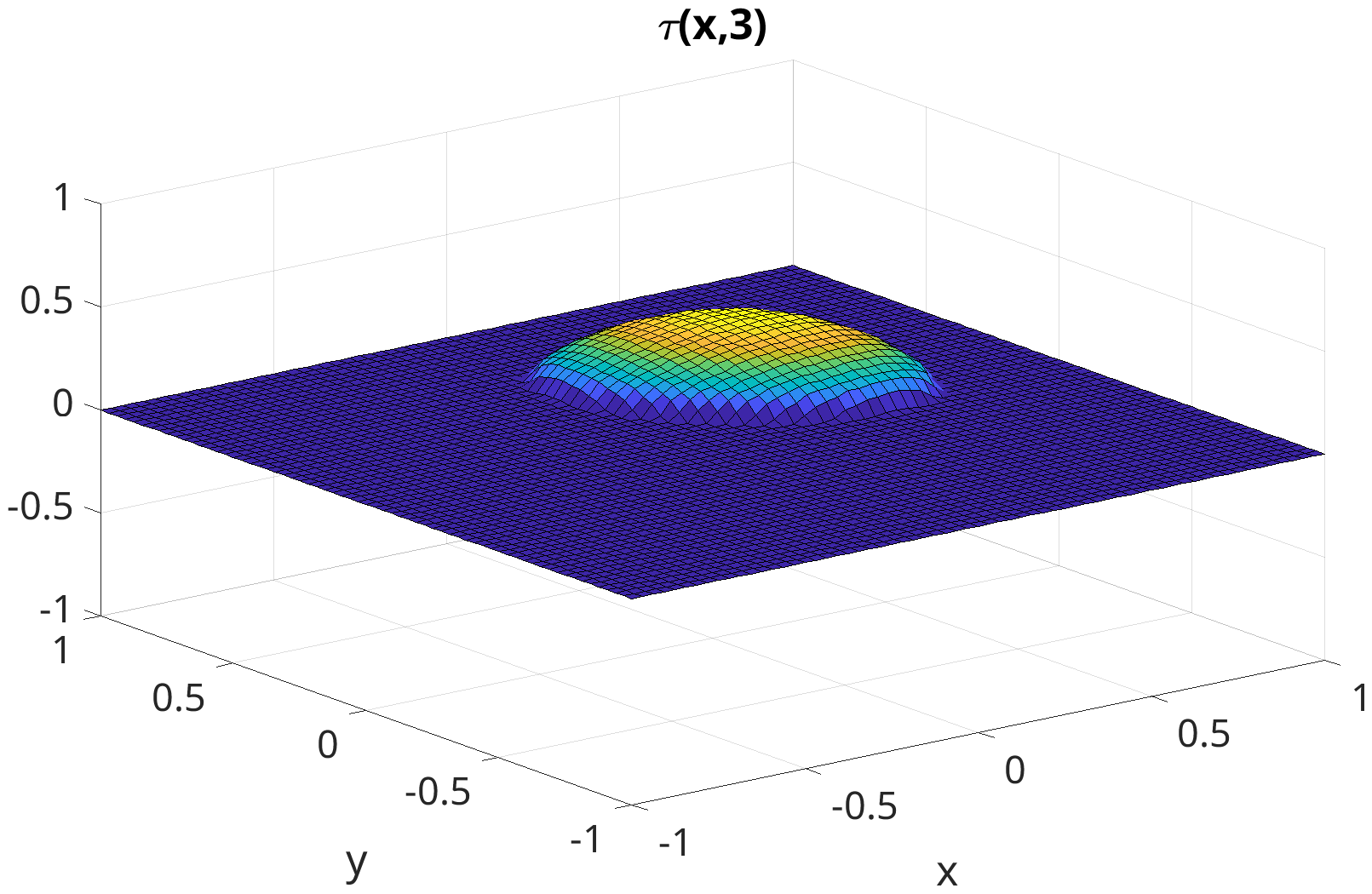}
		\includegraphics[width=4.8cm]{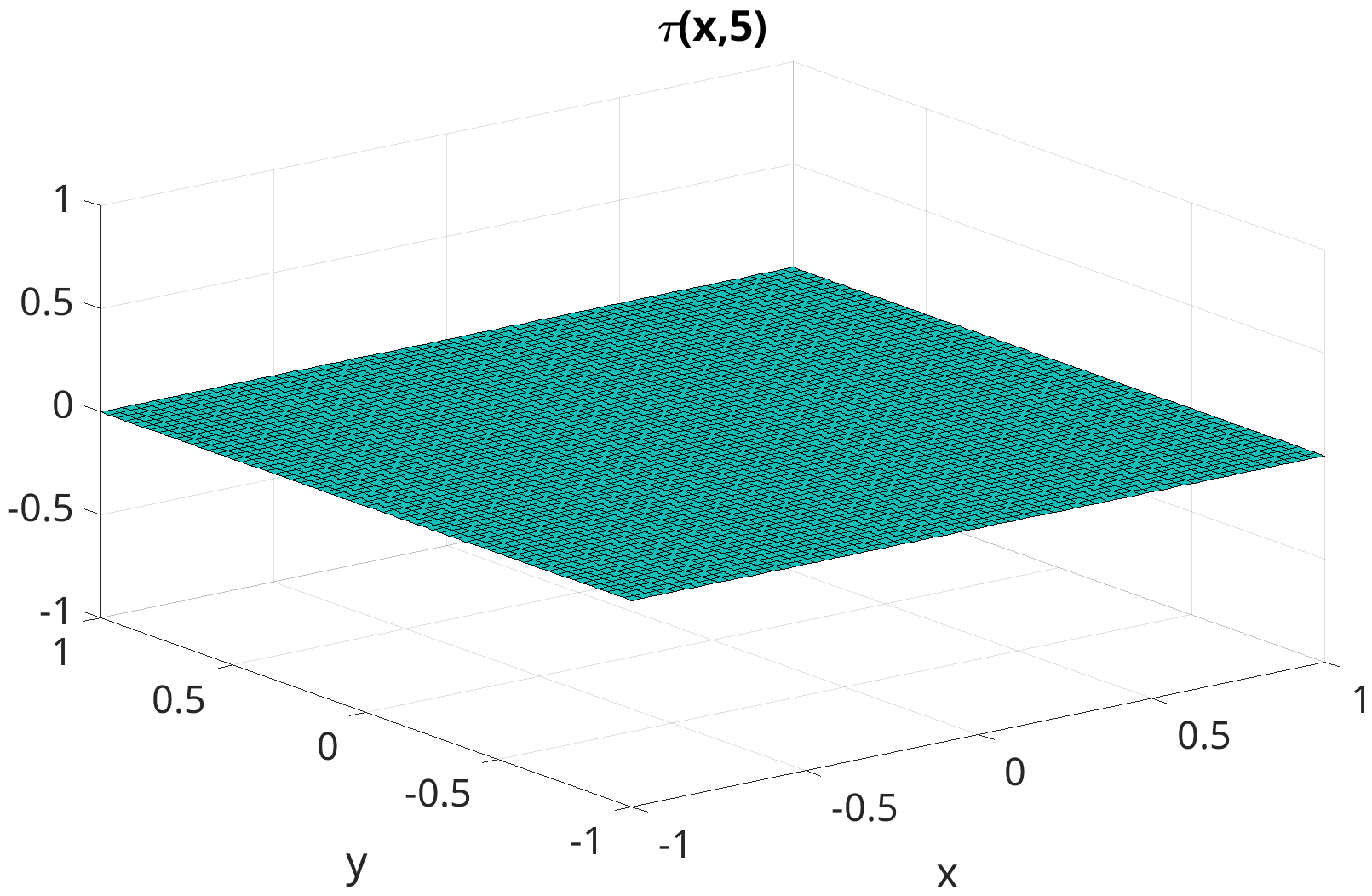}
		\includegraphics[width=4.8cm]{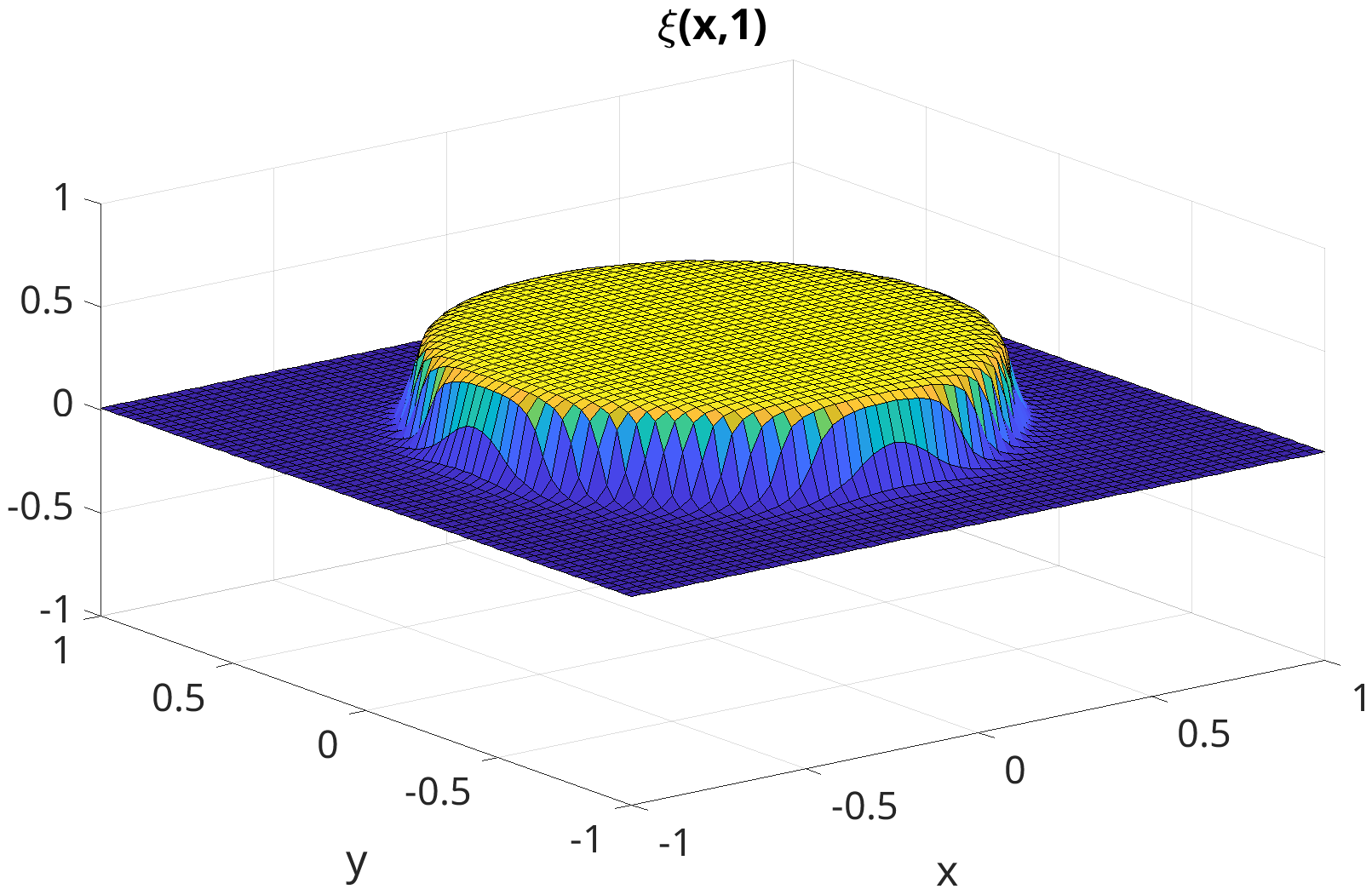}
		\includegraphics[width=4.8cm]{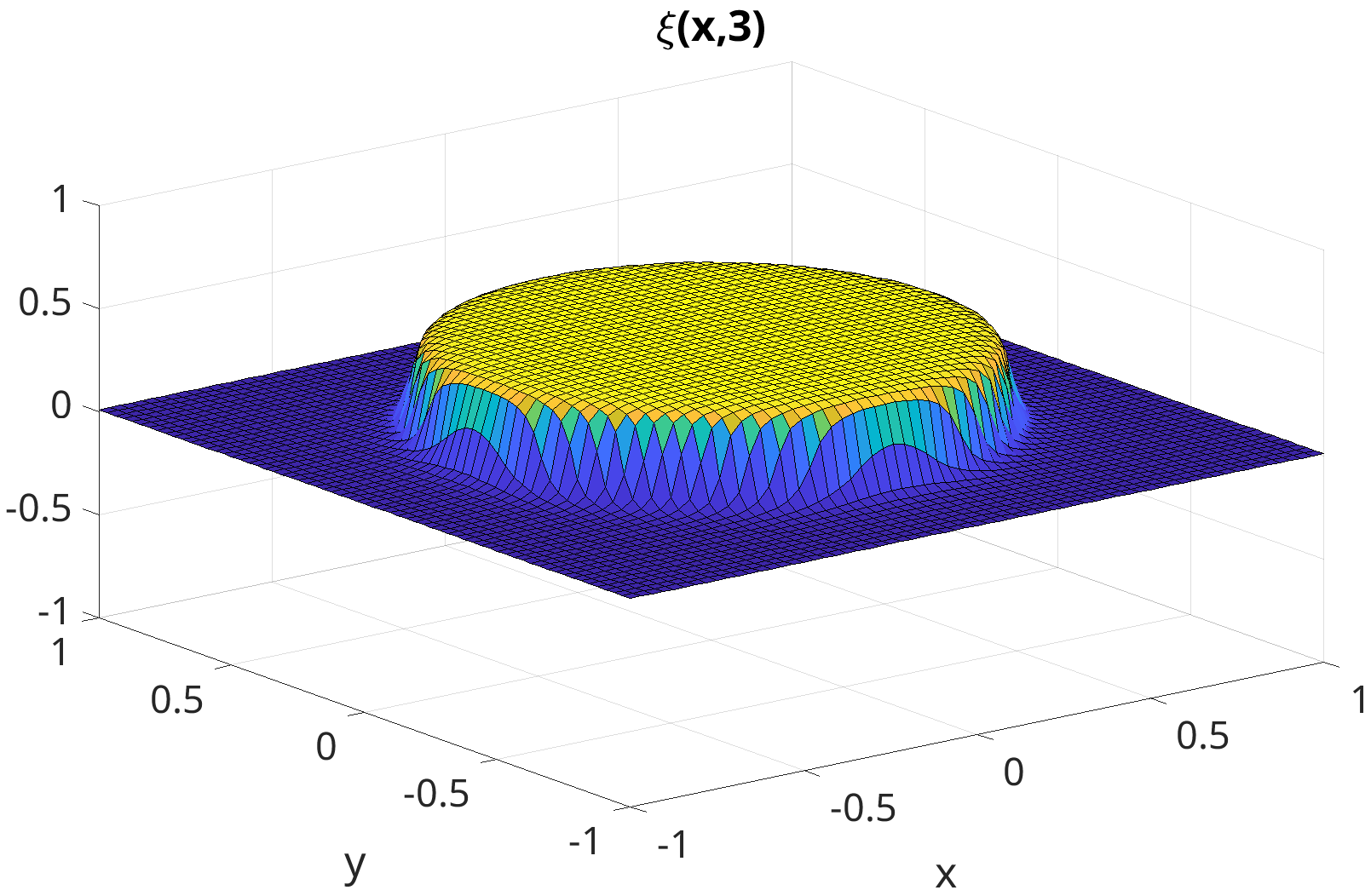}
		\includegraphics[width=4.8cm]{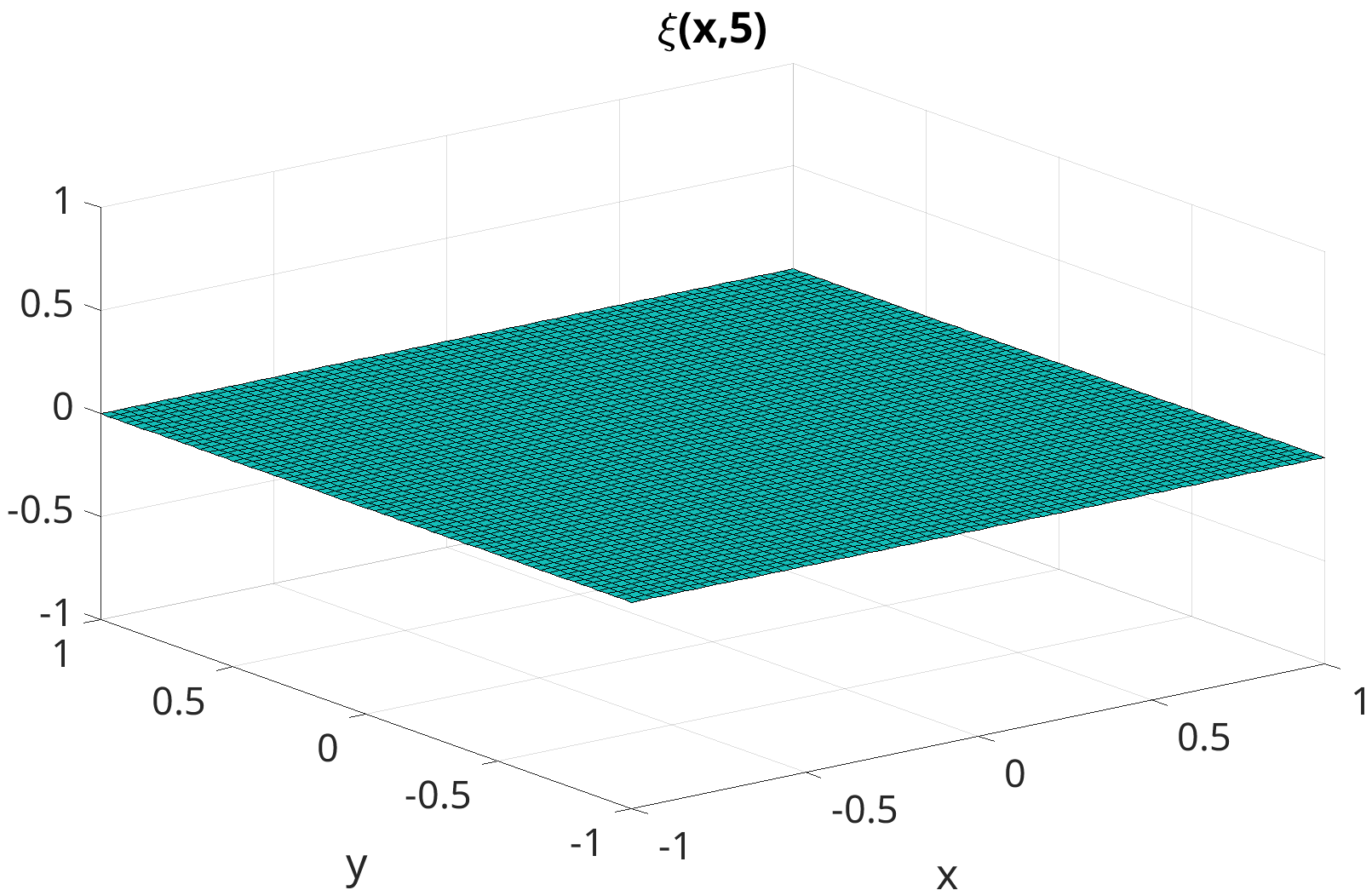}
		\caption{{\bf Case 2}, big control region. Space distributions of the optimal control variables $c(x, t)$ (first row), $\tau (x, t)$ (second row) and $\xi(x, t)$ (third row) at three time instants.}
		\label{fig_reg_contr_big_c_tau_xi}
	\end{center}
\end{figure}

\end{document}